\documentclass[]{interact}

\usepackage{epstopdf}
\usepackage[caption=false]{subfig}
\usepackage{algorithm}
\usepackage{algpseudocode}
\usepackage{graphicx}
\usepackage{wrapfig}

\usepackage{booktabs} 
\usepackage{array}
\usepackage{arydshln}
\usepackage[numbers,sort&compress]{natbib}
\bibpunct[, ]{[}{]}{,}{n}{,}{,}
\renewcommand\bibfont{\fontsize{10}{12}\selectfont}

\theoremstyle{plain}
\newtheorem{theorem}{Theorem}[section]
\newtheorem{lemma}[theorem]{Lemma}

\newtheorem{assum}{Assumption}

\theoremstyle{definition}

\theoremstyle{remark}

\numberwithin{equation}{section} 

\begin{document}


\title{Quasi-Newton Methods for Machine Learning: \\
Forget the Past, Just Sample}

\author{
\name{A.~S. Berahas\textsuperscript{a}\thanks{Email: albertberahas@gmail.com}, M. Jahani\textsuperscript{b}, P. Richt\'arik\textsuperscript{c} and M. Tak\'a\v c\textsuperscript{d}}
\affil{\textsuperscript{a}Department of Industrial and Operations Engineering,  University of Michigan; \textsuperscript{b}Department of Industrial and Systems Engineering, Lehigh University; \textsuperscript{c} Computer, Electrical and Mathematical Science and Engineering Division, KAUST; \textsuperscript{d} Mohamed bin Zayed University of Artificial Intelligence (MBZUAI)}
}

\maketitle

\begin{abstract}
We present two sampled quasi-Newton methods (sampled LBFGS and sampled LSR1) for solving empirical risk minimization problems that arise in machine learning. Contrary to the classical variants of these methods that sequentially build Hessian or inverse Hessian approximations as the optimization progresses, our proposed methods sample points randomly around the current iterate at every iteration to produce these approximations. As a result, the  approximations constructed make use of more reliable (recent and local) information, and do not depend on past iterate  information that could be significantly stale. Our proposed algorithms are efficient in terms of accessed data points (epochs) and have enough concurrency to take advantage of parallel/distributed computing environments. We provide convergence guarantees for our proposed methods. Numerical tests on a toy classification problem as well as on popular benchmarking binary classification and neural network training tasks reveal that the methods outperform their classical variants.
\end{abstract}


\section{Introduction}
\label{sec:intro}

In supervised machine learning, one seeks to minimize the empirical risk,
\begin{align}   \label{eq:opt_prob}
    \min_{w \in \mathbb{R}^d} F(w) := \frac{1}{n} \sum_{i=1}^n f(w;x^i,y^i) = \frac{1}{n} \sum_{i=1}^n f_i(w)
\end{align}
where $f: \mathbb{R}^d \rightarrow \mathbb{R}$ is the composition of a prediction function (parametrized by $w$) and a loss function, and $(x^i,y^i)$, for $i=1,\dots,n$, denote the training examples (samples). Difficulties arise in minimizing the function $F$ for three main reasons: $(1)$ the number of samples $n$ is large; $(2)$ the number of variables $d$ is large; and, $(3)$ the objective function is nonconvex.

In the last decades, much effort has been devoted to the development of stochastic first-order methods that have a low per-iteration cost, enjoy optimal complexity, are easy to implement, and that have proven to be effective for many machine learning applications. At present, the preferred method for large-scale applications is the stochastic gradient (SG) method \cite{robbins1951stochastic,bottou2004large}, and its variance-reduced \cite{johnson2013accelerating,schmidt2017minimizing,defazio2014saga,nguyen2017sarah} and adaptive variants \cite{duchi2011adaptive,kingma2014adam}. However, these methods have several issues: (1) they are highly sensitive to the choice of hyper-parameters (e.g., steplength and batch size) and tuning can be  cumbersome; (2) they suffer from ill-conditioning; and, (3) they often offer limited opportunities for parallelism; see \cite{berahas2017investigation,xu2017second,kylasa2018gpu,Roosta-Khorasani2018,bottou2018optimization}.

In order to alleviate these issues, stochastic Newton \cite{byrd2011use,martens2010deep,bollapragada2016exact,Roosta-Khorasani2018,xu2017newton} and stochastic quasi-Newton \cite{schraudolph2007stochastic,byrd2016stochastic,curtis2016self,mokhtari2015global,gower2016stochastic,berahas2016multi,keskar2016adaqn,berahas2017robust,jahani2019scaling,jahani2020sonia} methods have been proposed. These methods attempt to combine the speed of Newton's method and the scalability of first-order methods by incorporating curvature information in a judicious manner, and have proven to work well for several machine learning tasks \cite{berahas2017investigation,xu2017second}.

With the advances in distributed and GPU computing, it is now possible to go beyond stochastic Newton and quasi-Newton methods and use large batches, or even the full dataset, to compute function, gradient and Hessian vector products in order to train machine learning models. In the large batch regime, one can take advantage of parallel and distributed computing and fully utilize the capabilities of GPUs. However, researchers have observed that \textit{well-tuned} first-order methods (e.g., ADAM) are far more effective than full batch methods (e.g., LBFGS) for large-scale  applications \cite{hardt2016train,keskar2016large}.

Nevertheless, in this paper we focus on (full) batch methods that incorporate local second-order (curvature) information of the objective function. These methods mitigate the effects of ill-conditioning, avoid or diminish the need for hyper-parameter tuning, have enough concurrency to take advantage of parallel computing, and, due to requiring fewer iterations enjoy low communication costs in distributed computing environments. Specifically, we focus on quasi-Newton methods \cite{nocedal_book}; methods that construct curvature information using first-order (gradient) information. We propose two variants of classical quasi-Newton methods that sample a small number of random points at every iteration to build (inverse) Hessian approximations.

\begin{figure}
  \centering
    \includegraphics[width=0.9\textwidth]{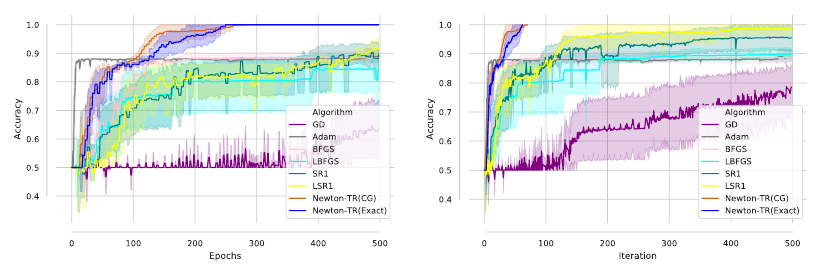}
  \caption{\small Performance of GD, ADAM, BFGS, LBFGS, SR1, LSR1, Newton-TR(CG, Exact) on a toy classification problem.}
  \label{MotivationFigure}
\end{figure}

We are motivated by the results presented in Figure \ref{MotivationFigure} that illustrate the performance (for 10 different starting points) of several stochastic and deterministic, first- and second-order methods on a toy neural network classification task, given budget; see Section \ref{sec:num_res} for details. As is clear from the results, first-order methods converge very slowly, and sometimes even fail to achieve 100\% accuracy. Similarly, classical quasi-Newton methods are also slow or stagnate. On the other hand, methods that use the true Hessian are able to converge in very few iterations from all starting points. This seems to suggest that for some neural network training tasks second-order information is important, and that the curvature information captured by classical quasi-Newton methods may not be adequate or useful.

The key idea of our proposed methods is to leverage the fact that quasi-Newton methods can incorporate second-order information using only gradient information at a reasonable cost, but at the same time to enhance the (inverse) Hessian approximations by using more reliable (recent and local) information. The fundamental component of our methods, and what differentiates them from the classical variants, is the manner in which the curvature pairs are constructed. To this end, we propose to \textit{forget} past curvature information and \textit{sample} new curvature pairs at every iteration.

\paragraph*{Contributions}Our contributions can be summarized as follows: 
\begin{itemize}
    \item We propose two novel quasi-Newton methods that use sampling to construct Hessian approximations. We analyze the convergence properties of both methods, and show that their theoretical guarantees match those of their classical limited memory counterparts.
    \item We derive expressions for the probability of accepting curvature pairs within our quasi-Newton frameworks. Moreover, we illustrate empirically the probability of accepting curvature pairs on a toy classification problem.
    \item We discuss the implementation costs of the sampled quasi-Newton methods and compare them to the classical variants, and illustrate the scaling properties of the methods compared to the SG method on distributed computing platforms on real large-scale network architectures.
    \item We illustrate the practical performance of the methods on a toy classification problem as well as on popular benchmarking binary classification and neural network training tasks, and show their advantages over the classical variants. We posit that this is the case since the (inverse) Hessian approximations constructed by our proposed methods capture better (more informative) curvature information. Moreover, the proposed methods are easily parallelizable and efficient in terms of iteration, epochs and communication.
\end{itemize}

The paper is organized as follows. We conclude this section with a literature review of quasi-Newton methods. We describe the classical (L)BFGS and (L)SR1 methods in Section \ref{sec:qn}, and in Section \ref{sec:samp_qn} we present our proposed sampled quasi-Newton variants. In Section \ref{sec:dist_comp}, we discuss the computational cost of the proposed methods and show their scaling properties. We show the theoretical properties of our proposed methods in Section \ref{sec:conv_analysis}. In Section \ref{sec:prob} we provide a theoretical and numerical analysis on the probabilities of accepting the sampled points within our proposed quasi-Newton frameworks. Numerical results on neural network training tasks are reported in Section \ref{sec:num_res}. Finally, in Section \ref{sec:fin_rem} we provide some final remarks and discuss several avenues for future work.

\paragraph*{Literature Review}
Quasi-Newton methods, such as BFGS \cite{broyden1967quasi, fletcher1970new,goldfarb1970family,shanno1970conditioning} and SR1 \cite{byrd1996analysis,conn1991convergence,khalfan1993theoretical} and their limited-memory variants LBFGS \cite{nocedal1980updating,liu1989limited} and LSR1 \cite{brust2017solving,lu1996study}, respectively, have been studied extensively in the deterministic nonlinear optimization literature. These methods incorporate curvature (second-order) information using only gradient (first-order) information, have good theoretical guarantees, and have proven to be effective in practice.

In the context of deep neural networks, both full batch and stochastic quasi-Newton methods seem to perform worse than (stochastic) first-order methods. Nevertheless, several stochastic quasi-Newton methods have been proposed; see e.g., \cite{schraudolph2007stochastic,byrd2016stochastic,berahas2017robust}. What distinguishes these methods from one another is the way in which curvature pairs are constructed. Our methods borrow some of the ideas proposed in  \cite{byrd2016stochastic,gower2016stochastic,liu2018acceleration}. Specifically, we use Hessian vector products in lieu of gradient displacements.

Possibly the closest works to ours are Block BFGS \cite{gao2018block} and its stochastic variant \cite{gower2016stochastic}. These methods construct multiple curvature pairs to update the quasi-Newton matrices. However, there are several key features that are different from our approach; in these works $(1)$ the Hessian approximation is not updated at every iteration, and $(2)$ they enforce that multiple secant equations hold simultaneously.

\section{Quasi-Newton Methods}
\label{sec:qn}

In this section, we review two classical quasi-Newton methods (BFGS and SR1) and their limited memory variants (LBFGS and LSR1). This will set the stage for our proposed sampled quasi-Newton methods.

\subsection{BFGS and LBFGS}
Let us begin by considering the BFGS method and then consider its limited memory version. At the $k$th iteration, the BFGS method computes a new iterate by the formula
\begin{align}   \label{eq:bfgs}
    w_{k+1} = w_k - \alpha_k H_k \nabla F(w_k),
\end{align}
where $\alpha_k$ is the step length, $\nabla F(w_k)$ is the gradient of \eqref{eq:opt_prob} and $H_{k}$ is the inverse BFGS
Hessian approximation that is updated at every iteration by means
of the formula
\begin{gather*}   
H_{k+1}=V_{k}^{T}H_{k}V_{k}+\rho_{k}s_{k}s_{k}^{T},\\ 
\rho_{k}=\tfrac{1}{y_{k}^{T}s_{k}}, \quad V_{k}=I-\rho_{k}y_{k}s_{k}^{T} ,
\end{gather*}
where  the curvature pairs $(s_k,y_k)$ are defined as
\begin{align}   \label{eq:def_sy}
    s_{k} = w_{k} - w_{k-1}, \quad y_{k} = \nabla F(w_{k}) - \nabla F(w_{k-1}).
\end{align}
As is clear, the curvature pairs \eqref{eq:def_sy} are constructed sequentially (at every iteration), and as such the inverse Hessian approximation at the $k$th iteration $H_k$ depends on iterate (and gradient) information from past iterations.

The inverse BFGS Hessian approximations are constructed to satisfy two conditions:
\begin{align*}
    H_{k+1}y_k = s_k, \quad  \text{and} \quad s_k^Ty_k >0,
\end{align*}
the secant and curvature conditions, respectively, as well as symmetry. Consequently, as a result, as long as the initial inverse Hessian approximation is positive definite, then all subsequent inverse BFGS Hessian approximations are also positive definite. Note, the new (inverse) Hessian approximation $H_{k+1}$ differs from the old approximation $H_k$ by a rank-2 matrix. 

In the limited memory version, the matrix $H_k$ is defined at each iteration as the result of applying $m$ BFGS updates to a multiple of the identity matrix using the set of $m$ most recent curvature pairs $\{s_i, y_i\}$ kept in storage. As a result, one need not store the dense inverse Hessian approximation, rather one can store two $m \times d$ matrices and compute the matrix-vector product in \eqref{eq:bfgs} via the two-loop recursion \cite{nocedal_book}. After the step has been computed, the oldest pair $(s_j, y_j)$ is discarded and the new curvature pair is stored.

\subsection{SR1 and LSR1}

Contrary to the BFGS updating formula, and as suggested by the name, the symmetric-rank-1 (SR1) updating formula allows one to satisfy the secant equation and maintain symmetry with a simpler rank-1 update. However, unlike BFGS, the SR1 update does not guarantee that the updated matrix maintains positive definiteness. As such, the SR1 method is usually implemented with a trust region; we introduce it in this way below.

At the $k$th iteration, the SR1 method computes a new iterate by the formula
\begin{align}   \label{eq:sr1}
    w_{k+1} = w_k + p_k,
\end{align}
where $p_k$ is the minimizer of the following subproblem
\begin{align}
    &{\min_{p}}  \; m_k(p) = F(w_k) + \nabla F(w_k)^Tp + \tfrac{1}{2} p^T B_k p , \label{eq:tr_obj}\\
    & \quad \text{s.t.} \qquad \| p \| \leq \Delta_k, \nonumber
\end{align}
$\Delta_k$ is the trust region and $B_k$ is the SR1 Hessian approximation computed as
\begin{align}   \label{eq:sr1_hess}
    B_{k+1} = B_k +  \tfrac{(y_k - B_ks_k)(y_k - B_ks_k)^T}{(y_k - B_ks_k)^Ts_k}. 
\end{align}
Similar to LBFGS, in the limited memory version of SR1 the matrix $B_k$ is defined at each iteration as the result of applying $m$ SR1 updates to a multiple of the identity matrix, using a set of $m$ correction pairs $\{s_i, y_i\}$ kept in storage.

\section{Sampled Quasi-Newton Methods}
\label{sec:samp_qn}

In this section, we describe our two proposed sampled quasi-Newton methods; S-LBFGS and S-LSR1. The main idea of these methods, and what differentiates them from the classical variants, is the way in which curvature pairs are constructed. At every iteration, a small number ($m$) of points are sampled around the current iterate and used to construct a new set of curvature pairs. In other words, contrary to the sequential nature of classical quasi-Newton methods, our proposed methods \textit{forget} all past curvature pairs and construct new curvature pairs from scratch via \textit{sampling}. 

\begin{figure}
	\centering
	\includegraphics[width=1\textwidth]{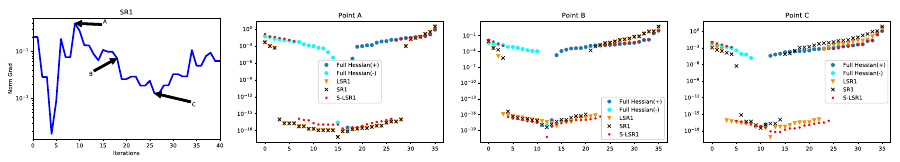}
	\caption{\small Comparison of the eigenvalues of (L)SR1 and S-LSR1 (@ A, B, C) for a toy classification problem (\texttt{small network}). The $x$-axis in the left-most plot (in Figures \ref{SmalllvdNet} and \ref{MediumvdNet}) depicts the iterations; the $x$-axis in the plots in columns 2-4 (Figures \ref{SmalllvdNet} and \ref{MediumvdNet}) depicts the indices of eigenvalues (sorted in ascending order) of the Hessian and Hessian approximations.
	\label{SmalllvdNet}}
\end{figure}

\begin{figure}
	\centering
	\includegraphics[width=1\textwidth]{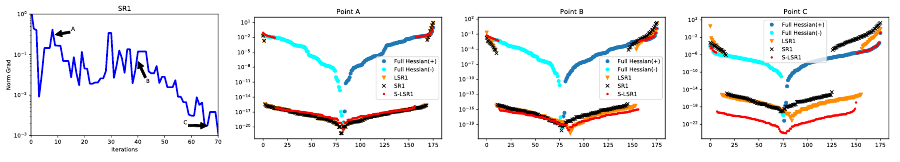}
	\caption{\small Comparison of the eigenvalues of (L)SR1 and S-LSR1 (@ A, B, C) for a toy classification problem (\texttt{medium network}).
	\label{MediumvdNet}}
\end{figure}

Our motivation stems from the following observation: by constructing Hessian approximations via sampling, one is able to better capture curvature information of the objective function. In Figures \ref{SmalllvdNet} and \ref{MediumvdNet}, we show the spectrum of the true Hessian, and compare it to the spectra of different SR1 Hessian approximations at several points for two toy classification problems.  As is clear from the results, the eigenvalues of the S-LSR1 Hessian approximations better match the eigenvalues of the true Hessian compared to the eigenvalues of the SR1 and LSR1 Hessian approximations. This is not surprising since S-LSR1 uses newly sampled local information, and unlike the classical variants does not rely on past information that could be significantly stale. Similar results were obtained for other problems; see Appendix \ref{sec:eigs_app} for details.

This, of course, does not come for free. The classical variants construct curvature pairs as the optimization progresses at no additional cost, whereas the sampled quasi-Newton methods require the construction of $m$ new curvature pairs at every iteration. We discuss implementation issues and the computational cost of the sampled quasi-Newton methods in Sections \ref{sec:samp_SYpairs} and \ref{sec:dist_comp}.

We now discuss the way in which curvature pairs are constructed, and then formally introduce our proposed sampled quasi-Newton methods.

\subsection{Sampling Curvature Pairs}
\label{sec:samp_SYpairs}

As mentioned above, the key component of our proposed algorithms is the way in which curvature pairs are constructed. A pseudo-code of our proposed sampling strategy and construction of the curvature pairs in given in Algorithm \ref{alg:calSY}. Let $S \in \mathbb{R}^{d \times m}$ and $Y \in \mathbb{R}^{d \times m}$ denote the matrices of all curvature pairs constructed during the $k$th iteration.

\begin{algorithm}
{ 
\small
\caption{Compute new $(S,Y)$ curvature pairs}
  \label{alg:calSY}
 {\bf Input:} $w$ (iterate), $m$ (memory), $r$ (sampling radius), $S = [\;]$, $Y =[\;]$ (curvature pair containers). 

  \begin{algorithmic}[1]
  \State Compute $\nabla F(w)$
  \For {$i=1,2,...,m$}
  		\State Sample a random direction $\sigma_i$
        \State Construct $\bar{w} = w + r\sigma_i $
        \State Set $s = w - \bar{w}$ and 
        
        $\qquad y = \begin{cases}
    \nabla F(w) - \nabla F(\bar{w}), \ \ &\text{Option I} \\
    \nabla^2 F(w)s,\ \ &\text{Option II}
\end{cases}$
        \State Set $S = [S \; s]$ and $Y = [Y \; y]$
    \EndFor
  \end{algorithmic}
  {\bf Output:} $S, Y$
  }
\end{algorithm}

Both S-LBFGS and S-LSR1 use the subroutine described in Algorithm \ref{alg:calSY}. At every iteration, given the current iterate and gradient, $m$ curvature pairs are constructed. The subroutine first samples points around the current iterate along a random directions $\sigma_i$ and sets the iterate displacement curvature pair ($s$), and then creates the gradient difference curvature pair ($y$) via gradient differences (Option I) or Hessian vector products (Option II). Note that the random directions $\sigma_i$ can be arbitrary; in the latter part of the paper (Sections \ref{sec:num_res} and \ref{sec:prob}), we make an explicit choise on the directions. 

Our theory holds for both options; however, in our numerical experiments we present results with Option II only for the following reasons. Option I requires $m$ gradient evaluations ($m$ epochs), and thus requires accessing the data $m$ times. On the other hand, Option II only requires a single Hessian matrix product which can be computed very efficiently on a GPU, as the $y$ curvature pairs can be constructed simultaneously, i.e., $Y = \nabla^2 F(w) S$, and thus only requires accessing the data once. Moreover, Option I requires choosing the sampling radius $r$, whereas Option II does not since it is scale invariant.

\begin{figure}[H]
	\centering
	\includegraphics[width=0.92\textwidth]{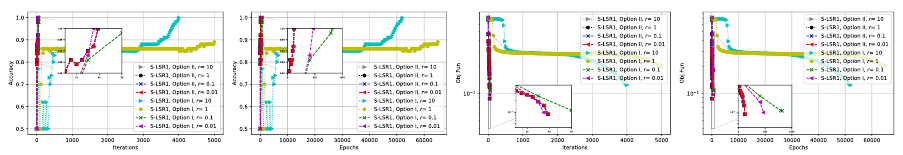}
	\caption{\small Comparison of algorithms with Option I and Option II on toy classification problem (\texttt{small network}).
	\label{opt1_opt2_smallNet}}
\end{figure}

\begin{figure}[H]
	\centering
	\includegraphics[width=0.92\textwidth]{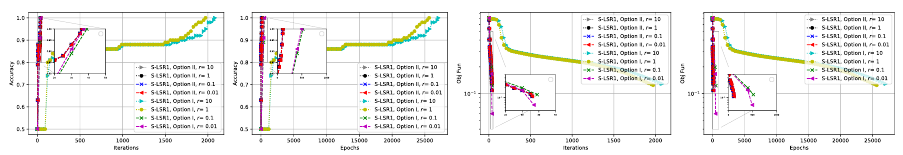}
	\caption{\small Comparison of algorithms with Option I and Option II on toy classification problem (\texttt{medium network}).
	\label{opt1_opt2_mediumNet}}
\end{figure}

Before we proceed with our presentation of the S-LBFGS and S-LSR1 methods, we empirically compare the performance of a methods that uses Option I and Option II. As is clear from Figures \ref{opt1_opt2_smallNet} and \ref{opt1_opt2_mediumNet}, the performance of the method that uses Option I is highly dependant on the choice of the sampling radius ($r$). If this parameter is not chosen appropriately, the performance of the method can be slow. This is not the case when Option II is utilized, and one can attribute this to the fact that Option II is scale invariant. Moreover, the benefits of using Option II can clearly be observed in the plots in terms of epochs. Again, this is due to the fact each iteration using Option I requires accessing the data at  $m$ times to construct the curvature pairs, whereas Option II required only a single access of the data to construct the curvature pairs.

\subsection{Sampled LBFGS (S-LBFGS)}

At the $k$th iteration, the S-LBFGS method computes a new iterate via \eqref{eq:bfgs}, where the inverse Hessian approximation is constructed using the curvature pairs sampled by Algorithm \ref{alg:calSY}. The S-LBFGS method is outlined in Algorithm \ref{alg:slbfgs}. 

 \begin{algorithm}
 {\small 
\caption{Sampled LBFGS (S-LBFGS)}
  \label{alg:slbfgs}
 {\bf Input:} $w_{0}$ (initial iterate), $m$ (memory), $r$ (sampling radius). 
 
  \begin{algorithmic}[1]
  \For {$k=0,1,2,...$}
    \State Compute new $(S_k,Y_k)$ pairs via Algorithm \ref{alg:calSY}
    \State Compute the search direction $p_{k}=-H_{k}\nabla F(w_k)$
    \State Choose the steplength $\alpha_{k} >0$
    \State Set  $w_{k+1}=w_k+\alpha_{k}p_{k}$ 
\EndFor
  \end{algorithmic}
  }
\end{algorithm}

Algorithm \ref{alg:slbfgs} is almost identical to the classical (L)BFGS algorithm \cite{nocedal_book}; however, it has two key differentiating elements: (1) the way in which  curvature pairs are created; and, (2) the location in the algorithm where the curvature pairs are constructed. Both elements can be interpreted as features of S-LBFGS. First, using a similar argument as that for the S-LSR1 method (Figure \ref{SmalllvdNet}), the inverse Hessian approximations constructed by this method better capture local curvature information of the objective function. Moreover, notice that the first set of curvature pairs is constructed before a single step is taken by the method (Line 2). This allows the method to take a quasi-Newton-type (well-scaled) step from the first iteration which is not the case for classical BFGS methods that usually take a gradient-type step in the first iteration, and in which imposing the correct scale can be an issue. This, possibly, is a more important implication of the method, as the first step taken by quasi-Newton methods can be of paramount importance. 

In order to fully specify the S-LBFGS method, we need to describe how the steplength is selected (Algorithm \ref{alg:slbfgs}, Step 4). We consider two variants of the method: $(1)$ constant steplength selection, and $(2)$ adaptive steplength selection. Our theory (Section \ref{sec:conv_analysis}, Theorems \ref{thm:const_strongly} and \ref{thm_non_MB}), explicitly defines the manner in which the steplength should be selected in order to ensure convergence. Of course, in practice, one can (potentially) use a larger steplength, and as such in this approach the steplength ($\alpha_k = \alpha$) is a tuneable parameter. We also consider an adaptive Armijo backtracking mechanism for selecting the steplength at every iteration. Given the current iterate $w_k$, the steplength is chosen to satisfy the following sufficient decrease condition
\begin{align}   \label{eq:armijo}
    F(w_{k} + \alpha_k p_k) \leq F(w_k) - c_1 \alpha_k \nabla F(w_k)^T H_k \nabla F(w_k)
\end{align}
where $c_1 \in (0,1)$. The mechanism works as follows. Given an initial steplength (say $\alpha_k = 1$), the function is evaluated at the trial point $w_{k} + \alpha_k p_k$ and condition \eqref{eq:armijo} is checked. If the trial point satisfies \eqref{eq:armijo}, then the step is accepted. If the trial point does not satisfy \eqref{eq:armijo}, the steplength is reduced (e.g., $\alpha_k = \tau \alpha_k$ for $\tau \in (0,1)$). This process is repeated until a steplength that satisfies \eqref{eq:armijo} is found. We should note that under reasonable assumptions on the function $F$ (see \cite{nocedal_book}) this procedure is well defined since the search direction uses the true gradient, $H_k$ is a positive definite matrix, and the true function is used in condition \eqref{eq:armijo}.

\subsection{Sampled LSR1}
At the $k$th iteration, the S-LSR1 method computes a new iterate via \eqref{eq:sr1}, where the Hessian approximation in \eqref{eq:tr_obj} is constructed using the curvature pairs sampled by Algorithm \ref{alg:calSY}. The S-LSR1 method is outlined in Algorithm  \ref{alg:slsr1}.

\begin{algorithm}[]
 {\small 
\caption{Sampled LSR1 (S-LSR1)}
  \label{alg:slsr1}
 {\bf Input:} $w_{0}$ (initial iterate),  
$m$ (memory), $r$ (sampling radius), $\Delta_0$ (initial trust region radius), $\eta_1 \in (0,1)$ (step acceptance parameter).

  \begin{algorithmic}[1]
  \For {$k=0,1,2,...$}
    \State Compute new $(S_k,Y_k)$ pairs via Algorithm \ref{alg:calSY}
    \State Compute $p_k$ by solving the subproblem \eqref{eq:tr_obj}
    \State Compute $\rho_k = \frac{F(w_k) - F(w_k + p_k)}{m_k(0) - m_k(p_k)}$
    \If {$\rho_k \geq \eta_1$}
        \State Set $w_{k+1} = w_k + p_k$
    \Else
        \State Set $w_{k+1} = w_k$
    \EndIf
    \State $\Delta_{k+1} = \texttt{adjustTR}(\Delta_{k},\rho_k)$ [see Appendix \ref{sec:tr_alg}]
\EndFor
  \end{algorithmic}
  }
\end{algorithm}

The S-LSR1 method has the same key features as S-LBFGS that differentiates it from the classical SR1 methods. The subroutine \texttt{adjustTR} (Step 10, Algorithm \ref{alg:slsr1}) adjusts the trust-region based on the progress made by the method. For brevity we omit the details of this subroutine, and refer the reader to  Appendix \ref{sec:impl_lsr1} for the details.

\section{Distributed Computing and Computational Cost}
\label{sec:dist_comp}

In this section, we show the scalability of the sampled quasi-Newton methods as compared to the SG method, and compare the computational cost to the classical variants. 

\subsection{Distributed Computing}

\begin{figure}[H]
    \centering
    \includegraphics[width=0.95\textwidth]{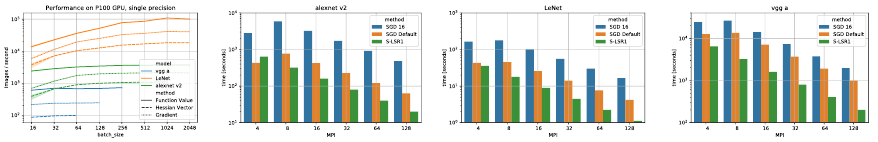}
\caption{\small Performance (Images/second) as a function of batch size for different DNN models and operations on a single P100 GPU (left). Time (seconds) to complete 1 epoch of SG and to perform 1 iteration of S-LSR1 on a dataset with 1M images using varying number of MPI processes (bar plots).}
\label{fig:scaling}
\end{figure}

Recently, there has been a huge effort to scale SG-type algorithms to solve \texttt{Imagenet} using hundreds of GPUs; see e.g., \cite{goyal2017accurate,jia2018highly,akiba2017extremely,you2018imagenet}. In Figure \ref{fig:scaling} (left), we show how the batch size affects the number of images processed per second to compute the function, gradient and Hessian vector products on a NVIDIA Tesla P100 GPU  for various deep neural networks\footnote{The structure of the deep neural network is taken from: \texttt{https://github.com/tensorflow/models/tree/master/research/slim}.}; see Table~\ref{tbl:networks}.

\begin{table}[H]
\caption{\small Deep Neural Networks used in the experiments.}
\label{tbl:networks}
\centering
\begin{small}
\begin{tabular}{@{}lrcr@{}}\toprule
{\bf model} & {\bf$d$\quad} & {\bf input} & {\bf \# classes}
 \\ \midrule
{\bf LeNet} & 3.2M & $28\times28\times3$ & 10
\\ \hdashline  
{\bf alexnet v2} & 50.3M & $224\times224\times3$ & 1,000
\\ \hdashline
{\bf vgg a} & 132.8M & $224\times224\times3$ & 1,000
\\
\bottomrule 
\end{tabular}
\end{small}
\end{table}

As is clear, by using small batch sizes one is not able to fully utilize the power of GPUs. On the other hand, using larger batches in conjunction with SG-type algorithms does not necessarily reduce training time \cite{das2016distributed,takac2013mini}.  Another observation that can be extracted from Figure \ref{fig:scaling} is that the cost of computing function values, gradients and Hessian vector products appears to be comparable for these networks.

In Figure \ref{fig:scaling} (bar plots), we compare the time to perform 1 epoch of the SG method (assuming we have 1M images) with the time to perform 1 iteration of S-LSR1. For SG, we show results for different batch sizes on each GPU\footnote{{ Each GPU has 1 MPI process that is used for communicating updates. Note, we are running 4 MPI processes for each physical node, i.e., each node has 4 P100 GPUs}}: (1) batch size 16 (SGD 16); and, (2) batch size 32, 64 and 128 for vgg a, LeNet and alexnet v2, respectively (SGD Default). The reason that there is no significant benefit when using more GPUs for the SG method is that the cost is dominated by the communication. For S-LSR1, that is not the case; as we scale up the number of MPI processes, we get good performance improvements since there is much less communication involved. See Appendix \ref{sec:A:Communication}  for more details.

\subsection{Cost, Storage and Parallelization}
The cost per iteration of the different quasi-Newton methods can be deconstructed as follows:
\begin{align}
\text{Cost} = {\substack{\text{Cost of gradient} \\ \text{computation}}} + {\substack{\text{Cost of forming/} \\ \text{taking step.}}}
\end{align}
Note, motivated by the results in Figure \ref{fig:scaling}, we assume that the cost computing a function value, gradient and Hessian vector product is comparable and is $\mathcal{O}(nd)$. The cost of computing the gradient is common for each method, whereas the search directions are computed differently for BFGS-type methods and SR1-type methods. More specifically, for BFGS methods we employ a line search  and for SR1 method we use a trust region and solve the subproblem \eqref{eq:tr_obj} using CG \cite{nocedal_book}. We denote the number of line search iterations and CG iterations as $\kappa_{ls}$ and $\kappa_{tr}$, respectively. Table \ref{tbl:cost_QN} summarizes the computational cost and storage for the different quasi-Newton methods.

As is clear from Table \ref{tbl:cost_QN}, the proposed sampled quasi-Newton methods do not have a significantly higher cost per iteration than the classical limited memory variants of the methods. In the regime where $m \ll n,d$, the computational cost of the methods are $\mathcal{O}(nd)$. Moreover, the storage requirements for the sampled quasi-Newton methods is the same as that of limited-memory quasi-Newton methods. We should also note, that several computations that are required in our proposed methods are easily parallelizeable. These computations are the gradient evaluations, the function evaluations and the construction of the gradient displacement curvature pairs $y$.
\begin{table}[]
\caption{\small Summary of Computational Cost and Storage (per iteration) for different Quasi-Newton methods.}
\label{tbl:cost_QN}
\centering
\begin{small}
\begin{tabular}{lcc}
\toprule
\textbf{method} & \textbf{computational cost} & \textbf{storage}  \\  \midrule
\textbf{BFGS} & $nd + d^2 + \kappa_{ls}nd$ & $d^2$   \\ \hdashline
\textbf{LBFGS} & $nd + 4md + \kappa_{ls}nd$ & $2md$  \\ \hdashline
\textbf{S-LBFGS} & $nd + mnd + 4md + \kappa_{ls}nd$ & $2md$   \\ \hline
\textbf{SR1} & $nd + d^2 + nd + \kappa_{tr}d^2$ & $d^2$   \\ \hdashline
\textbf{LSR1} & $nd + nd +\kappa_{tr}md$ & $2md$   \\
\hdashline
\textbf{S-LSR1} & $nd + mnd + nd +\kappa_{tr}md$ & $2md$  
\\ \bottomrule
\end{tabular}
\end{small}
\end{table}

\section{Convergence Analysis}
\label{sec:conv_analysis}

In this section, we present convergence analyses for the sampled quasi-Newton methods. 

\subsection{Sampled LBFGS}

We derive convergence results for the sampled LBFGS method with fixed step sizes and adaptive step sizes for strongly convex and nonconvex functions.

\subsubsection{Strongly Convex Functions}

We make the following standard assumptions.
\begin{assum} \label{ass:diff} $F$ is twice continuously differentiable.
\end{assum}
\begin{assum} \label{ass:strong_conv} There exist positive constants $\mu$ and $L$ such that
\begin{align*}
    \mu I \preceq \nabla^2F(w) \preceq L I, \quad \text{for all } w \in \mathbb{R}^d.
\end{align*}
\end{assum}

First, we show that the inverse Hessian approximations $H_k$ generated by the sampled
LBFGS method have eigenvalues that are uniformly bounded above and away from zero. The proof
technique is an adaptation of that in \cite{byrd2016stochastic,berahas2016multi}; however, modifications are necessary since in our approach the inverse Hessian approximations are constructed using information only from the current iterate, and not constructed sequentially.

\begin{lemma} \label{lem:1} If Assumptions \ref{ass:diff} and \ref{ass:strong_conv} hold, there exist constants $0 < \mu_1 \leq \mu_2$ such that the inverse Hessian approximations $\{ H_k\}$ generated by Algorithm \ref{alg:slbfgs} satisfy,
\begin{align}   \label{eq:bnd_Hess1}
\mu_1 I \preceq H_k \preceq \mu_2 I,\qquad \text{for } k=0,1,2,\dots.
\end{align}
\end{lemma}

\begin{proof} Instead of analyzing the inverse Hessian approximation $H_k$, we study the direct Hessian approximation $B_k = H_k^{-1}$. In this case, the sampled LBFGS updating formula is given as follows. 
At the $k$th iteration, given a set of curvature pairs $(s_{k,j},y_{k,j})$, for $j=1,\dots,m$
\begin{enumerate}
\item Set $B_k^{(0)}=\frac{y_{k,l}^Ty_{k,l}}{s_{k,l}^Ty_{k,l}}I$, where $l$ is chosen uniformly at random from $\{1,\dots,m \}$. 
\item For $i=1,\dots,m$ compute
\begin{align*}    
B_k^{(i)}=B_k^{(i-1)}-\frac{B_k^{(i-1)}s_{k,i}s_{k,i}^TB_k^{(i-1)}}{s_{k,i}^TB_k^{(i-1)}s_{k,i}} + \frac{y_{k,i}y_{k,i}^T}{y_{k,i}^Ts_{k,i}}.
\end{align*}
\item Set $B_{k+1} = B_k^{(m)}$. 
\end{enumerate}
In our algorithm (Algorithm \ref{alg:calSY}), there are two options for constructing the curvature pairs $s_{k,j}$ and $y_{k,j}$. At the current iterate $w_k$ we sample points $\bar{w}_j$ for $j=1,\dots,m$ and set
\begin{align}
    s_{k,j} = w_{k}-\bar{w}_j, \qquad y_{k,j}= \nabla F(w_k) - \nabla F(\bar{w}_j) \qquad \text{Option I},& \label{eq:curv_upd_opt1}\\
    s_{k,j} = w_{k}-\bar{w}_j, \qquad y_{k,j}=\nabla^2 F(w_k)s_k \qquad \text{Option II}.& \label{eq:curv_upd_opt2}
\end{align}
We now derive an upper and lower bound for $\frac{\|y_{k,j}\|^2}{y_{k,j}^Ts_{k,j}}$, for all $j = 1,\dots, m$, for both options.

\textbf{Option I:} A consequence of Assumption \ref{ass:strong_conv} is that the eigenvalues of the Hessian matrix are bounded above and away from zero. Utilizing this fact, the convexity of the objective function and the definitions  \eqref{eq:curv_upd_opt1}, we have
\begin{align} \label{eq:expr1}
y_{k,j}^Ts_{k,j} \geq \frac{1}{L}\|y_{k,j}\|^2  \quad & \Rightarrow \quad \frac{\|y_{k,j}\|^2 }{y_{k,j}^Ts_{k,j}} \leq L.
\end{align}
On the other hand, strong convexity of the functions, the consequence of Assumption \ref{ass:strong_conv} and definitions \eqref{eq:curv_upd_opt1}, provide a lower bound,
\begin{align} \label{eq:expr2}
y_{k,j}^Ts_{k,j} \leq \frac{1}{\mu}\|y_{k,j}\|^2  \quad & \Rightarrow \quad \frac{\|y_{k,j}\|^2 }{y_{k,j}^Ts_{k,j}} \geq \mu.
\end{align} 
Combining the upper and lower bounds \eqref{eq:expr1} and \eqref{eq:expr2}  
\begin{align}    \label{eq:bound}
\mu \leq \frac{\|y_{k,j}\|^2}{y_{k,j}^Ts_{k,j}} \leq L.
\end{align}

\textbf{Option II:} 
A consequence of Assumption \ref{ass:strong_conv} is that the eigenvalues of the Hessian matrix are bounded above and away from zero. Utilizing this fact and the definitions \eqref{eq:curv_upd_opt1}, we have
\begin{align} \label{eq:expr1_II}
    \mu \| s_{k,j} \|^2 \leq y_{k,j}^T s_{k,j} = s_{k,j}^T \nabla^2 F(w_k)s_{k,j} \leq L \| s_{k,j} \|^2.
\end{align}
We have that,
\begin{align}
    \frac{\|y_{k,j}\|^2}{y_{k,j}^T s_{k,j}} = \frac{s_{k,j}^T \nabla^2 F(w_k)^2 s_{k,j}}{s_{k,j}^T \nabla^2 F(w_k) s_{k,j}},
\end{align}
and since $\nabla^2 F(w_k)$ is symmetric and positive definite, it has a square root and so
\begin{align}    \label{eq:bound_II}
\mu \leq \frac{\|y_{k,j}\|^2}{y_{k,j}^Ts_{k,j}} \leq L.
\end{align}

The bounds on $\frac{\|y_{k,j}\|^2}{y_{k,j}^Ts_{k,j}}$ prove that for any $l$ chosen uniformly at random from $\{1,\dots,m \}$ the eigenvalues of the matrices $B_k^{(0)}=\frac{y_{k,l}^Ty_{k,l}}{s_{k,l}^Ty_{k,l}}I$ at the start of the sampled LBFGS update cycles are bounded above and away from zero, for all $k$ and $l$. We now use a Trace-Determinant argument to show that the eigenvalues of $B_k$ are bounded above and away from zero. 

Let $Tr(B)$ and $\det(B)$ denote the trace and determinant of matrix $B$, respectively. The trace of the matrix $B_{k+1}$ can be expressed as, 
\begin{align} \label{eq:trace}
Tr(B_{k+1}) &= Tr(B_k^{(0)}) - Tr\sum_{i=1}^{m}\left(\frac{B_k^{(i-1)}s_{k,i}s_{k,i}^TB_k^{(i-1)}}{s_{k,i}^TB_k^{(i-1)}s_{k,i}}\right) + Tr\sum_{i=1}^{m} \frac{y_{k,i}y_{k,i}^T}{y_{k,i}^Ts_{k,i}}\nonumber\\
&\leq Tr(B_k^{(0)}) + \sum_{i=1}^{m} \frac{\|y_{k,i}\|^2}{y_{k,i}^Ts_{k,i}}\nonumber\\
& \leq Tr(B_k^{(0)}) + mL \leq C_1, 
\end{align}
for some positive constant $C_1$, where the inequalities above are due to \eqref{eq:bound}, and the fact that the eigenvalues of the initial L-BFGS matrix $B_k^{(0)}$ are bounded above and away from zero.

Using a result due to Powell \cite{powell1976some}, the determinant of the matrix  $B_{k+1}$ generated by the sampled LBFGS method can be expressed as, 
\begin{align}	\label{eq:det}
\det (B_{k+1}) &= \det (B_{k}^{(0)}) \prod_{i=1}^{m} \frac{y_{k,i}^Ts_{k,i}}{s_{k,i}^TB_{k}^{(i-1)}s_{k,i}} \nonumber\\
& =  \det (B_{k}^{(0)}) \prod_{i=1}^{m} \frac{y_{k,i}^Ts_{k,i}}{s_{k,i}^Ts_{k,i}}  \frac{s_{k,i}^Ts_{k,i}}{s_{k,i}^TB_{k}^{(i-1)}s_{k,i}}\nonumber\\
& \geq \det (B_{k}^{(0)}) \Big( \frac{\mu}{C_1} \Big)^{m} \geq C_2,
\end{align}
for some positive constant $C_2$, where the above inequalities are due to the fact that the largest eigenvalue of $B_{k}^{(i)}$ is less than $C_1$, Assumption \ref{ass:strong_conv}, and the fact that $\frac{\mu}{C_1}<1$.

The trace \eqref{eq:trace} and determinant \eqref{eq:det} inequalities derived above imply that largest eigenvalues of all matrices $B_k$ are bounded above, uniformly, and that the smallest eigenvalues of all matrices $B_k$ are bounded away from zero, uniformly.
\end{proof}

\paragraph*{Constant Step Length} Utilizing Lemma \ref{lem:1}, we show that the sampled LBFGS method with a constant step length converges linearly. 

\begin{theorem} \label{thm:const_strongly}
Suppose that Assumptions \ref{ass:diff} and \ref{ass:strong_conv} hold, and let $F^{\star} = F(w^{\star})$, where $w^{\star}$ is the minimizer of $F$. Let $\{w_k\}$ be the iterates generated by Algorithm~\ref{alg:slbfgs}, where $0 <  \alpha_k = \alpha  \leq \frac{\mu_1}{\mu_2^2 L}$, and $w_0$ is the starting point. Then, for all $k\geq 0$, 
\begin{align*}   
 F(w_k) - F^{\star}& \leq \big( 1-\alpha \mu \mu_1 \big)^k \left[ F(w_0) - F^{\star} \right].
\end{align*}
\end{theorem}

\begin{proof} We have that
\begin{align} 
F(w_{k+1}) & = F(w_k -\alpha H_k \nabla F(w_k))  \nonumber\\
 & \leq F(w_k) + \nabla F(w_k)^T (-\alpha H_k \nabla F(w_k)) + \frac{L}{2}\| \alpha H_k \nabla F(w_k)\|^2  \nonumber\\
 & \leq F(w_k) - \alpha \nabla F(w_k)^T  H_k \nabla F(w_k) + \frac{\alpha^2 \mu_2^2 L}{2} \| \nabla F(w_k)\|^2  \nonumber\\
 & \leq F(w_k) - \alpha \mu_1 \| \nabla F(w_k) \|^2  + \frac{\alpha^2 \mu_2^2 L}{2} \| \nabla F(w_k)\|^2 \nonumber\\
 & = F(w_k) - \alpha \left(\mu_1 - \alpha \frac{\mu_2^2 L}{2} \right)\| \nabla F(w_k) \|^2  \label{eq:proof_step_1}\\
 & \leq F(w_k) - \alpha \frac{\mu_1}{2} \| \nabla F(w_k) \|^2,  \label{eq:proof_step}
\end{align}
where the first inequality is due to Assumption \ref{ass:strong_conv}, the second and third inequalities arise as a consequence of Lemma \ref{lem:1} and the last inequality is due to the choice of the steplength. By strong convexity, we have $2\mu(F(w)-F^\star) \leq \| \nabla F(w)\|^2$, and thus
\begin{align*}
    F(w_{k+1}) \leq F(w_k) - \alpha \mu \mu_1 (F(w_k)-F^\star).
\end{align*}
Subtracting $F^\star$ from both sides, 
\begin{align*}
    F(w_{k+1}) - F^\star \leq (1 - \alpha \mu \mu_1) (F(w_k)-F^\star).
\end{align*}
Recursive application of the above inequality yields the desired result.
\end{proof}

Theorem \ref{thm:const_strongly} shows that the S-LBFGS method  converges to the optimal solution at a linear rate. This result is similar in nature to the result for LBFGS \cite{liu1989limited}. We should also mention that unlike first-order methods (e.g., Gradient Descent, $H_k = I$), the step length range and the rate of convergence of the S-LBFGS method depends on $\mu_1$ and $\mu_2$, the smallest and largest eigenvalues of the S-LBFGS Hessian approximation. In the worst-case, the presence of the matrix $H_k$ can make the results in Theorem \ref{thm:const_strongly} significantly worse than that of the first-order variant if the update has been unfortunate and generates ill-conditioned matrices. We should note, however, such worst-case behavior is almost never observed in practice for BFGS updating.

\paragraph*{Adaptive Step Length} We now show a similar result for the case where the step length is chosen by an Armijo backtracking line search \eqref{eq:armijo}. 

\begin{theorem} \label{thm:linesearch_strongly}
Suppose that Assumptions \ref{ass:diff} and \ref{ass:strong_conv} hold. Let $\{w_k\}$ be the iterates generated by Algorithm~\ref{alg:slbfgs}, where $\alpha_k$ is the maximum value in $\{ \tau^{-j}:j=0,1,\dots\}$ satisfying \eqref{eq:armijo} with $0<c_1<1$, and $w_0$ is the starting point.
Then for all $k\geq 0$,
\begin{align*}   
 F(w_k) - F^{\star}& \leq \left( 1-\frac{4\mu \mu_1^2 c_1(1-c_1)\tau}{\mu_2^2L} \right)^k \left[ F(w_0) - F^{\star} \right].
\end{align*}
\end{theorem}

\begin{proof} Starting with \eqref{eq:proof_step_1} we have
\begin{align*} 
F(w_{k} - \alpha_k H_k \nabla F(w_k)) 
 & \leq F(w_k) - \alpha_k \left(\mu_1 - \alpha_k \frac{\mu_2^2 L}{2} \right)\| \nabla F(w_k) \|^2.  
\end{align*}
From the Armijo backtracking condition \eqref{eq:armijo}, we have \begin{align}
    F(w_{k} - \alpha H_k \nabla F(w_k)) &\leq F(w_k) - c_1 \alpha_k \nabla F(w_k)^T H_k \nabla F(w_k) \nonumber\\
    & \leq F(w_k) - c_1 \mu_1 \alpha_k \|\nabla F(w_k) \|^2. \label{eq:proof_armijo1}
\end{align}
Looking at \eqref{eq:proof_step_1} and \eqref{eq:proof_armijo1}, it is clear that the Armijo condition is satisfied when 
\begin{align} \label{eq.step_length_ls}
    \alpha_k \leq \frac{2\mu_1(1-c_1)}{\mu_2^2 L}.
\end{align}
Thus, any $\alpha_k$ that satisfies \eqref{eq.step_length_ls} is guaranteed to satisfy the Armijo condition \eqref{eq:armijo}. Since we find $\alpha_k$ using a constant backtracking factor of $\tau < 1$, we have
that 
\begin{align}\label{eq.eq.step_length_ls_lower}
    \alpha_k \geq \frac{2\mu_1(1-c_1)\tau}{\mu_2^2 L}.
\end{align}
Therefore, from \eqref{eq:proof_step_1} and by \eqref{eq.step_length_ls} and \eqref{eq.eq.step_length_ls_lower} we have
\begin{align}   \label{eq:proof_step_linesearch}
    F(w_{k+1}) &\leq F(w_k) - \alpha_k \left(\mu_1 - \alpha_k \frac{\mu_2^2 L}{2} \right)\| \nabla F(w_k) \|^2\nonumber\\
    & \leq F(w_k) - \alpha_k c_1 \mu_1 \| \nabla F(w_k) \|^2 \nonumber\\
    & \leq F(w_k) - \frac{2\mu_1^2c_1(1-c_1)\tau}{\mu_2^2L}\|\nabla F(w_k) \|^2.
\end{align}
By strong convexity, we have $2\mu(F(w)-F^\star) \leq \| \nabla F(w)\|^2$, and thus
\begin{align}   \label{eq:sc_linesearch}
    F(w_{k+1})
    & \leq F(w_k) - \frac{4\mu \mu_1^2 c_1(1-c_1)\tau}{\mu_2^2L}(F(w)-F^\star).
\end{align}
Subtracting $F^\star$ from both sides, and applying \eqref{eq:sc_linesearch} recursively yields the desired result.  
\end{proof}

Theorem \ref{thm:linesearch_strongly} shows that the sampled LBFGS method with an adaptive backtracking line search converges to the optimal solution at a linear rate. We should note that this result is worse than the constant step length result (Theorem \ref{thm:const_strongly}), i.e., the rate constant is larger. This is not surprising since this is a worst-case result; however, in practice, such an approach performs significantly better and does not require tuning the steplength parameter.

\subsubsection{Nonconvex Functions}
For nonconvex functions, the BFGS method is known fail \cite{mascarenhas2004bfgs,dai2002convergence}. Even for LBFGS, which makes only a finite number of updates at each iteration, one cannot guarantee that the (inverse) Hessian approximations have eigenvalues that are uniformly bounded above and away from zero. To establish convergence of the BFGS method in the nonconvex setting several techniques have been proposed including: $(i)$ \textit{cautious} updating \cite{li2001global}; $(ii)$ \textit{modified} updating \cite{li2001modified}; and $(iii)$ \textit{damping} \cite{powell1978algorithms}. Here we employ a cautious strategy that is well suited to our particular algorithm; at the $k$th iteration, we update the (inverse) Hessian approximation using only the set of curvature pairs that satisfy
\begin{align}   \label{curv}
    s^Ty > \epsilon \| s\|^2,
\end{align}
where $\epsilon > 0$ is a predetermined constant. If no curvature pairs satisfy \eqref{curv}, then the new (inverse) Hessian approximation is set to $H_{k} = I$. Using said mechanism we prove that the eigenvalues of the (inverse) Hessian approximations generated by the S-LBFGS method are bounded above and away from zero. For this analysis, we make the following assumptions in addition to Assumption \ref{ass:diff}.

\begin{assum} \label{ass:boundedf} The function $ F(w)$ is bounded below by a scalar $\widehat F$.
\end{assum}

\begin{assum} \label{ass:lip} The gradients of $F$ are $L$-Lipschitz continuous for all $w \in \mathbb{R}^d$.
\end{assum}

\begin{lemma}		\label{lemma3}
Suppose that Assumptions \ref{ass:diff} and \ref{ass:lip}  hold. Let $\{H_k \}$ be the inverse Hessian approximations generated by Algorithm~\ref{alg:slbfgs}, with the modification that the inverse  approximation update is performed using only curvature pairs that satisfy \eqref{curv},
for some $\epsilon >0$, and $H_{k} = I$ if no curvature pairs satisfy \eqref{curv}. Then, there exist constants $0<\mu_1\leq \mu_2$ such that 
\begin{align}    \label{eq:bnd_Hess2}
\mu_1 I \preceq H_k \preceq \mu_2 I,\qquad \text{for } k=0,1,2,\dots.
\end{align}
\end{lemma}

\begin{proof} Note, that in the nonconvex setting, there is a chance that no curvature pairs are selected in Algorithm \ref{alg:calSY}. In this case, the inverse Hessian approximation is $H_k = I$, and thus $\mu_1 = \mu_2 = 1$ and condition \eqref{eq:bnd_Hess2} is satisfied.

Similar to the proof of Lemma \ref{lem:1}, we study the direct Hessian approximation $B_k = H_k^{-1}$. In our algorithm, there are two options for updating the curvature pairs $s_{k,j}$ and $y_{k,j}$:
\begin{align}
    s_{k,j} = w_{k}-\bar{w}_j, \qquad y_{k,j}= \nabla F(w) - \nabla F(\bar{w}_j) \qquad \text{Option I},& \label{eq:curv_upd_opt1_2}\\
    s_{k,j} = w_{k}-\bar{w}_j, \qquad y_{k,j} =\nabla^2 F(w_k)s_k \qquad \text{Option II},& \label{eq:curv_upd_opt2_2}
\end{align}
for $j=1,\dots,m$. Let $\tilde{m}_k \in \{1,...,m \}$ denote the number of curvature pairs that satisfy \eqref{curv} at the $k$th iteration, where $m$ is the memory. At the $k$th iteration, given a set of curvature pairs $(s_{k,j},y_{k,j})$, for $j=1,\dots,\tilde{m}_k$ we update the Hessian approximation recursively (using the procedure described in the proof of Lemma \ref{lem:1}, and set $B_{k+1} = B_{k}^{\tilde{m}_k}$.

In this setting, 
the skipping mechanism \eqref{curv} provides both an upper and lower bound on the quantity $\frac{\|y_{k,j}\|^2 }{y_{k,j}^Ts_{k,j}}$, for both Options, which in turn ensures that the initial sampled LBFGS Hessian approximation is bounded above and away from zero. 

The lower bound is attained by repeated application of Cauchy's inequality to condition \eqref{curv}. We have from \eqref{curv} that
\begin{align*}
	\epsilon \| s_{k,j} \|^2 & < y_{k,j}^Ts_{k,j} \leq  \| y_{k,j} \| \| s_{k,j} \| \quad \Rightarrow \quad \| s_{k,j} \| < \frac{1}{\epsilon} \| y_{k,j} \|.
\end{align*} 
It follows that 
\begin{align}  \label{eq:lower}
	s_{k,j}^Ty_{k,j} \leq \| s_{k,j} \| \| y_{k,j} \| < \frac{1}{\epsilon} \| y_{k,j} \|^2 \quad \Rightarrow \quad \frac{\| y_{k,j} \|^2}{s_{k,j}^Ty_{k,j}} > \epsilon.
\end{align}

The upper bound is attained by the Lipschitz continuity of gradients,
\begin{align}
	 y_{k,j}^Ts_{k,j} & > \epsilon \| s_{k,j} \|^2 \nonumber\\
	&\geq  \epsilon \frac{ \| y_{k,j} \|^2}{L}  \quad \Rightarrow \quad \frac{\| y_{k,j} \|^2}{s_{k,j}^Ty_{k,j}} < \frac{L^2}{\epsilon} \label{eq:upper}.
\end{align} 

Combining \eqref{eq:lower} and \eqref{eq:upper}, we have
\begin{align*}
 \epsilon < \frac{\|y_{k,j}\|^2 }{y_{k,j}^Ts_{k,j}} < \frac{L^2}{\epsilon}.
\end{align*}

The bounds on $\frac{\|y_{k,j}\|^2}{y_{k,j}^Ts_{k,j}}$ prove that for any $l$ chosen uniformly at random from $\{1,\dots,\tilde{m}_k \}$ the eigenvalues of the matrices $B_k^{(0)}=\frac{y_{k,l}^Ty_{k,l}}{s_{k,l}^Ty_{k,l}}I$ at the start of the sampled LBFGS update cycles are bounded above and away from zero, for all $k$ and $l$. The rest of the proof follows the same trace-determinant argument as in the proof of Lemma \ref{lem:1}, the only difference being that the last inequality in \ref{eq:det} comes as a result of the cautious update strategy. 
\end{proof}

\paragraph*{Constant Step Length} Utilizing Lemma \ref{lemma3}, we show that the sampled LBFGS with a cautious updating strategy and a constant step length converges. 

\begin{theorem} \label{thm_non_MB}
Suppose that Assumptions \ref{ass:diff},  \ref{ass:boundedf} and \ref{ass:lip} hold. Let $\{w_k\}$ be the iterates generated by Algorithm~\ref{alg:slbfgs}, with the modification that the inverse Hessian approximation update is performed using only curvature pairs that satisfy \eqref{curv},
for some $\epsilon >0$, and $H_{k} = I$ if no curvature pairs satisfy \eqref{curv},  where 
$0 <  \alpha_k = \alpha  \leq \frac{\mu_1}{\mu_2^2  L}$,
and $w_0$ is the starting point. Then, 
\begin{align}   \label{eq:lim_res_non_const}
    \lim_{k\rightarrow \infty} \| \nabla F(w_k) \| = 0,
\end{align}
and, moreover, for any $T > 1$,
\begin{align*}	
\frac{1}{T}\sum_{k=0}^{T-1} \| \nabla F(w_k) \|^2  &  \leq \frac{2[ F(w_0) - \hat{F}]}{\alpha \mu_1 T } \xrightarrow[]{T \rightarrow \infty}0.
\end{align*}
\end{theorem}

\begin{proof} We start with \eqref{eq:proof_step} \begin{align*} 
F(w_{k+1})  \leq F(w_k) - \alpha \frac{\mu_1}{2} \| \nabla F(w_k) \|^2.  
\end{align*}
Summing both sides of the above inequality from $k=0$  to $T -1$,
\begin{align*} 
	\sum_{k=0}^{T-1}(F(w_{k+1}) - F(w_k)) \leq  - \sum_{k=0}^{T-1}\alpha \frac{\mu_1}{2} \| \nabla F(w_k) \|^2. 
\end{align*}
The left-hand-side of the above inequality is a telescopic sum and thus,
\begin{align*}
	\sum_{k=0}^{T-1} \left[F(w_{k+1})-F(w_k) \right]  &= F(w_{T})-F(w_0)
	\geq \widehat{F} -F(w_0),
\end{align*}
where the inequality is due to $\hat{F} \leq F(w_T)$ (Assumption \ref{ass:boundedf}). Using the above, we have
\begin{align}   \label{eq:int_nonconvex}
    \sum_{k=0}^{T-1} \| \nabla F(w_k) \|^2  &  \leq \frac{2[ F(w_0) - \widehat{F}]}{\alpha \mu_1}.
\end{align}
Taking limits we obtain,
\begin{align*}
    \lim_{T \rightarrow \infty} \sum_{k=0}^{T-1} \| \nabla F(w_k) \|^2 < \infty,
\end{align*}
which implies \eqref{eq:lim_res_non_const}. Dividing \eqref{eq:int_nonconvex} by $T$ we conclude
\begin{align*}
    \frac{1}{T}\sum_{k=0}^{T-1} \| \nabla F(w_k) \|^2  &  \leq \frac{2[ F(w_0) - \widehat{F}]}{\alpha \mu_1 T}.
\end{align*}
\end{proof}

Theorem \ref{thm_non_MB} shows that, if a small enough step length is chosen, 
the S-LBFGS method converges to a stationary point. 

\paragraph*{Adaptive Step Length} We now show a similar result for the case where the step length is chosen by an Armijo backtracking line search \eqref{eq:armijo}.

\begin{theorem} \label{thm_non_MB_linesearch}
Suppose that Assumptions \ref{ass:diff},  \ref{ass:boundedf} and \ref{ass:lip} hold. Let $\{w_k\}$ be the iterates generated by Algorithm~\ref{alg:slbfgs}, with the modification that the inverse Hessian approximation update is performed using only curvature pairs that satisfy \eqref{curv},
for some $\epsilon >0$, and $H_{k} = I$ if no curvature pairs satisfy \eqref{curv},  where $\alpha_k$ is the maximum value in $\{ \tau^{-j}:j=0,1,\dots\}$ satisfying \eqref{eq:armijo} with $0<c_1<1$, and where $w_0$ is the starting point. Then, 
\begin{align}   \label{eq:lim_res_non}
    \lim_{k\rightarrow \infty} \| \nabla F(w_k) \| = 0,
\end{align}
and, moreover, for any $T > 1$,
\begin{align*}	
\frac{1}{T}\sum_{k=0}^{T-1} \| \nabla F(w_k) \|^2  &  \leq \frac{ \mu_2^2 L[ F(w_0) - \widehat{F}]}{ 2\mu_1^2c_1(1-c_1)\tau T} \xrightarrow[]{\tau\rightarrow \infty}0.
\end{align*}
\end{theorem}

\begin{proof} We start with \eqref{eq:proof_step_linesearch} 
\begin{align*} 
F(w_{k+1}) & \leq F(w_k) - \frac{2\mu_1^2c_1(1-c_1)\tau}{\mu_2^2L}\|\nabla F(w_k) \|^2.
\end{align*}

Summing both sides of the above inequality from $k=0$  to $T -1$,
\begin{align*} 
	\sum_{k=0}^{T-1}(F(w_{k+1}) - F(w_k)) \leq  - \sum_{k=0}^{T-1}\frac{2\mu_1^2c_1(1-c_1)\tau}{\mu_2^2L} \| \nabla F(w_k) \|^2. 
\end{align*}
The left-hand-side of the above inequality is a telescopic sum and thus,
\begin{align*}
	\sum_{k=0}^{T-1} \left[F(w_{k+1})-F(w_k) \right]  &= F(w_{T})-F(w_0)
	\geq \widehat{F} -F(w_0),
\end{align*}
where the inequality is due to $\hat{F} \leq F(w_T)$ (Assumption \ref{ass:boundedf}). Using the above, we have
\begin{align}   \label{eq:int_nonconvex1}
    \sum_{k=0}^{T-1} \| \nabla F(w_k) \|^2  &  \leq \frac{\mu_2^2 L[ F(w_0) - \widehat{F}]}{2\mu_1^2c_1(1-c_1)\tau}.
\end{align}
Taking limits we obtain,
\begin{align*}
    \lim_{\tau \rightarrow \infty} \sum_{k=0}^{\tau-1} \| \nabla F(w_k) \|^2 < \infty,
\end{align*}
which implies \eqref{eq:lim_res_non}. Dividing \eqref{eq:int_nonconvex1} by $T$ we conclude
\begin{align*}
    \frac{1}{T}\sum_{k=0}^{T-1} \| \nabla F(w_k) \|^2  &  \leq \frac{\mu_2^2 L[ F(w_0) - \widehat{F}]}{ 2\mu_1^2c_1(1-c_1)\tau T}.
\end{align*}
\end{proof}

Theorem \ref{thm_non_MB_linesearch} shows that, the S-LBFGS method that employs an Armijo backtracking linesearch \eqref{eq:armijo}  converges to a stationary point. 

\paragraph*{Adaptive Step Length} We now show a similar result for the case where the step length is chosen by an Armijo backtracking line search \eqref{eq:armijo}.

\begin{theorem} \label{thm_non_MB_linesearch}
Suppose that Assumptions \ref{ass:diff},  \ref{ass:boundedf} and \ref{ass:lip} hold. Let $\{w_k\}$ be the iterates generated by Algorithm~\ref{alg:slbfgs}, with the modification that the inverse Hessian approximation update is performed using only curvature pairs that satisfy \eqref{curv},
for some $\epsilon >0$, and $H_{k} = I$ if no curvature pairs satisfy \eqref{curv},  where $\alpha_k$ is the maximum value in $\{ \tau^{-j}:j=0,1,\dots\}$ satisfying \eqref{eq:armijo} with $0<c_1<1$, and where $w_0$ is the starting point. Then, 
\begin{align}   \label{eq:lim_res_non}
    \lim_{k\rightarrow \infty} \| \nabla F(w_k) \| = 0,
\end{align}
and, moreover, for any $T > 1$,
\begin{align*}
\frac{1}{T}\sum_{k=0}^{T-1} \| \nabla F(w_k) \|^2  &  \leq \frac{ \mu_2^2 L[ F(w_0) - \widehat{F}]}{ 2\mu_1^2c_1(1-c_1)\tau T} \xrightarrow[]{\tau\rightarrow \infty}0.
\end{align*}
\end{theorem}

\begin{proof} We start with \eqref{eq:proof_step_linesearch} 
\begin{align*} 
F(w_{k+1}) & \leq F(w_k) - \frac{2\mu_1^2c_1(1-c_1)\tau}{\mu_2^2L}\|\nabla F(w_k) \|^2.
\end{align*}

Summing both sides of the above inequality from $k=0$  to $T -1$,
\begin{align*} 
	\sum_{k=0}^{T-1}(F(w_{k+1}) - F(w_k)) \leq  - \sum_{k=0}^{T-1}\frac{2\mu_1^2c_1(1-c_1)\tau}{\mu_2^2L} \| \nabla F(w_k) \|^2. 
\end{align*}
The left-hand-side of the above inequality is a telescopic sum and thus,
\begin{align*}
	\sum_{k=0}^{T-1} \left[F(w_{k+1})-F(w_k) \right]  &= F(w_{T})-F(w_0)
	\geq \widehat{F} -F(w_0),
\end{align*}
where the inequality is due to $\hat{F} \leq F(w_T)$ (Assumption \ref{ass:boundedf}). Using the above, we have
\begin{align}   \label{eq:int_nonconvex1}
    \sum_{k=0}^{T-1} \| \nabla F(w_k) \|^2  &  \leq \frac{\mu_2^2 L[ F(w_0) - \widehat{F}]}{2\mu_1^2c_1(1-c_1)\tau}.
\end{align}
Taking limits we obtain,
\begin{align*}
    \lim_{\tau \rightarrow \infty} \sum_{k=0}^{\tau-1} \| \nabla F(w_k) \|^2 < \infty,
\end{align*}
which implies \eqref{eq:lim_res_non}. Dividing \eqref{eq:int_nonconvex1} by $T$ we conclude
\begin{align*}
    \frac{1}{T}\sum_{k=0}^{T-1} \| \nabla F(w_k) \|^2  &  \leq \frac{\mu_2^2 L[ F(w_0) - \widehat{F}]}{ 2\mu_1^2c_1(1-c_1)\tau T}.
\end{align*}
\end{proof}

Theorem \ref{thm_non_MB_linesearch} shows that, the S-LBFGS method that employs an Armijo backtracking linesearch \eqref{eq:armijo}  converges to a stationary point. 

\subsection{Sampled LSR1}

We derive convergence results for the sampled SR1 method for general nonconvex objective functions. 

In order to establish convergence results one needs to ensure that the SR1 Hessian update equation \eqref{eq:sr1_hess} is well defined. To this end, we employ a cautious updating mechanism  that is well suited to our particular algorithm. At the $k$th iteration, we update the Hessian approximation using only the set of curvature pairs that satisfy
\begin{align}   \label{eq:cond_slsr1}
    | s^T(y-Bs) |  > \epsilon \| s\|^2,
\end{align}
where $\epsilon > 0$ is a predetermined constant. If no curvature pairs satisfy \eqref{eq:cond_slsr1}, then the new Hessian approximation is set to $B_{k} = I$. It is not trivial to test this condition in practice without explicitly constructing $d \times d$ matrices. We discuss this in detail in Section \ref{sec:num_res}; see Appendix \ref{sec:impl_lsr1} for more details.

For the analysis in this section, we make the following  assumption in addition to \ref{ass:diff}, \ref{ass:boundedf} and \ref{ass:lip}.
\begin{assum} \label{ass:cauchy_dec} For all $k$,
\begin{align*}
    m_k(0) - m_k(p_k) \geq \xi \| \nabla F(w_k)\| \min \left\{ \tfrac{\| \nabla F(w_k)\|}{\beta_k}, \Delta_k\right\},
\end{align*}
where $\xi \in (0,1)$ and $\beta_k = 1 + \| B_k\|$.
\end{assum}
Assumption \ref{ass:cauchy_dec} ensures that at every iteration we solve the trust-region subproblem sufficiently accurately.

We prove that the Hessian approximations $B_k$ generated by the S-LSR1 method are uniformly bounded from above. The proof technique is an adaptation of that in \cite{lu1996study}; however, modifications are necessary since the Hessian approximations are constructed using information only from the current iterate, and not constructed sequentially.

\begin{lemma} \label{lem:sr1} Suppose that Assumptions \ref{ass:diff}, \ref{ass:lip} and \ref{ass:cauchy_dec} hold. Let $\{ B_k\}$ be the Hessian approximations generated by Algorithm \ref{alg:slsr1}, with the modification that the approximation update is performed using only curvature pairs that satisfy \eqref{eq:cond_slsr1}, for some $\epsilon>0$, and $B_k = I$ if no curvature pairs satisfy \eqref{eq:cond_slsr1}. Then, there exists a constant $\nu_2 > 0$ such that
\begin{align}   \label{eq:bnd_Hess3}
    \|B_k\| \leq \nu_2, \qquad \text{for } k=0,1,2,\dots.  
\end{align}
\end{lemma}

\begin{proof} As in the proof of Lemma \ref{lemma3}, note that there is a chance that no curvature pairs are selected in Algorithm \ref{alg:calSY}. In this case, the Hessian approximation is $B_k = I$, and thus $ \nu_2 = 1$ and condition \eqref{eq:bnd_Hess3} is satisfied.

We now consider the case where at least one curvature pair is selected by Algorithm \ref{alg:calSY}. In this case, the sampled LSR1 updating formula is given as follows. Let $\tilde{m}_k \in \{1,...,m \}$ denote the number of curvature pairs that satisfy \eqref{eq:cond_slsr1} at the $k$th iteration, where $m$ is the memory. At the $k$th iteration, given a set of curvature pairs $(s_{k,j},y_{k,j})$, for $j=1,\dots,\tilde{m}_k$
\begin{enumerate}
\item Set $B_k^{(0)} = \gamma_kI$, where $0\leq \gamma_k < \gamma$. 
\item For $i=1,\dots,\tilde{m}_k$ compute
\begin{align*}    
    B_{k}^{(i)} = B_k^{(i-1)} +  \frac{(y_{k,i} - B_k^{(i-1)}s_{k,i})(y_{k,i} - B_k^{(i-1)}s_{k,i})^T}{(y_{k,i} - B_k^{(i-1)}s_{k,i})^Ts_{k,i}}.
\end{align*}
\item Set $B_{k+1} = B_k^{(\tilde{m}_k)}$. 
\end{enumerate}

In our algorithm (Algorithm \ref{alg:calSY}), there are two options for constructing the curvature pairs $s_{k,j}$ and $y_{k,j}$. At the current iterate $w_k$ we sample points $\bar{w}_j$ for $j=1,\dots,m$ and set
\begin{align}
    s_{k,j} = w_{k}-\bar{w}_j, \qquad y_{k,j}= \nabla F(w_k) - \nabla F(\bar{w}_j) \qquad \text{Option I},& \label{eq:curv_upd_opt1_sr1}\\
    s_{k,j} = w_{k}-\bar{w}_j, \qquad y_{k,j}=\nabla^2 F(w_k)s_k \qquad \text{Option II}.& \label{eq:curv_upd_opt2_sr1}
\end{align}
Given a set of $\tilde{m}_k$ curvature pairs that satisfy \eqref{eq:cond_slsr1}, we now prove an upper bound for $\| B_k\|$. We first prove the bound for a given iteration $k$ and for all updates to the Hessian approximation $i=0,1,\dots,\tilde{m}_k$ ($\| B_k^{i}\|$), and then get an upper bound for all $k$ ($\| B_k\|$). 

For a given iteration $k$, we prove a bound on $\| B_k^{i}\|$ via induction, and show 
\begin{align}   \label{eq:sr1_ind}
    \| B_k^{(i)}\| \leq \left(1 + \frac{1}{\epsilon} \right)^i \gamma_k + \left[ \left( 1 + \frac{1}{\epsilon} \right)^i - 1\right] \bar{\gamma}_k,
\end{align}
where $\bar{\gamma}_k$ is such that $\|\nabla^2 F(w_k)\| \leq \bar{\gamma}_k$, and whose existence follows from Assumption \ref{ass:lip} ($\bar{\gamma}_k \leq L < \infty$).
For $i=0$, the bound holds trivially since $B_k^{(0)} = \gamma_kI$. Now assume that \eqref{eq:sr1_ind} holds true for some $i\geq 0$. Note that all the curvature pairs that are used in the update of the Hessian approximation satisfy \eqref{eq:cond_slsr1}. By the definition of the SR1 updates, we have for some index $i+1$ that
\begin{align*}
    B_{k}^{(i+1)} = B_k^{(i)} +  \frac{(y_{k,i+1} - B_k^{(i)}s_{k,i+1})(y_{k,i+1} - B_k^{(i)}s_{k,i+1})^T}{(y_{k,i+1} - B_k^{(i)}s_{k,i+1})^Ts_{k,i+1}},
\end{align*}
and thus
\begin{align*}
    \| B_{k}^{(i+1)}\| &\leq \| B_k^{(i)} \| +  \left\Vert\frac{(y_{k,i+1} - B_k^{(i)}s_{k,i+1})(y_{k,i+1} - B_k^{(i)}s_{k,i+1})^T}{(y_{k,i+1} - B_k^{(i)}s_{k,i+1})^Ts_{k,i+1}} \right\Vert,\\
    &\leq \| B_k^{(i)} \| + \frac{\| (y_{k,i+1} - B_k^{(i)}s_{k,i+1})(y_{k,i+1} - B_k^{(i)}s_{k,i+1})^T\|}{\epsilon \| y_{k,i+1} - B_k^{(i)}s_{k,i+1}\| \|s_{k,i+1} \|}\\
    &\leq \| B_k^{(i)} \| + \frac{\|y_{k,i+1} - B_k^{(i)}s_{k,i+1}\|}{\epsilon \|s_{k,i+1} \|}\\
    &\leq \| B_k^{(i)} \| + \frac{\|y_{k,i+1}\|}{\epsilon \|s_{k,i+1} \|} + \frac{\|B_k^{(i)}s_{k,i+1}\|}{\epsilon \|s_{k,i+1} \|}\\
    &\leq \| B_k^{(i)} \| + \frac{\|y_{k,i+1}\|}{\epsilon \|s_{k,i+1} \|} + \frac{\|B_k^{(i)}\|}{\epsilon}\\
    & = \left( 1 + \frac{1}{\epsilon}\right)\| B_k^{(i)} \| + \frac{\bar{\gamma_k}}{\epsilon}
\end{align*}
where the first inequality is due to the application of the triangle inequality, the second inequality is due to condition \eqref{eq:cond_slsr1}, the fourth inequality is due to the application of the triangle inequality, and the fifth inequality is due to application of Cauchy's inequality and in the last inequality we used that $\bar{\gamma_k} \geq \bar{\gamma}_{k,i+1} = \frac{\|y_{k,i+1}\|}{ \|s_{k,i+1} \|}>0$. Substituting \eqref{eq:sr1_ind}, 
\begin{align*}
    \| B_{k}^{(i+1)}\| & \leq \left( 1 + \frac{1}{\epsilon}\right) \left[ \left(1 + \frac{1}{\epsilon} \right)^i \gamma_k + \left[ \left( 1 + \frac{1}{\epsilon} \right)^i - 1\right] \bar{\gamma}_k\right] + \frac{\bar{\gamma_k}}{\epsilon}\\
    & = \left(1 + \frac{1}{\epsilon} \right)^{i+1} \gamma_k + \left[ \left( 1 + \frac{1}{\epsilon} \right)^{i+1} - 1\right] \bar{\gamma}_k
\end{align*}
which completes the inductive proof. Thus, for any $k$ we have an upper bound on the Hessian approximation. Therefore, since $B_{k+1} = B_k^{(\tilde{m}_k)}$, the sampled SR1 Hessian approximation constructed at the $k$th iteration satisfies
\begin{align*}
    \| B_{k+1}\| & \leq \left(1 + \frac{1}{\epsilon} \right)^{\tilde{m}_k} \gamma_k + \left[ \left( 1 + \frac{1}{\epsilon} \right)^{\tilde{m}_k} - 1\right] \bar{\gamma}_k.
\end{align*}

Now we generalize the result for all iterations $k$. For $k=0$, the bound holds trivially, since the first step of the sampled LSR1 method is a gradient method ($B_0 = I$). For $k\geq 1$, we assume that $\gamma_k \leq \gamma < \infty$ and $\bar{\gamma}_k \leq \bar{\gamma} \leq L  < \infty$ for all $k$, and thus
\begin{align*}
    \| B_{k+1}\| & \leq \left(1 + \frac{1}{\epsilon} \right)^{\tilde{m}_k} \gamma_k + \left[ \left( 1 + \frac{1}{\epsilon} \right)^{\tilde{m}_k} - 1\right] \bar{\gamma}_k\\
    & \leq \left(1 + \frac{1}{\epsilon} \right)^{\tilde{m}_k} \gamma + \left[ \left( 1 + \frac{1}{\epsilon} \right)^{\tilde{m}_k} - 1\right] \bar{\gamma} \leq \nu_2,
\end{align*}
for some $\nu_2>0$. This completes the proof.
\end{proof}

Utilizing Lemma \ref{lem:sr1}, we show that the S-LSR1 with a cautious updating strategy converges. In order to prove the following result, we make use of well-known results for Trust-Region methods; see \cite{conn2000trust}. As such, the proof is identical to  \cite[Theorem 6.4.5]{conn2000trust}; to keep the paper self contained and due to the notation differences we include the proof.

\begin{theorem} \label{thm2:sr1} Suppose that Assumptions \ref{ass:diff}, \ref{ass:boundedf},  \ref{ass:lip} and \ref{ass:cauchy_dec} hold. Let $\{ w_k\}$ be the iterates generated by Algorithm \ref{alg:slsr1}, with the modification that the Hessian approximation update is performed using only curvature pairs that satisfy (\ref{eq:cond_slsr1}), for some $\epsilon>0$, and $B_k = I$ if no curvature pairs satisfy \eqref{eq:cond_slsr1}.  Then,  
\begin{align*}
{\lim_{k\rightarrow \infty}} \| \nabla F(w_k) \| = 0.
\end{align*}
\end{theorem}

\begin{proof} Assume, for the purpose of establishing a contradiction, that there is a subsequence of successful iterations (where $\rho_k > \eta_1$, Line 6, Algorithm \ref{alg:slsr1}), indexed by $t_i \subseteq \mathcal{S}$ where $\mathcal{S} = \{ k \geq 0 | \rho_k \geq \eta_1 \}$, such that
\begin{align}   \label{eq:tr_sub}
    \| \nabla F(w_{t_i})\| \geq 2\delta > 0
\end{align}
for some $\epsilon > 0$ and for all $i$. Theorem 6.4.5 from \cite{conn2000trust} then ensures the existence for each $t_i$ of a first successful iteration $\ell(t_i) > t_i$ such that
\begin{align*}
    \| \nabla F(w_{\ell(t_i)})\| < \delta > 0.
\end{align*}
Let $\ell_i = \ell(t_i)$, we thus obtain that there is anotehr subsequence of $\mathcal{S}$ indexed by $\{ \ell_i\}$ such that 
\begin{align}   \label{eq:tr1}
    \| \nabla F(w_k)\| \geq \delta, \quad \text{for} \quad t_i \leq k < \ell_i \quad \text{and} \quad \| \nabla F(w_{\ell_i})\| < \delta.
\end{align}
We now restrict our attention to the subsequence of successful iterations whose indices are in the set
\begin{align*}
    \mathcal{K} = \{k \in \mathcal{S} |  t_i \leq k < \ell_i \},
\end{align*}
where $t_i$ and $\ell_i$ belong to the subsequences $\mathcal{S}$ and $\mathcal{K}$, respectively.

Using Assumption \ref{ass:cauchy_dec}, the fact that $\mathcal{K} \subseteq \mathcal{S}$ and \eqref{eq:tr1}, we deduce that for $k \in \mathcal{K}$
\begin{align}   \label{eq:tr2}
    F(w_k) - F(w_k) \geq \eta_1[m_k(0) - m_k(p_k)] \geq \xi \delta \eta_1 \min \left[\frac{\delta}{\nu_2 + 1}, \Delta_k\right]
\end{align}
where we used the result of Lemma \ref{lem:sr1}. Since the sequence $\{ F(w_k)\}$ is monotonically decreasing and bounded below (Assumption \ref{ass:boundedf}), it is convergent, and the left-hand-side of \eqref{eq:tr2} must tend to zero as $k \rightarrow \infty$. Thus, 
\begin{align}
    \lim_{k \rightarrow \infty, \; k \in \mathcal{K}} \Delta_k = 0.
\end{align}

As a consequence, the term containing $\Delta_k$ is the dominant term in the $\min$ \eqref{eq:tr2} and we have, for $k \in \mathcal{K}$ sufficiently large,
\begin{align}
    \Delta_k \leq \frac{F(w_k) - F(w_{k+1})}{(\nu_2 + 1)\delta \eta_1}.
\end{align}
From this bound, we deduce that, for $i$ sufficiently large
\begin{align}   \label{eq:tr3}
    \|w_{t_i} - w_{\ell_i}\| \leq \sum_{j = t_i, \; j \in \mathcal{K}}^{\ell_i -1}  \|w_{j} - w_{j+1}\| \leq \sum_{j = t_i, \; j \in \mathcal{K}}^{\ell_i -1} \Delta_j \leq \frac{F(w_{t_i}) - F(w_{\ell_i})}{(\nu_2 + 1)\delta \eta_1}.
\end{align}
As a consequence of Assumption \ref{ass:boundedf} and the monotonicity of the sequence $\{ F(w_k)\}$, we have that the right-hand-side of \eqref{eq:tr3} must converge to zero, and thus $\|w_{t_i} - w_{\ell_i}\| \rightarrow 0$ as $i \rightarrow \infty$. 

By continuity of the gradient (Assumption \ref{ass:diff}), we thus deduce that $\| \nabla F(w_{t_i}) - \nabla F(w_{\ell_i})\| \rightarrow 0$. However, this is impossible because of the definitions of $\{ t_i\}$ and $\{ \ell_i\}$, which imply that  $\| \nabla F(w_{t_i}) - \nabla F(w_{\ell_i})\| \geq \delta$. Hence, no subsequence satisfying \eqref{eq:tr_sub} can exist, and the theorem is proved.
\end{proof}

Theorem \ref{thm2:sr1} shows that the sampled SR1 method converges to a stationary point. This result is similar in nature to that of the LSR1 method; see \cite{lu1996study}.

\section{Probabilistic Bounds on Sampled Quasi-Newton Methods}\label{sec:prob}

Since our proposed methods randomly select $m$ curvature pairs ($\{s,y\}$) at every iteration and we require the pairs satisfy certain conditions (\eqref{curv} and \eqref{eq:cond_slsr1} for S-LBFGS and S-LSR1, respectively), a fair question to ask is how many pairs are accepted and used to construct Hessian approximations at every iteration. Alternatively, the question can be posed as what is the probability that a given random $\{s,y\}$ pair satisfies the required conditions and is used in the quasi-Newton Hessian approximations.

In this section, we present probabilistic bounds that illustrate the probability of accepting a given $\{s,y\}$ pair. To do this, we leverage the form of Option II for computing the $y$ vectors (given a vector $s$) and the fact that $s$ can be any random vector. We will assume throughout this section that $s$ is uniformly sampled on a unit sphere, i.e., $s \sim \mathcal{U}(\mathcal{S}(0,1))$. We also illustrate the probabilities of accepting pairs empirically for synthetic problems with different dimensions and acceptance tolerances, and on two toy classification problem.

\subsection{Probabilistic Bounds for S-LBFGS}\label{sec.prob_SLBFGS}

In this section, we present results that illustrate the probability that the pairs generated within the S-LBFGS method satisfy \eqref{curv}. We first derive an expression for the probability of accepting a pair $\{s,y\}$, and then provide some empirical evidence to show the probability of accepting pairs for different problems.

By Option II, \eqref{curv} can be expressed as 
\begin{align*}
	\frac{s^T y}{\| s\|^2} = \frac{s^T \nabla^2 F(w)s}{\| s\|^2} > \epsilon,
\end{align*}
for any $w \in \mathbb{R}^d$. Notice that the middle term above is the \emph{Raleigh quotient} of a random vector $s$ with respect to the Hessian matrix. Thus, for any $w \in \mathbb{R}^d$ and any given random vector $s$, we are interested in the following probability,
\begin{align*}
	\mathbb{P}\left[ \frac{s^T \nabla^2 F(w)s}{\| s\|^2} > \epsilon \right] = 1 - \mathbb{P}\left[ \frac{s^T \nabla^2 F(w)s}{\| s\|^2} \leq \epsilon \right].
\end{align*}

The following theorem gives an expression for the probability of accepting the pair $\{s,y\}$.
\begin{theorem} \label{thm:probSLBFGS} Let $\lambda = (\lambda_1, \lambda_2, \dots, \lambda_d)$ be the eigenvalues of the true Hessian at some point $w \in \mathbb{R}^d$ ($\nabla^2 F(w)$), $s \in \mathbb{R}^d$ be a random vector uniformly distributed on a sphere, and $\epsilon>0$ be a prescribed tolerance. Then, 
\begin{align}   \label{cdfSLBFGS}
    \mathbb{P}\left[ \frac{s^T\nabla^2F(w) s}{\| s\|^2} > \epsilon \right] = \frac{1}{2} +  
\frac{1}{\pi}\int_0^{\infty} \frac{	\sin \left(\frac{1}{2}  \sum_{l=1}^d \tan^{-1} \left((\lambda_l - \epsilon )u \right)  \right)}{u \prod_{l=1}^d \left(1+(\lambda_l-\epsilon)^2 u^2\right)^{\frac{1}{4}}}du.
\end{align}
\end{theorem}

\begin{proof} We refer interested readers to \cite[Theorem 9]{boman1999infeasibility} for the proof of this theorem. 
\end{proof}

\begin{figure}[ht]
\centering
\subfloat[\label{fig.prob_1}]{%
\resizebox*{0.24\textwidth}{!}{\includegraphics[trim=0.6cm 7cm 3cm 7cm,width=\textwidth]{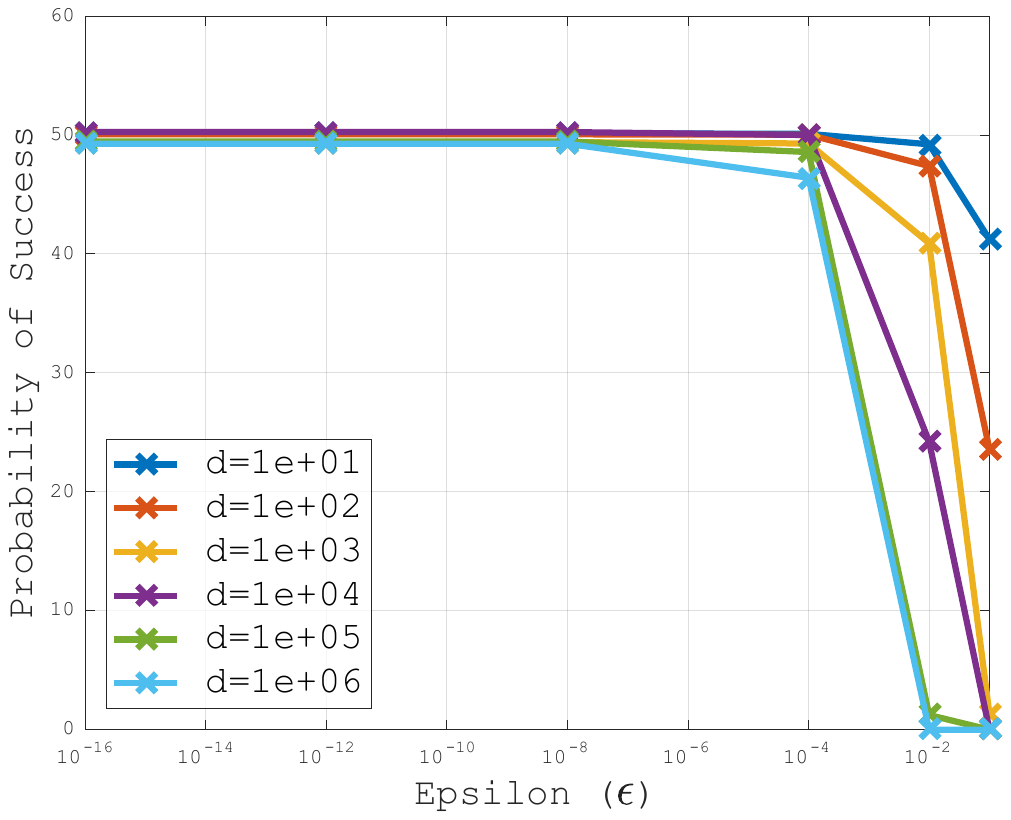}}}
\hspace{2pt}
\subfloat[\label{fig.prob_2}]{%
\resizebox*{0.24\textwidth}{!}{\includegraphics[trim=0.6cm 7cm 3cm 7cm,width=\textwidth]{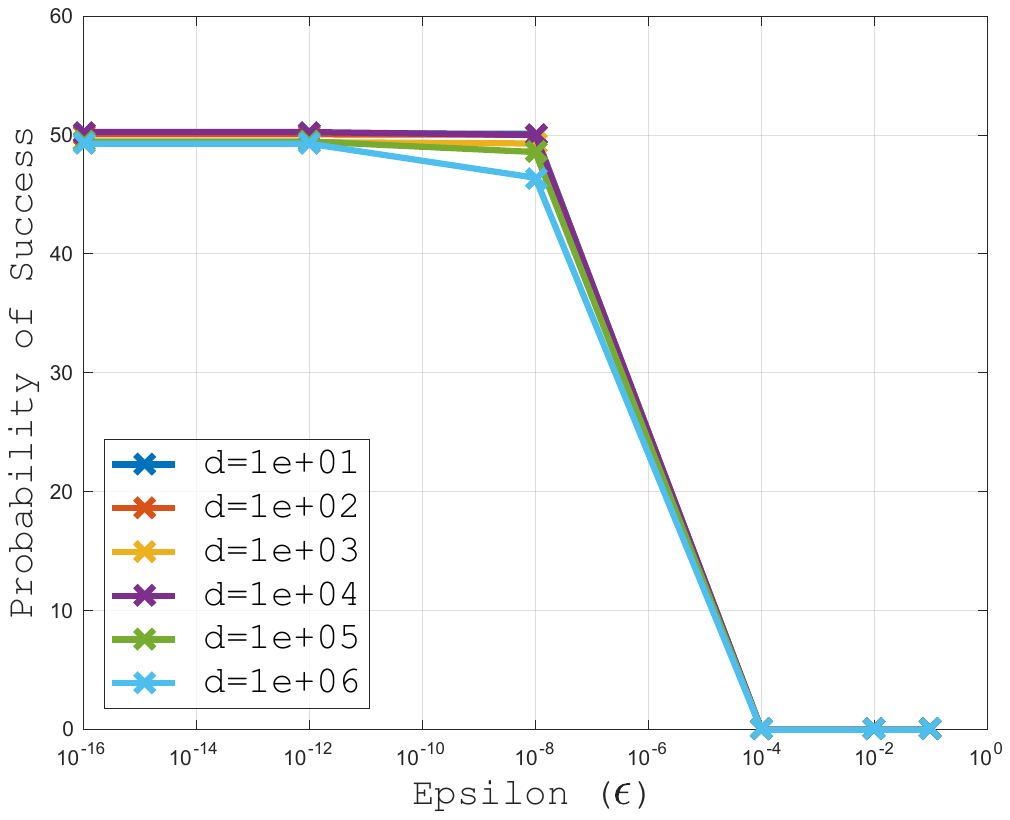}}}
\hspace{2pt}
\subfloat[\label{fig.prob_3}]{%
\resizebox*{0.24\textwidth}{!}{\includegraphics[trim=0.6cm 7cm 3cm 7cm,width=\textwidth]{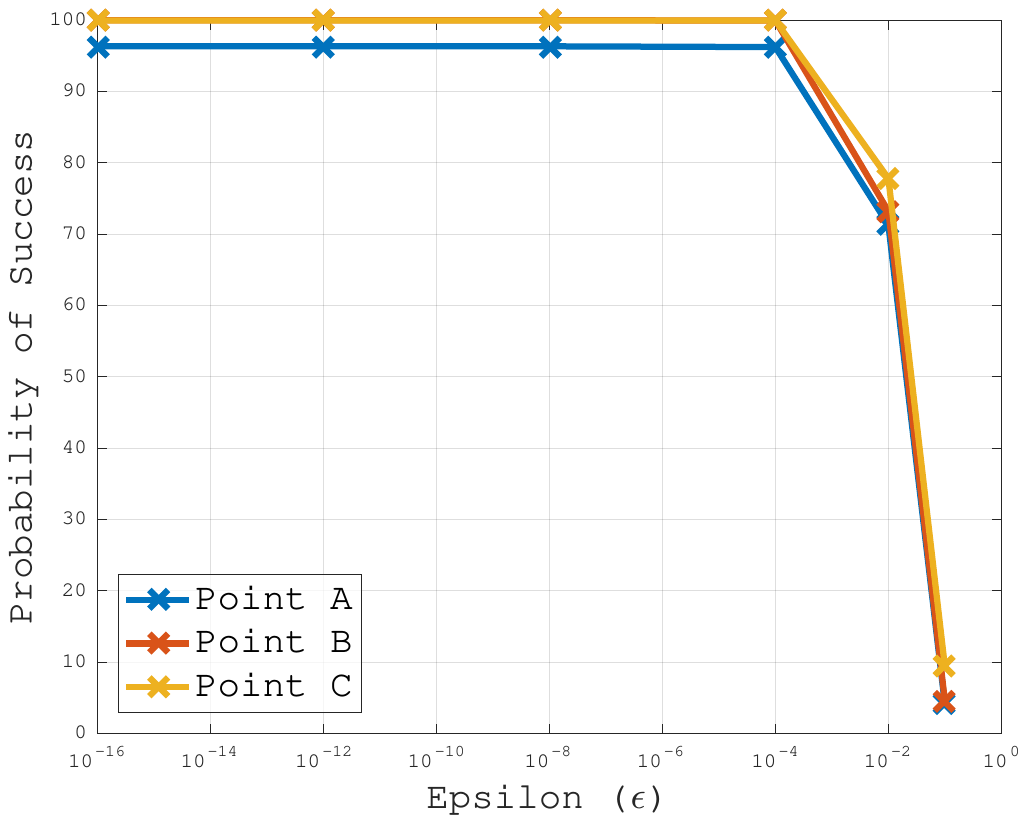}}}
\hspace{2pt}
\subfloat[\label{fig.prob_4}]{%
\resizebox*{0.24\textwidth}{!}{\includegraphics[trim=0.6cm 7cm 3cm 7cm,width=\textwidth]{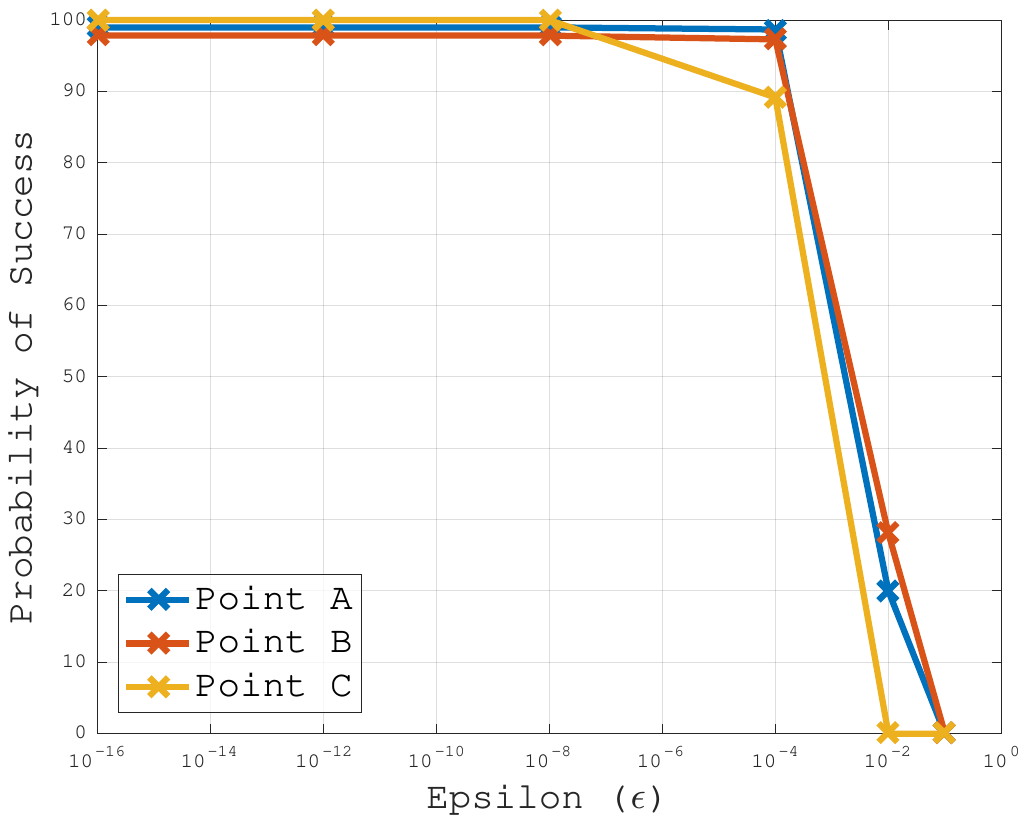}}}
\caption{\small Empirical investigation of accepting curvature pairs for different problems (S-LBFGS).\label{fig.prob_results_slbfgs}}
\end{figure}

Although the result of Theorem \ref{thm:probSLBFGS} is interesting, a reasonable criticism is that it requires knowledge of the eigenvalues of the Hessian matrix, something that is prohibitively expensive to compute for many deep learning training problems. However, we show numerically that for certain neural network problems, the probability of accepting pairs is relatively high. Specifically, we present empirical results that illustrate the probability that a given pair $\{s,y\}$ is accepted for different problems and different $\epsilon$. The problems considered are summarized in Table \ref{tbl:prob_details}, and the results are given in Figure \ref{fig.prob_results_slbfgs}.

\begin{table}[H]
\caption{\small Problem details for empirical evaluation of probabilities.}
\label{tbl:prob_details}
\centering
\begin{small}
\begin{tabular}{cc}
\toprule
\textbf{problem} & \textbf{figure}  \\ \midrule
$\lambda = (\underbrace{-1,\dots,-1}_{\text{$d/2$ -1's}},\underbrace{1,\dots,1}_{\text{$d/2$ 1's}})$   & Figures \ref{fig.prob_1} \& \ref{fig.prob_1_sr1}      \\ \hdashline
$\lambda = 10^{-4}(\underbrace{-1,\dots,-1}_{\text{$d/2$ -1's}},\underbrace{1,\dots,1}_{\text{$d/2$ 1's}})$        & Figures \ref{fig.prob_2} \& \ref{fig.prob_2_sr1} \\ \hdashline
Toy Problem: (\texttt{small network})     & Figures \ref{fig.prob_3} \& \ref{fig.prob_3_sr1}    \\ \hdashline
Toy Problem: (\texttt{medium network})     & Figures \ref{fig.prob_4} \& \ref{fig.prob_4_sr1}       \\ \bottomrule
\end{tabular}
\end{small}
\end{table}


The first two problems (Table \ref{tbl:prob_details}) have synthetic eigenvalue distributions. The goal of these problems is to investigate the effect of $n$ and $\epsilon$ on the probability. As is clear, for the first problem, the probability of accepting a random $\{s,y\}$ pair is around $50\%$ and the probability decreases last as $\epsilon \rightarrow 1$, which is not surprising due to the eigenvalue structure. For the second problem, where the eigenvalues are smaller, the probability becomes almost zero for $\epsilon \leq 10^{-4}$. For the last two problems, we considered the structures of the toy classification problems (small and medium) for different points in parameter space. Note, points A, B and C are the same as those used in Figures \ref{SmalllvdNet} and \ref{MediumvdNet}. As is clear from the Figures \ref{fig.prob_3} and \ref{fig.prob_4}, the probability of accepting a random curvature pair is very high as long as $\epsilon$ is not too large. The main takeaway from these numerical results is that the probability of accepting curvature pairs is relatively large as long as the tolerance is not chosen to be too large (in practice $\epsilon \approx 10^{-4}-10^{-8}$).

\subsection{Probabilistic Bounds on S-LSR1}

In this section, we present results that illustrate the probability that the pairs generated within the S-LSR1 method satisfy \eqref{eq:cond_slsr1}. We first derive an expression for the probability of accepting a pair $\{s,y\}$, and then provide some empirical evidence to show the probability of accepting pairs for different problems.

By Option II, \eqref{eq:cond_slsr1} can be expressed as 
\begin{align} \label{eq:sr1_prob}
	\frac{|s^T(y - Bs)|}{\| s\|^2} = \frac{|s^T( \nabla^2 F(w) - B)s|}{\| s\|^2} > \epsilon,
\end{align}
for any $w \in \mathbb{R}^d$, where the matrix $B$ is some SR1 Hessian approximation. Clearly, the acceptance of a new pair $\{s,y\}$ depends on the matrix $B$. To be more precise, for some $w \in \mathbb{R}^d$, given a new pair $\{s_j,y_j\}$ for $j=1,\dots,m$, \eqref{eq:sr1_prob} can be expressed as
\begin{align} \label{eq:sr1_prob_1}
	 \frac{|s_j^T( \nabla^2 F(w) - B^{(j-1)})s_j|}{\| s_j\|^2} > \epsilon,
\end{align}
where $B^{(j)}$ is an SR1 Hessian approximation constructed using all the pairs $\{s_i,y_i \}_{i<j}$, and $B^{(0)}$ is the initial SR1 Hessian approximation (potentially $B^{(0)} = 0$). As is clear, the acceptance of the new pair $\{s_j,y_j\}$ depends recursively on all previously accepted curvature pairs. Similar to the S-LBFGS case, the left-hand-side of \eqref{eq:sr1_prob_1} is a \emph{Raleigh quotient} of a random vector $s_j$ with respect to the Hessian matrix and the matrix $B^{(j-1)}$. 

Thus, for any $w \in \mathbb{R}^d$, $B^{(j-1)} \in \mathbb{R}^{d \times d}$ and any given random vector $s_j \in \mathbb{R}^d$, we are interested in the following probability,
\begin{align*}
	    \mathbb{P}\left[ \frac{|s_j^T (\nabla^2 F(w) - B^{(j-1)})s_j|}{\| s_j\|^2} > \epsilon \right] &= 1 - \mathbb{P}\left[ \frac{|s_j^T (\nabla^2 F(w) - B^{(j-1)})s_j|}{\| s_j\|^2} \leq \epsilon \right] \\
	    &= 1 - \mathbb{P}\left[ - \epsilon \leq \frac{s_j^T (\nabla^2 F(w) - B^{(j-1)})s_j}{\| s_j\|^2} \leq \epsilon \right] \\
	    & = 1 - \left( \mathbb{P}\left[ \frac{s_j^T (\nabla^2 F(w) - B^{(j-1)})s_j}{\| s_j\|^2} \leq \epsilon \right] \right. \\
	    & \qquad - \left.\mathbb{P}\left[ \frac{s_j^T (\nabla^2 F(w) - B^{(j-1)})s_j}{\| s_j\|^2} \leq -\epsilon \right] \right).
\end{align*}

The following theorem gives an expression for the probability of accepting the pair $\{s_j,y_j\}$.
\begin{theorem} \label{thm:probSLSR1} Let $\tilde{\lambda}^{j-1} = (\tilde{\lambda}_1^{j-1}, \tilde{\lambda}_2^{j-1}, \dots, \tilde{\lambda}_d^{j-1})$ be the eigenvalues of the matrix $\nabla^2 F(w)- B^{(j-1)}$, $s_j \in \mathbb{R}^d$ be a random vector uniformly distributed on a sphere, and $\epsilon>0$ be a prescribed tolerance. Then, for all $j \in \{1,\dots,m \}$,
\begin{align*}  
    \mathbb{P}\left[ \frac{|s_j^T (\nabla^2 F(w) - B^{(j-1)})s_j|}{\| s_j\|^2} > \epsilon \right] &= 1 + \frac{1}{\pi}\int_0^{\infty} \frac{	\sin \left(\frac{1}{2}  \sum_{l=1}^d \tan^{-1} \left((\tilde{\lambda}_l^{j-1} - \epsilon )u\right)  \right)}{u \prod_{l=1}^d \left(1+(\tilde{\lambda}_l^{j-1}-\epsilon)^2 u^2\right)^{\frac{1}{4}}}du\\
    & \qquad - \frac{1}{\pi}\int_0^{\infty} \dfrac{	\sin \left(\frac{1}{2}  \sum_{l=1}^d \tan^{-1} \left((\tilde{\lambda}_l^{j-1} + \epsilon )u \right)  \right)}{u \prod_{l=1}^d \left(1+(\tilde{\lambda}_l^{j-1} + \epsilon)^2 u^2\right)^{\frac{1}{4}}}du.
\end{align*}
\end{theorem} 

\begin{proof}
The proof of this theorem is an adaptation of \cite[Theorem 9]{boman1999infeasibility}. Note that $\mathbb{P}[|X| \leq \eta] = \mathbb{P}[-\eta \leq X \leq \eta] = \mathbb{P}[X \leq \eta] - \mathbb{P}[X \leq -\eta]$, and $\mathbb{P}[|X| > \eta] = 1 - \mathbb{P}[|X| \leq \eta]$.
\end{proof}

As in the case for S-LBFGS, we now illustrate the probability of accepting pairs empirically. We conducted the same set of experiments as in Section \ref{sec.prob_SLBFGS}; see Table \ref{tbl:prob_details} for details. The probability of accpeting pairs for the synthetic problems is larger for S-LSR1 than S-LBFGS. This is due to the fact that negative values of the Rayleigh quotient are also accepted, as long as they are large enough in magnitude. That being said, the relative performance when the eigenvalues are chosen to be smaller is similar to the S-LBFGS method. For the toy classification problems, as for the S-LBFGS method, the probability of accepting curvature pairs is close to 100\% as long as $\epsilon$ is chosen to be small.


\begin{figure}[ht]
\centering
\subfloat[\label{fig.prob_1_sr1}]{%
\resizebox*{0.24\textwidth}{!}{\includegraphics[trim=0.6cm 7cm 3cm 7cm,width=\textwidth]{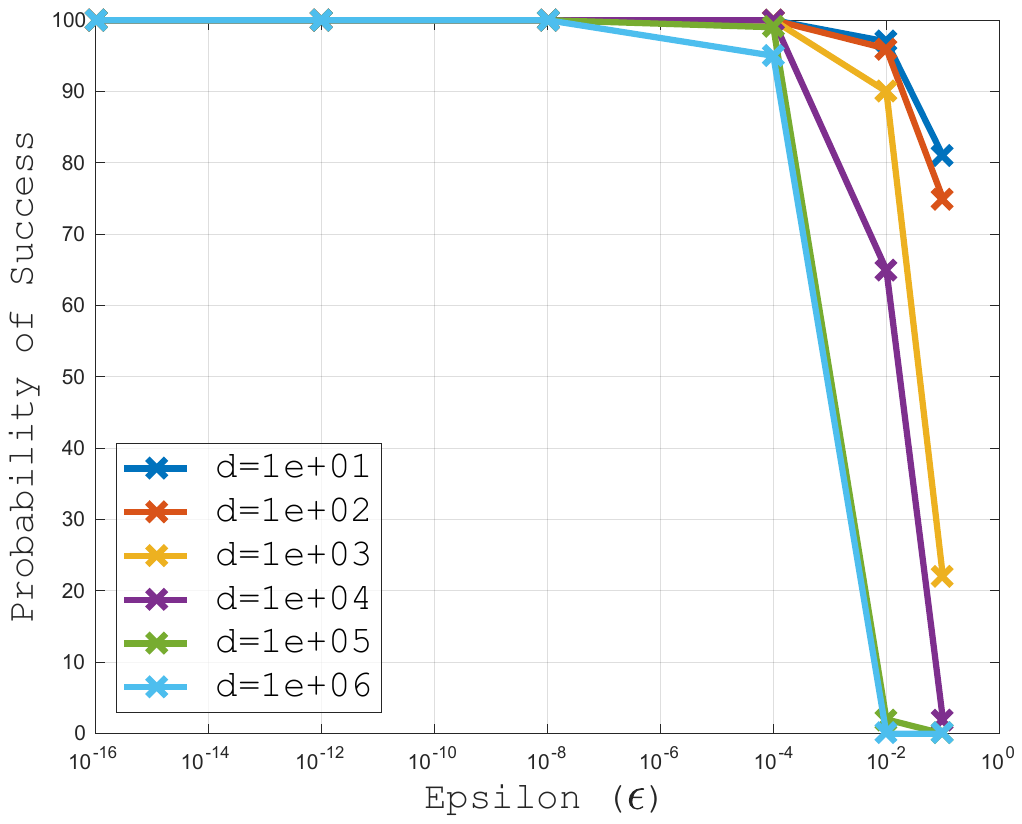}}}
\hspace{2pt}
\subfloat[\label{fig.prob_2_sr1}]{%
\resizebox*{0.24\textwidth}{!}{\includegraphics[trim=0.6cm 7cm 3cm 7cm,width=\textwidth]{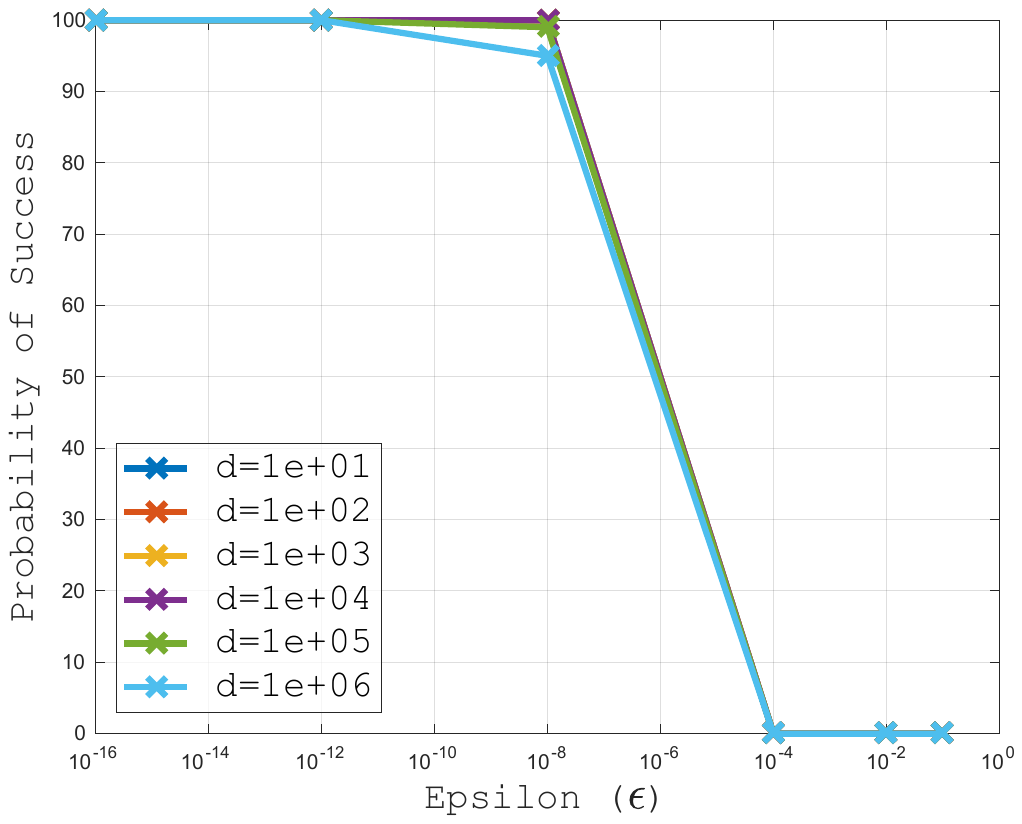}}}
\hspace{2pt}
\subfloat[\label{fig.prob_3_sr1}]{%
\resizebox*{0.24\textwidth}{!}{\includegraphics[trim=0.6cm 7cm 3cm 7cm,width=\textwidth]{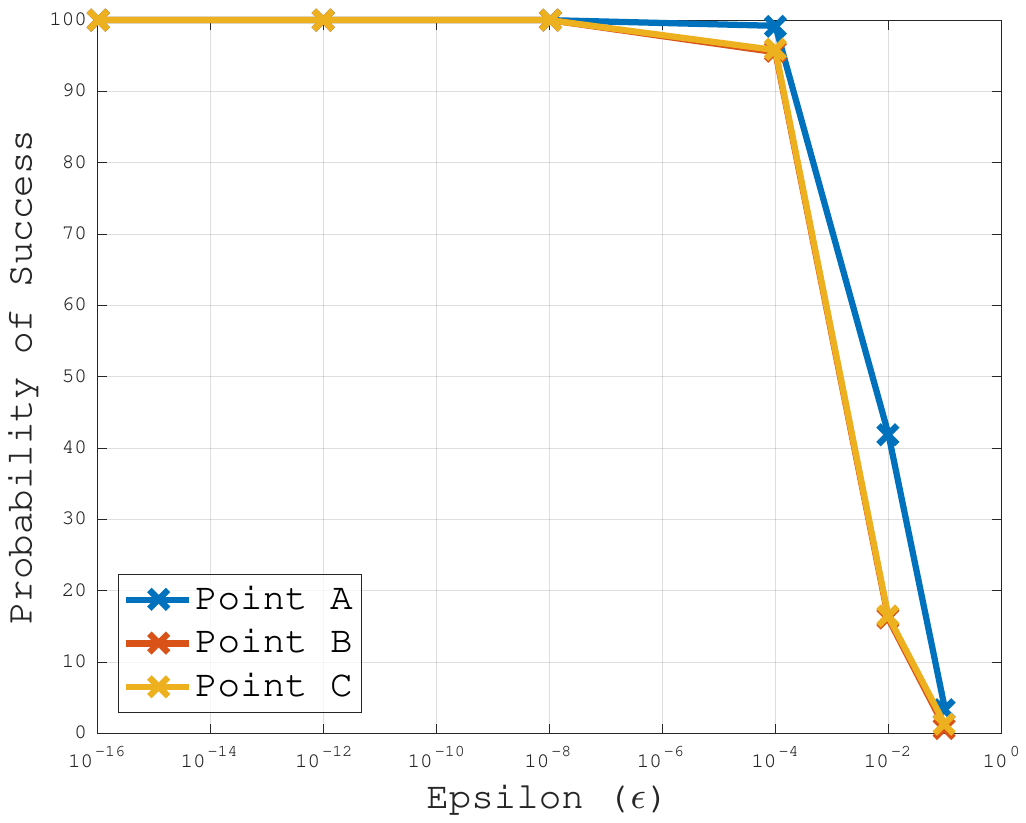}}}
\hspace{2pt}
\subfloat[\label{fig.prob_4_sr1}]{%
\resizebox*{0.24\textwidth}{!}{\includegraphics[trim=0.6cm 7cm 3cm 7cm,width=\textwidth]{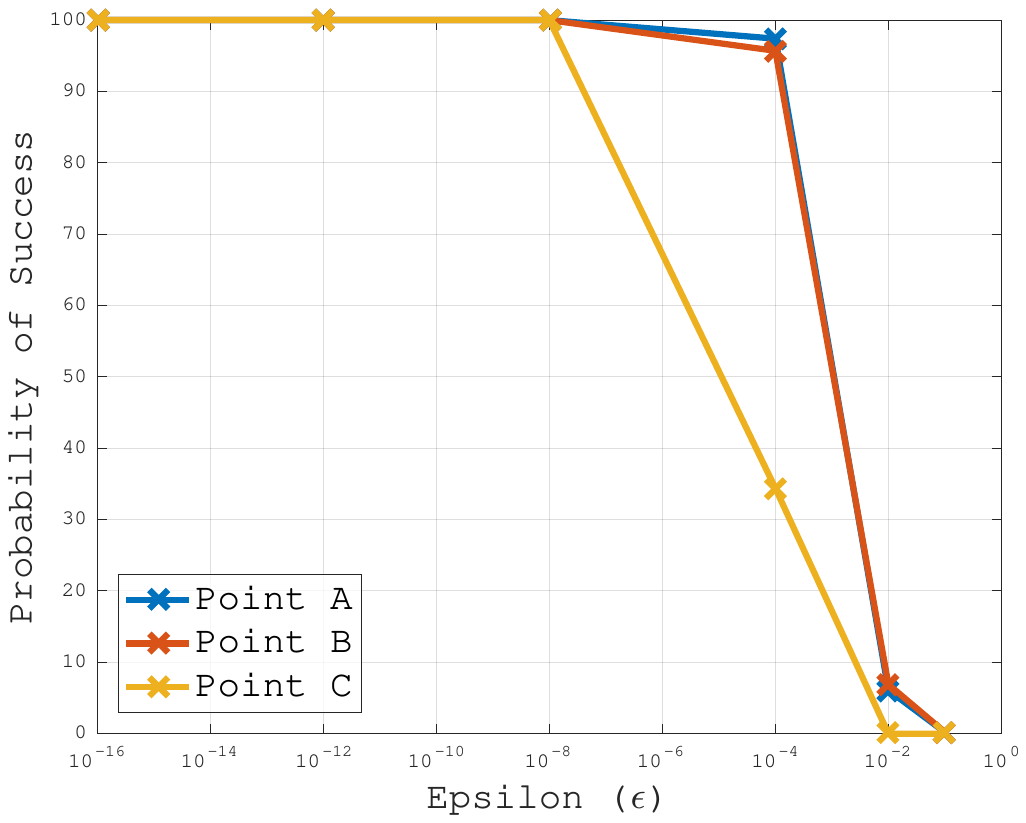}}}
\caption{\small Empirical investigation of accepting curvature pairs for different problems (S-LSR1).\label{fig.prob_results_slsr1}}
\end{figure}

\section{Numerical Experiments}
\label{sec:num_res}

In this section, we present numerical experiments on a toy classification problem as well as on popular benchmarking binary classification and neural network training tasks in order to illustrate the performance of our proposed sampled quasi-Newton methods.

\subsection{Method Specifications and Details}
Before we present the numerical results\footnote{All codes to reproduce results presented in this section are available at: \texttt{https://github.com/OptMLGroup/SQN}.}, we discuss the implementation details for all the methods. For ADAM \cite{kingma2014adam}, we tuned the steplength and batch size for each problem independently. For GD and BFGS-type methods, we computed the steplength using a backtracking Armijo line search \cite{nocedal_book}. For SR1-type methods, we solved the subproblems \eqref{eq:tr_obj} using CG-Steihaug \cite{nocedal_book}. For BFGS and SR1, we constructed the full (inverse) Hessian approximations explicitly, whereas for the limited-memory we never constructed the full matrices. For limited-memory BFGS methods we used the two-loop recursion to get the search direction \cite{nocedal_book}. Implementing the limited memory SR1 methods is not trivial; we made use of the compact representations of the SR1 matrices \cite{byrd1994} and computed the steps dynamically; see Appendix \ref{sec:impl_lsr1} for details.

\subsection{Toy Classification Problem}

\begin{wrapfigure}{r}{0.34\textwidth}

    \includegraphics[width=0.35\textwidth,
    trim={1cm 0.5cm 1cm 1cm},clip]{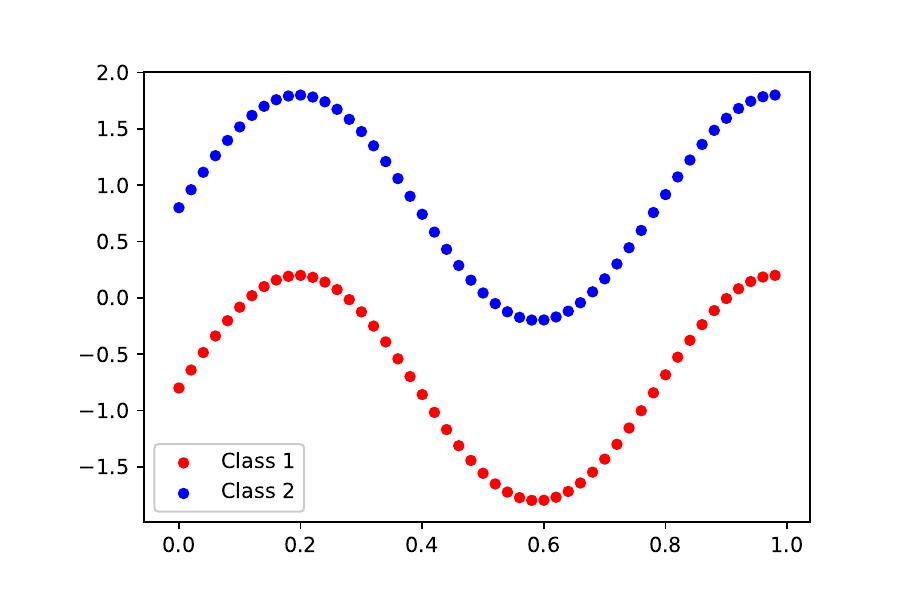}
    \vspace{-0.9cm}
  \caption{\small Toy Classification Problem}
  \label{toy_problem}
  \vspace{-0.4cm}
\end{wrapfigure}
Consider the following simple classification problem, illustrated in Figure \ref{toy_problem}, consisting of two classes (red and blue) each with 50 data points. The goal of this classification task is to find a nonlinear decision boundary that separates the two classes. We trained three fully connected neural networks--\texttt{small}, \texttt{medium} and \texttt{large}--with sigmoid activation functions and 4 hidden layers. The details of the three networks are summarized in Table~\ref{tbl:toy_details}.

\begin{table}[h]
\caption{\small Toy Classification Problem: Neural Network Details}
\label{tbl:toy_details}
\centering
\begin{small}
\begin{tabular}{@{}lccr@{}}\toprule
\textbf{network} & \textbf{structure} & \textbf{$d$}  \\ \midrule

\texttt{small} & 2-2-2-2-2-2 & $36$ 
\\ \hdashline  

\texttt{medium} & 2-4-8-8-4-2 & $176$ 
\\ \hdashline  

\texttt{large} & 2-10-20-20-10-2 & $908$

\\
\bottomrule 
\end{tabular}
\end{small}
 \end{table}
 
\begin{figure}[h]
\centering

\includegraphics[width=1\textwidth]{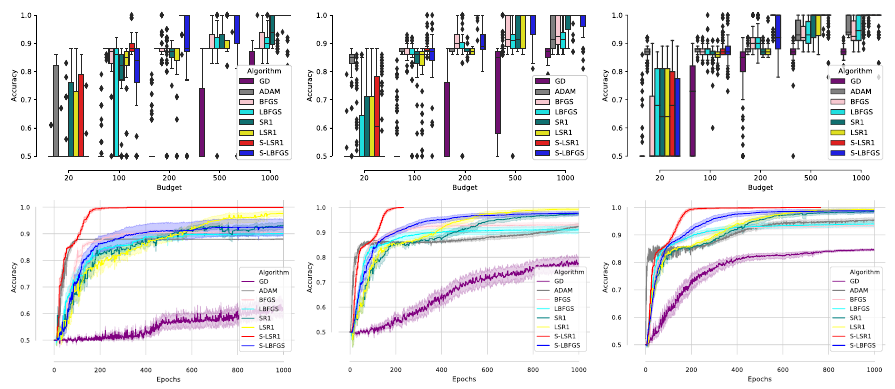}

%
%
    
    \caption{\small Performance of GD, ADAM, BFGS, LBFGS, SR1, LSR1, S-LSR1 and S-LBFGS on toy classification problems. Networks: \texttt{small} (left); \texttt{medium} (center); \texttt{large} (right).}
	\label{toy_boxplot}
\end{figure}

\begin{figure}[]
    \centering
    
    \includegraphics[width=1\textwidth]{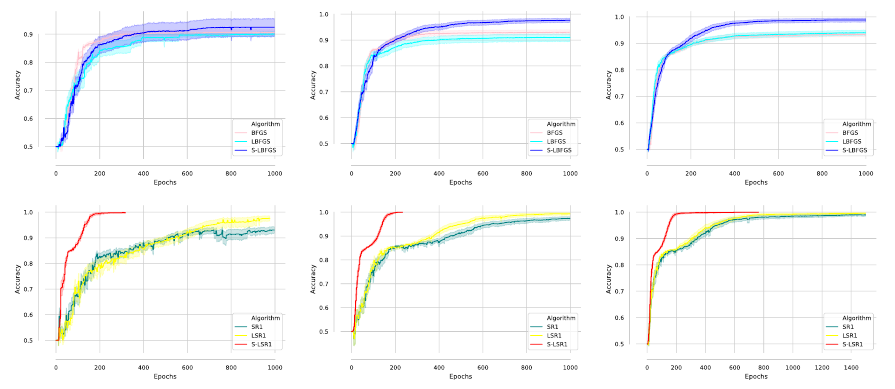}

    \caption{Performance of: BFGS, LBFGS, S-LBFGS (top); SR1, LSR1 and S-LSR1 (bottom), on toy classification problems. Networks: \texttt{small} (left); \texttt{medium} (center); \texttt{large} (right).}
	\label{toy_acc_fam}
\end{figure}

 For this problem, we ran each method 100 times starting from different initial points and show the results for different budget levels. The results are summarized in Figure \ref{toy_boxplot}. In order to better visualize the relative performance of our proposed sampled quasi-Newton methods compared to the classical variants, we show accuracy vs. epochs  plots in Figure \ref{toy_acc_fam}. As is clear from the figures, the proposed methods outperform their classical variants as well as the first-order methods. See Appendix \ref{sec:app_toy} for more results. 

For this problem, we ran each method 100 times starting from different initial points and show the results for different budget levels. The results are summarized in Figure \ref{toy_boxplot}. In order to better visualize the relative performance of our proposed sampled quasi-Newton methods compared to the classical variants, we show accuracy vs. epochs  plots in Figure \ref{toy_acc_fam}. As is clear from the figures, the proposed methods outperform their classical variants as well as the first-order methods. See Appendix \ref{sec:app_toy}  for more results.

The toy classification problem is inherently complex. As is clear from the results, first-order methods (GD and ADAM) are not competitive with other reported methods, as they require a significant computational budget in order to achieve low classification error. It is worth noting that as we increase the size of the neural networks (over-parameterization), the performance of these methods becomes better. On the other hand, quasi-Newton methods have better performance in this complex, albeit small, problem, primarily due to the use of curvature information. Amongst the reported results, our sampled quasi-Newton methods significantly outperform the classical methods. We posit that this is the case due to the use of more recent and local curvature information in the updates.

\subsection{Logistic Regression}\label{sec.log_reg}

\begin{figure}[]
	\centering
	
	\includegraphics[width=0.8\textwidth]{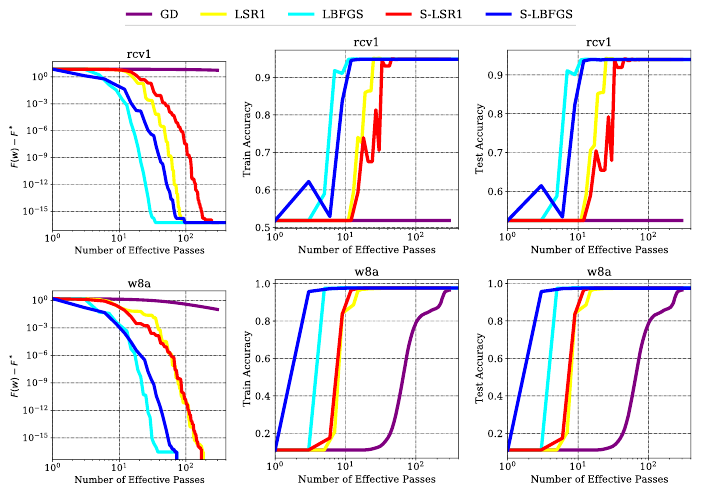}

	\caption{Performance of GD, LBFGS, LSR1, S-LSR1 and S-LBFGS on Logistic Regression problems; \texttt{rcv1} dataset (first row) and \texttt{w8a}  dataset (second row).}

	\label{log_regression}

\end{figure}

Next we consider $\ell_2$-regularized logistic regression problems of the form
\begin{align*}
    F(w) = \frac{1}{n}\sum_{i=1}^n \log  \left(1+e^{-y_ix_i^Tw} \right) + \frac{\lambda}{2}\|w\|^2,
\end{align*}
where $(x_i,y_i)_{i=1}^n$ denote the training examples and $\lambda>0$ is the regularization parameter. We present results on two popular machine learning datasets (\texttt{rcv1} and \texttt{w8a}; \cite{chang2011libsvm}); see Appendix \ref{sec:logistic}  for dataset details and more results. We compared the performance of the proposed sampled quasi-Newton methods with gradient descent (GD) and classical quasi-Newton methods (LSR1 and LBFGS). Figures \ref{log_regression} illustrates the performance of the methods in terms of optimality gap (training loss), training accuracy and testing accuracy. As is clear from Figure \ref{log_regression}, the sample quasi-Newton methods are competitive with the classical variants in terms of all three metrics. One can also observed that in the initial stages of the optimization, it appears that the sampled quasi-Newton methods outperform their classical counterparts.

%
%

\subsection{Nonlinear Least Squares}
\begin{figure}[]
	\centering
	
	\includegraphics[width=0.8\textwidth]{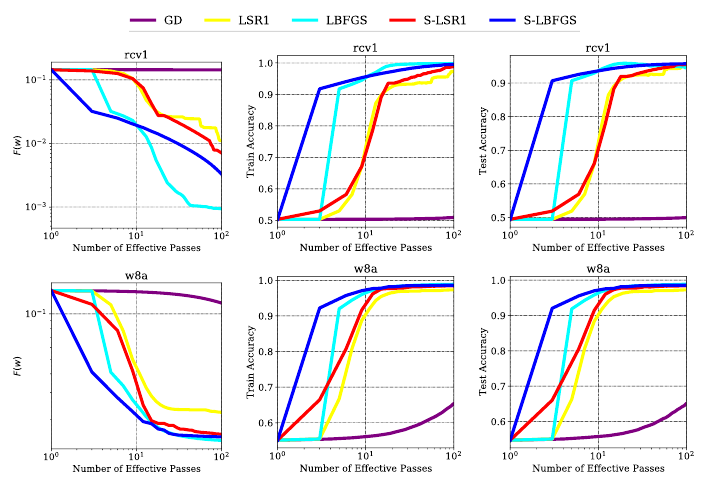}

	\caption{Performance of GD, LBFGS, LSR1, S-LSR1 and S-LBFGS on Nonlinear Least Squares problems; \texttt{rcv1} dataset (first row) and \texttt{w8a}  dataset (second row).}

	\label{nonlinear_least_sqs}
\end{figure}

In this section, we consider nonlinear least squares problems \cite{xu2017second} of the form
\begin{align*}
    F(w) = \frac{1}{n}\sum_{i=1}^n   \left(y_i - \frac{1}{1+e^{-x_i^Tw}}  \right)^2,
\end{align*}
where $(x_i,y_i)_{i=1}^n$ denote the training examples. We present results on the same datasets and compare against the same methods and using the same metrics as in Section \ref{sec.log_reg}. As is clear from Figure \ref{nonlinear_least_sqs}, the sampled quasi-Newton methods outperform their classical counterparts across the board. This is consistent with the results on other datasets; see Appendix \ref{sec:nonLS}.  We posit that this is due to the fact that sampling curvature pairs at every iterations allows for the method to incorporate more recent, local and reliable curvature information, a feature that can be indispensable in the nonconvex setting.

\subsection{Neural Network Training: MNIST and CIFAR10}

\begin{figure}[H]
	\centering
	
	\includegraphics[width=1\textwidth]{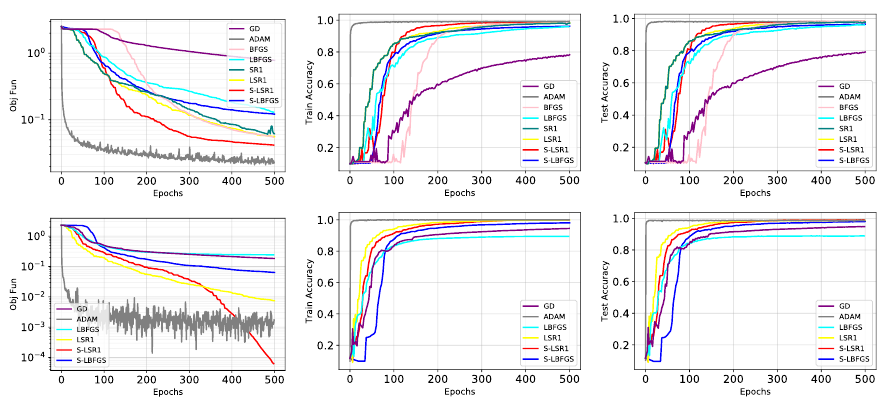}

	\caption{Performance of GD, ADAM, BFGS, LBFGS, SR1, LSR1, S-LSR1 and S-LBFGS on \texttt{MNIST} problems on \texttt{Net1} (first row) and \texttt{Net2} (second row).}

	\label{accuarcyMNISTNew}

\end{figure}

\begin{figure}[H]
	\centering

	\includegraphics[width=1\textwidth]{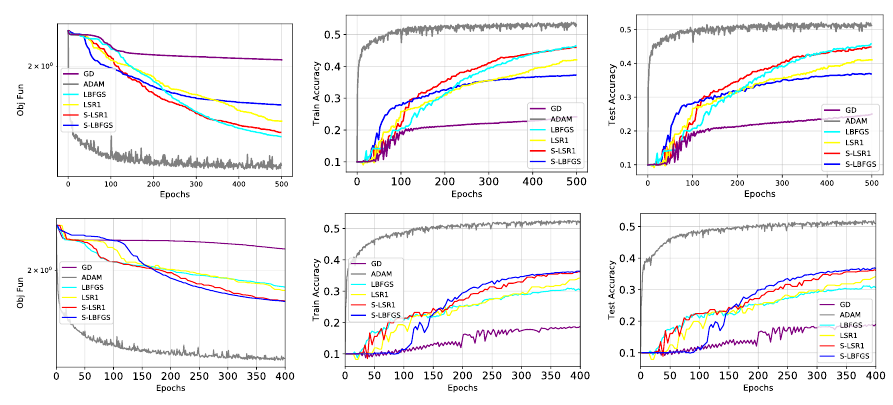}

	\caption{Performance of GD, ADAM, BFGS, LBFGS, SR1, LSR1, S-LSR1 and S-LBFGS on \texttt{CIFAR10} problems on \texttt{Net3} (first row) and \texttt{Net4} (second row).}

	\label{accuarcyCIFAR10New}

\end{figure}

We illustrate the performance of the sampled quasi-Newton methods on standard benchmarking neural network training tasks: \texttt{MNIST} \cite{lecun1998gradient} and \texttt{CIFAR10} \cite{krizhevsky2009learning}. The details of the problems are given in Table \ref{tbl:MNIST_CIFAR10_new}. For these problems we used sigmoid and softplus activation functions and softmax cross-entropy loss. For the memory-variant algorithms, we considered the memory from the set $m \in \{4,16,64,256\}$, and report the best performance with respect to the different memory sizes. The results of these experiments are given in Figures \ref{accuarcyMNISTNew} and \ref{accuarcyCIFAR10New}.

\begin{table}[H]
\caption{\small Details for \texttt{MNIST} and \texttt{CIFAR10} Problems.}
\label{tbl:MNIST_CIFAR10_new}
\centering
\begin{small}
\begin{tabular}{@{}lcccr@{}}\toprule
\textbf{problem} &\textbf{network} & \textbf{structure} & \textbf{$d$}  \\ \midrule
\texttt{MNIST} &\texttt{Net1} & $784-C_{5,3}-C_{5,5}-10-10$ & $1378$ 
\\ \hdashline  
\texttt{MNIST} & \texttt{Net2} & $784-C_{5,6}-C_{5, 16}-C_{4,120}-82-10$ & $44164$ 
\\ \midrule
\texttt{CIFAR10} &\texttt{Net3} & $1024,3-C_{5,3}-C_{5,4}-10-10$ & $1652$ 
\\ \hdashline  
\texttt{CIFAR10} & \texttt{Net4} & $1024,3-C_{5,3}-C_{5,5}-16-32-10$ & $2312$ 
\\ \midrule
\end{tabular}

$C_{k,ch}$: convolution with kernel $k$
and $ch$ output channels.\\
Net2 is equivalent to LeNet structure with 1-channel input.
\end{small}

\end{table}

Overall, the sampled quasi-Newton methods outperform their classical variants. We should note that the goal of these experiments is not to perform better than ADAM, rather the performance of ADAM can be viewed as a benchmark. The reasons for this are two-fold. First, ADAM is a stochastic algorithm while the other reported methods are deterministic. Second, we report results for the best hyper-parameter settings for ADAM (\textit{well-tuned}); see Appendix \ref{sec:app_MNIST}, while the other methods do not require tuning or they are insensitive to the choice of hyper-parameters. 

For the \texttt{MNIST} problems, the S-LSR1 method is able to achieve comparable accuracy to that of \textit{well-tuned} ADAM, after a lot more epochs. That being said, in a distributed setting, the time to perform one iteration (one epoch) of S-LSR1 is significantly smaller than the time to perform one epoch of ADAM, and as such in terms of Wall Clock Time, the proposed method could be more efficient. With regards to the \texttt{CIFAR10} problems, one can observe that our proposed sampled methods perform on par if not better than classical quasi-Newton methods. We posit that the reason that S-LSR1 has better performance than S-LBFGS is due to the possible utilization of negative curvature in the updates.

\section{Final Remarks and Future Work}
\label{sec:fin_rem}

This paper describes two novel quasi-Newton methods; S-LBFGS and S-LSR1. Contrary to classical quasi-Newton methods, these methods \textit{forget} past curvature information and \textit{sample} new curvature information at every iteration. Numerical results show that the methods are efficient in practice, and the convergence guarantees of the methods match those of the classical variants.

Our algorithms can be extended to the stochastic setting where gradients and/or Hessians are computed inexactly. Moreover, the algorithms could be made adaptive following the ideas from \cite{NIPS2016_6262,
jahani2018efficient}. Furthermore, stronger theoretical (e.g., superlinear convergence) results could be proven for some variants of the sampled quasi-Newton methods. Finally, a large-scale numerical investigation of the method would test the limits of these methods.

\section*{Acknowledgements}

This work was partially
supported by the U.S. National Science Foundation, under award numbers NSF:CCF:1618717 and  NSF:CCF:1740796, DARPA Lagrange award HR-001117S0039,
and XSEDE Startup grant IRI180020.

\appendix

\section{Additional Numerical Experiments and Method Details}

In this section, we present additional numerical results and expand on some details about the methods.

\subsection{Motivation Figure}

In this section, we present more motivating plots showing the accuracy vs. iterations and accuracy vs. epochs for a toy classification problem. In the following experiments, we ran each method from 10 different initial points.

\begin{figure}[h!]
	\centering
	\includegraphics[width=0.45\textwidth]{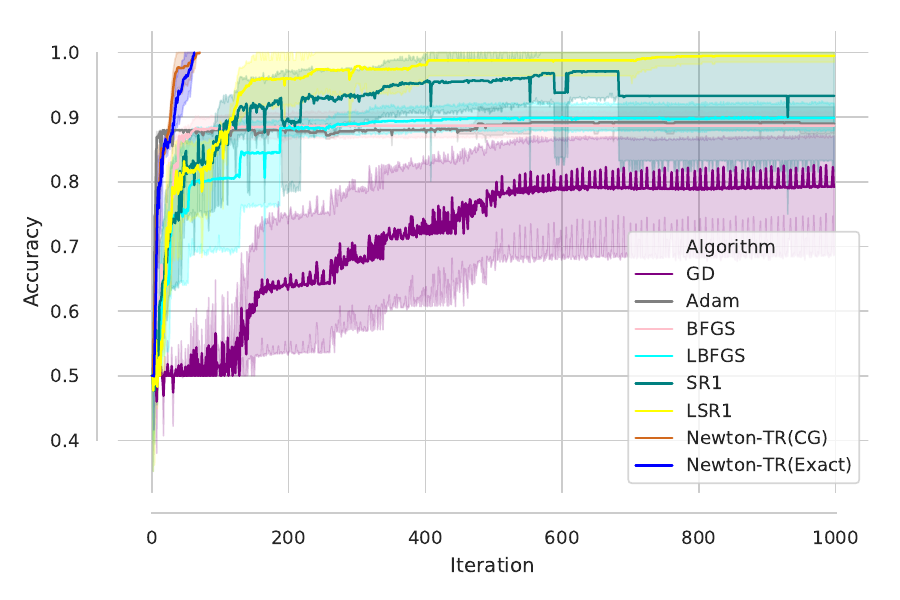}
	\includegraphics[width=0.45\textwidth]{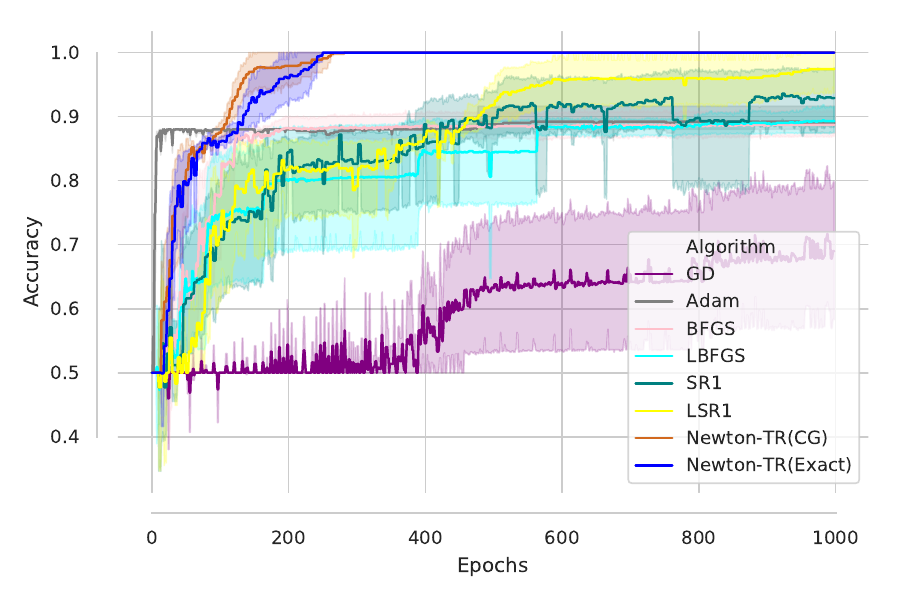}
	
	\includegraphics[width=0.45\textwidth]{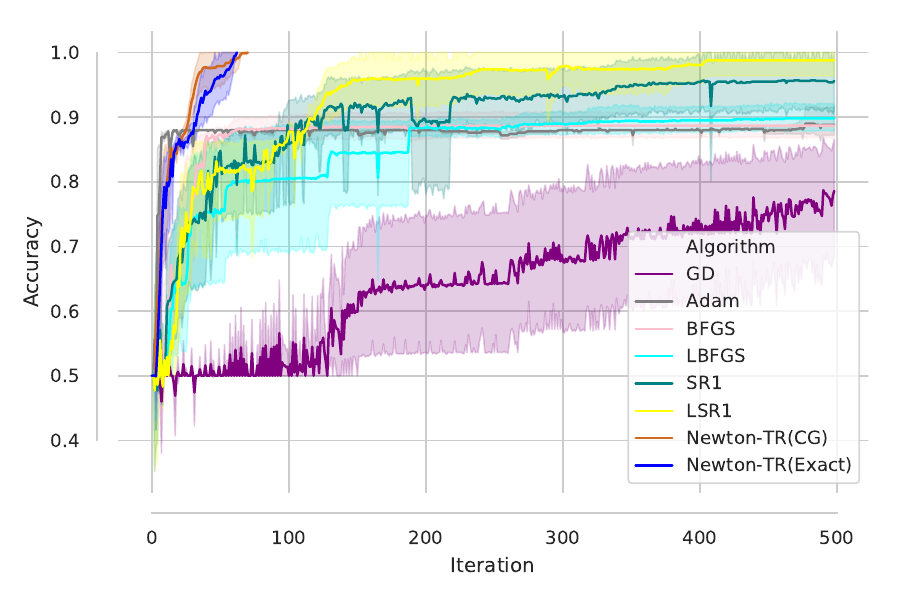}
	\includegraphics[width=0.45\textwidth]{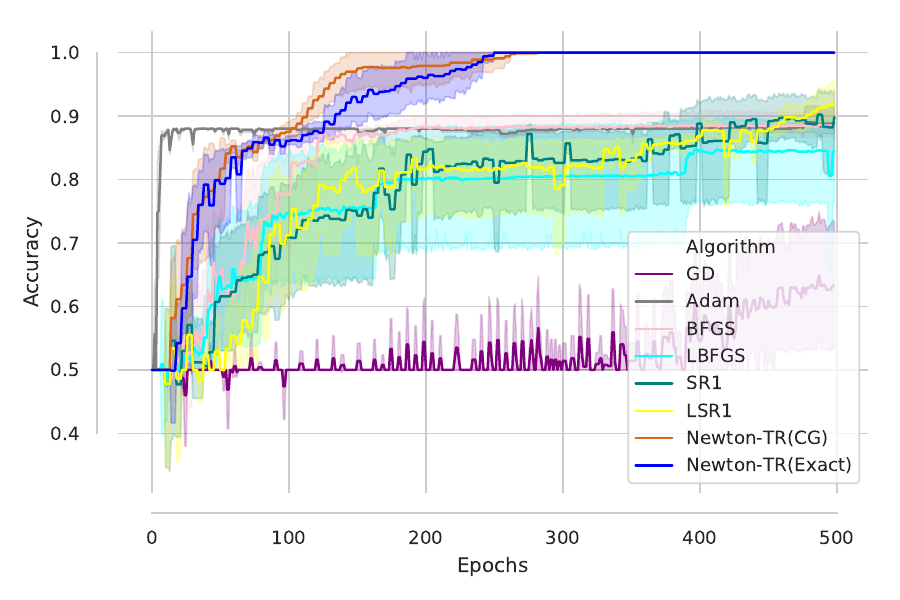}
	\caption{Comparison of Gradient Descent (GD), ADAM, BFGS, LBFGS, SR1, LSR1, Newton-TR (Exact), Newton-TR (CG) on a toy classification problem in terms of iterations and epochs.}
	\label{MotivationFigure_App}
\end{figure}

\clearpage

\subsection{Eigenvalue Figures}
\label{sec:eigs_app}

\begin{figure}[ht]
	\centering

	\includegraphics[width=0.45\textwidth]{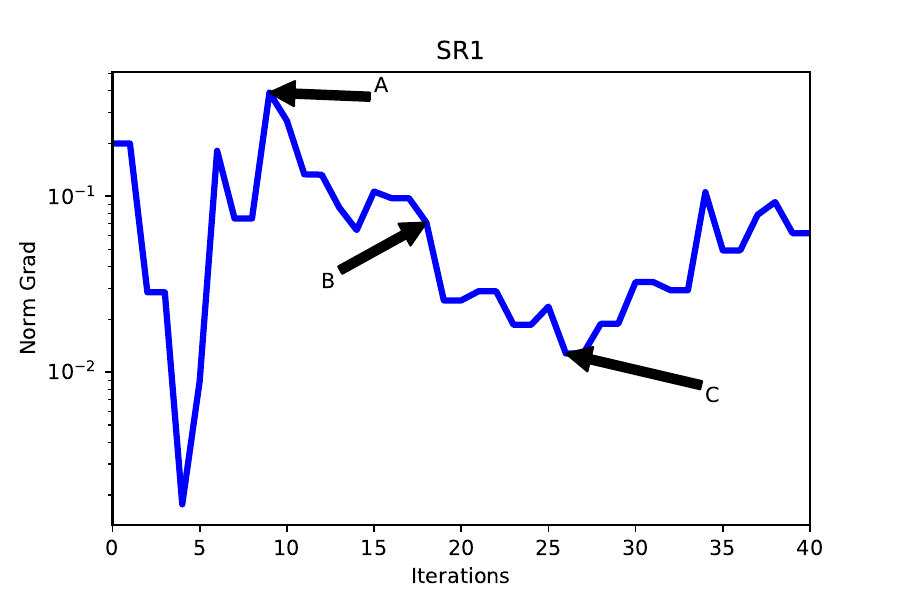}
	\includegraphics[width=0.45\textwidth]{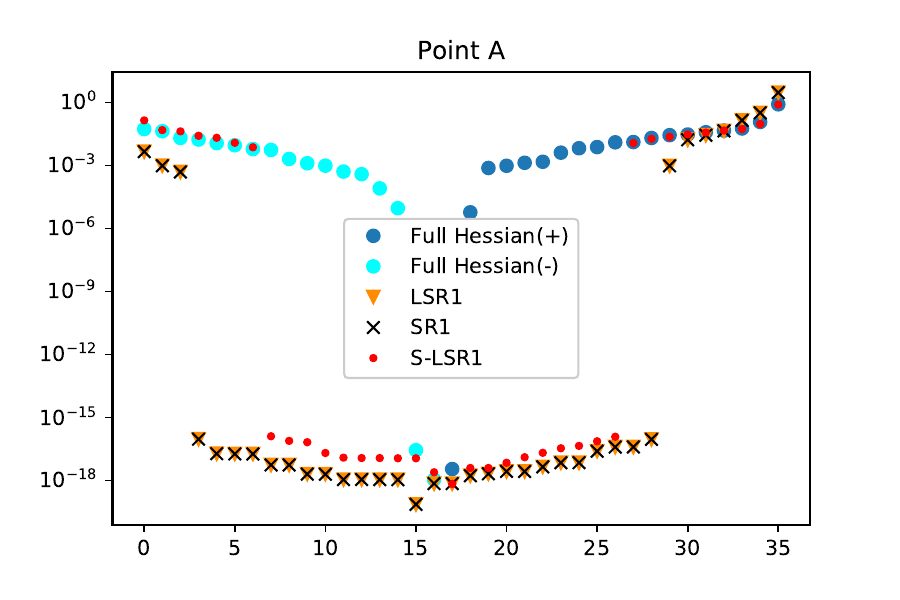}
	
	\includegraphics[width=0.45\textwidth]{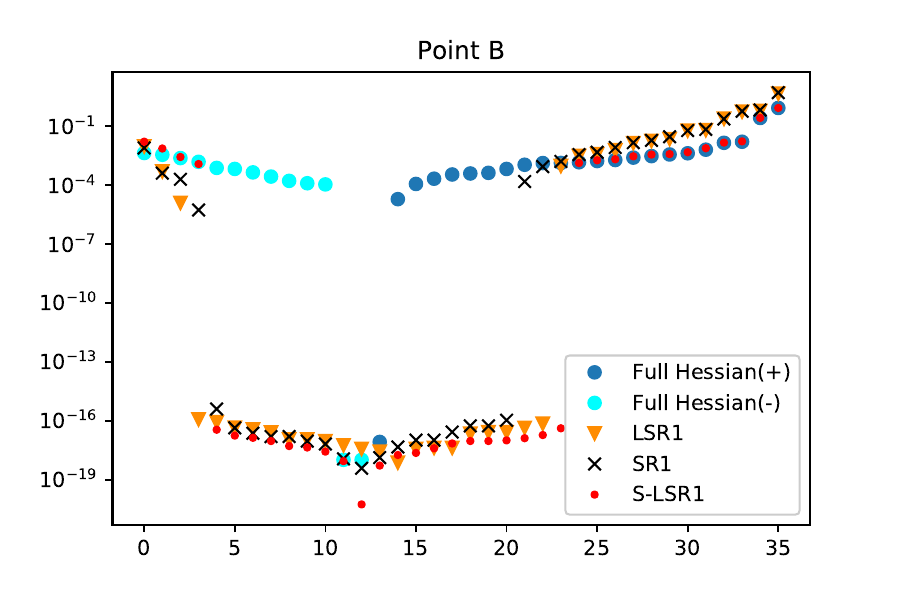}
	\includegraphics[width=0.45\textwidth]{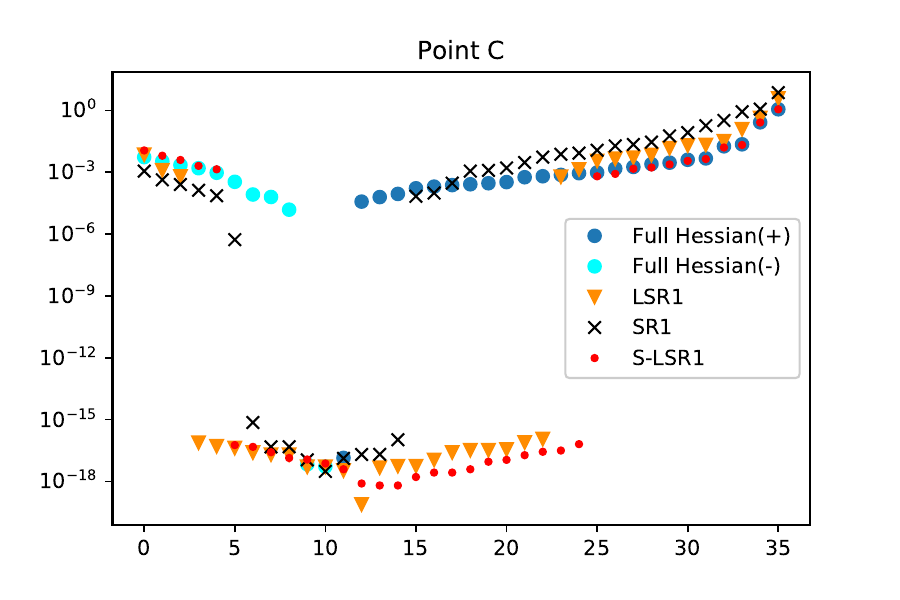}

	\caption{Comparison of the eigenvalues of SR1, LSR1 and S-LSR1 at different points for a toy classification problem. \label{fig:eig36}}
\end{figure}

\begin{figure*}[ht]
	\centering

	\includegraphics[width=0.45\textwidth]{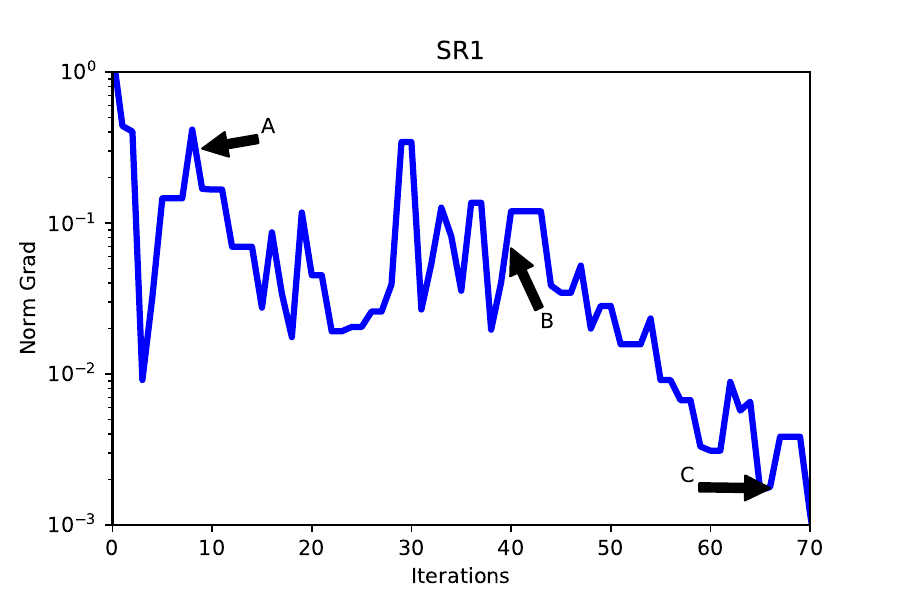}
	\includegraphics[width=0.45\textwidth]{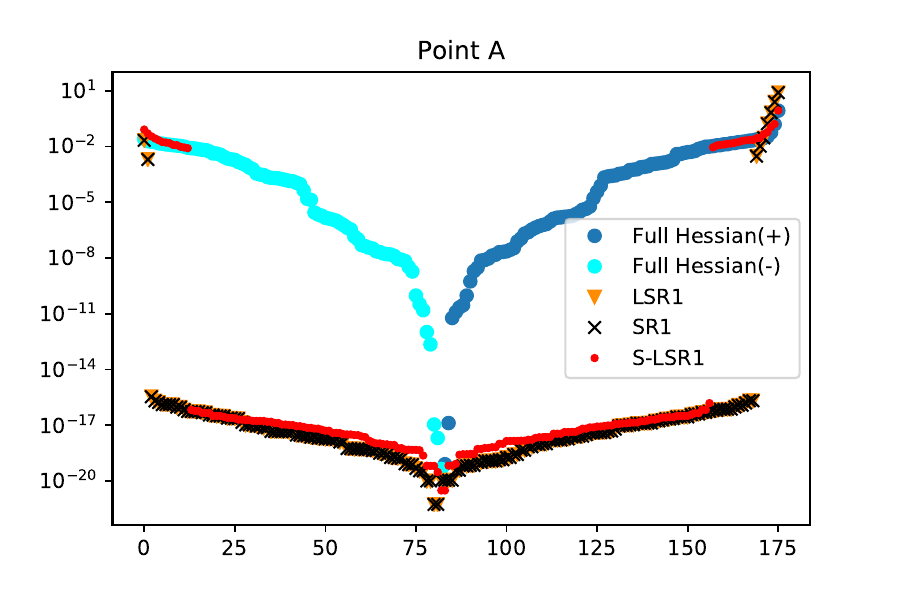}
	
	\includegraphics[width=0.45\textwidth]{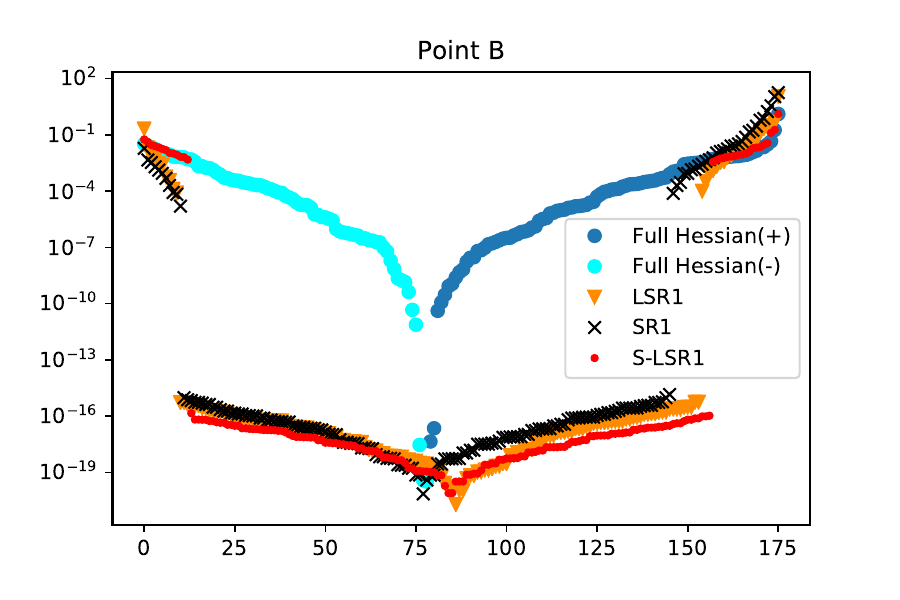}
	\includegraphics[width=0.45\textwidth]{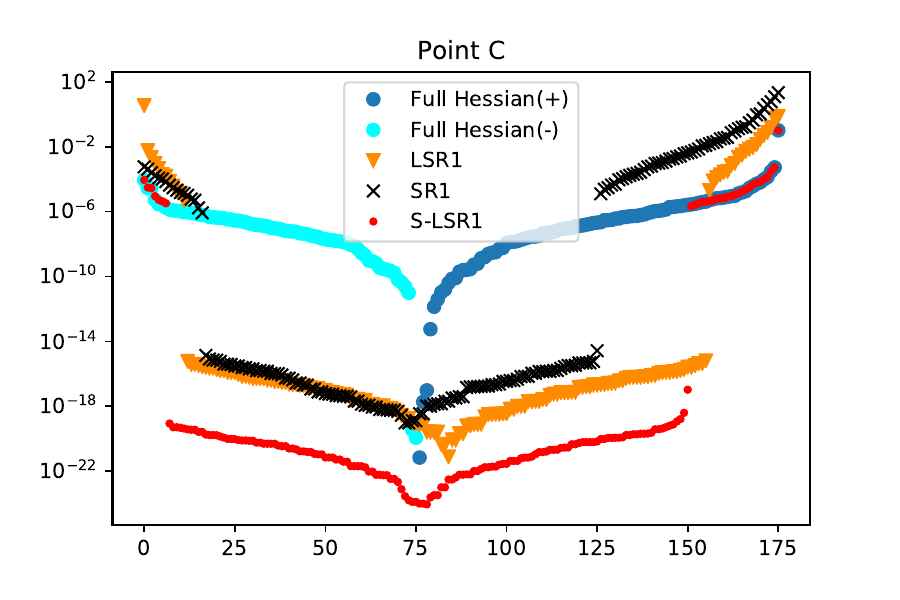}

	\caption{Comparison of the eigenvalues of SR1, LSR1 and S-LSR1 at different points for a toy classification problem. \label{fig:eig176}}
\end{figure*}

In this section, we describe the procedure in which Figures \ref{fig:eig36} and \ref{fig:eig176} (Figures 2 \& 3) were constructed. We plot the same figures below for ease of exposition.


To calculate the eigenvalues for SR1, LSR1 and S-LSR1 we used the following procedure.
\begin{enumerate}
    \item We ran the SR1 method for T iterations on a toy classification problem. During the optimization, we computed the eigenvalues of the SR1 Hessian approximation at several points (e.g., A, B and C); black $\times$ marks on plots.
    \item We stored all the curvature pairs $\{s_k,y_k \}_{k=1}^T$ and the iterates $\{ w_k\}_{k=1}^T$.
    \item We constructed the true Hessian at all iterations and computed the eigenvalues of the true Hessian; dark blue $\bullet$ (positive eigenvalues)  and light blue $\bullet$ (negative eigenvalues)  marks on plots. 
    \item We then computed the limited-memory SR1 Hessian approximations at several points (e.g., A, B and C) using the $m$ most recent pairs and computed the eigenvalues of the approximations; orange $\blacktriangledown$ marks on plots.
    \item Finally, we used the iterate information at points A, B and C, sampled $m$ points at random around those iterates with sampling radius $r$, constructed the sampled LSR1 Hessian approximations and computed the eigenvalues of the approximations; red $\bullet$ marks on plots. 
\end{enumerate}

Note: for Figure \ref{fig:eig36} we used $T = 40$, $m =16$ and $r = 0.01$, and for Figure \ref{fig:eig176} we used $T = 70$, $m =32$ and $r = 0.01$.

As is clear, the eigenvalues of the sampled LSR1 Hessian approximations better match the eigenvalues of the true Hessian. Similar results were obtained for other problems and for different parameters $m$ and $r$.


\subsection{Trust-Region Management Subroutine}
\label{sec:tr_alg}

In this section we present, in detail, the Trust-Region management subroutine ($\Delta_{k+1} = \texttt{adjustTR}(\Delta_{k},\rho_k)$) that is used in Algorithm 3. See \cite{nocedal_book} for further details.

\begin{center}
\begin{minipage}{.65\linewidth}
\begin{algorithm}[H]
\caption{$\Delta_{k+1} = \texttt{adjustTR}(\Delta_{k},\rho_k)$: Trust-Region management subroutine}
  \label{alg:tr_mgmt}
 {\bf Input:} $\Delta_k$ (current trust region radius), $0 \leq \eta_3 < \eta_2 < 1$, $\gamma_1 \in (0,1)$, $\zeta_1 > 1$, $\zeta_2 \in (0,1)$ (trust region parameters). 
  \begin{algorithmic}[1]
    \If {$\rho_k > \eta_2$}
        \If {$\|p_k\| \leq \gamma_1 \Delta_k$}
            \State Set $\Delta_{k+1} = \Delta_k$
        \Else
            \State Set $\Delta_{k+1} = \zeta_1 \Delta_k$
        \EndIf
    \ElsIf{$\eta_3 \leq \rho_k \leq \eta_2$}       \State Set         $\Delta_{k+1} =     \Delta_k$
    \Else
        \State $\Delta_{k+1} = \zeta_2 \Delta_k$
    \EndIf  
  \end{algorithmic}
\end{algorithm}
\end{minipage}
\end{center}

\clearpage

\subsection{Hessian-Free Implementation of Limited-Memory SR1 Methods}
\label{sec:impl_lsr1}
In this section, we discuss the practical implementation of limited-memory SR1 methods where we need not construct the Hessian approximation $B_k$ explicitly. For the purpose of this discuss we focus on the S-LSR1 method, but a similar approach can be used for the LSR1 method too. To do so, we utilize the compact representation of the Hessian approximation discussed in \cite{byrd1994} which is equivalent to \eqref{eq:sr1_hess}, i.e., 
\begin{align*}
    B_{k+1} = B_k +  \tfrac{(y_k - B_ks_k)(y_k - B_ks_k)^T}{(y_k - B_ks_k)^Ts_k}. 
\end{align*}
The compact representation can be expressed as follows:
\begin{align}   \label{eq:compactForm}
 B_{k+1} = B_k^{(0)} + (Y_k- B_k^{(0)}S_k)\Big(D_k+L_k +L_k^T -S_k^TB_k^{(0)}S_k\Big)^{-1}(Y_k-B_k^{(0)}S_k)^T,
\end{align}

where $S_k=[s_{k,1},s_{k,2},\dots, s_{k,m}] \in \mathbb{R}^{d\times m}$ and $Y_k=[y_{k,1},y_{k,2},\dots, y_{k,m}] \in \mathbb{R}^{d\times m}$, and $B_k^{(0)}$ is a symmetric positive definite initial Hessian approximation, which for the purpose of this discussion we assume has the form $B_k^{(0)} = \gamma_k I$ ($0 \leq \gamma_k < \gamma < \infty$). In \eqref{eq:compactForm}, $D_k$ and $L_k$ are two $m\times m$ matrices that are defined as follows, 
\begin{align}   \label{eq:hessFreeD}
 D_k = {diag}[s_{k,1}^Ty_{k,1}, \dots ,s_{k,m}^T y_{k,m}]
\end{align}
\begin{align}   \label{eq:hessFreeL}
 (L_k)_{i,j} = \begin{cases}
  s_{k,i-1}^T y_{k,j-1} & \text{if i  $>$ j} \\
  0 & \text{otherwise}
\end{cases}
\end{align}
The curvature pairs in the matrices $S_k$ and $Y_k$ are pairs that satisfy the condition given in \eqref{eq:cond_slsr1}, i.e., 
\begin{align}   \label{eq:cond_slsr1}
    | s^T(y-Bs) |  > \epsilon \| s\|^2,
\end{align}
where $\epsilon > 0$ is a predetermined constant.

In large-scale applications, it is not memory-efficient, or even possible for some applications, to store a $d\times d$ Hessian approximation matrix $B_{k+1}$. Instead, we can calculate the Hessian vector product $B_{k+1}v$, for some $v \in \mathbb{R}^d$, by leveraging the compact form of $B_k$ in \eqref{eq:compactForm} as follows:
\begin{align}   \label{eq:compactFormVec}
     B_{k+1}v = B_k^{(0)} v+ (Y_k- B_k^{(0)}S_k)\Big(D_k+L_k +L_k^T -S_k^TB_k^{(0)}S_k\Big)^{-1}(Y_k-B_k^{(0)}S_k)^Tv
\end{align}
The above, $B_{k+1}v$, is very efficient in terms of memory, and even more importantly efficient to compute;  the complexity of computing $B_{k+1}v$ is $\mathcal{O}(m^2d)$. 

In  Algorithm 3, we need to compute and use $B_{k+1}v$ in the following parts: (1) in checking the condition \eqref{eq:cond_slsr1}; (2) as part of the computation of solving the subproblem \eqref{eq:tr_obj}, i.e., 
\begin{align*}
    &{\min_{p}}  \; m_k(p) = F(w_k) + \nabla F(w_k)^Tp + \tfrac{1}{2} p^T B_k p ,\\
    & \quad \text{s.t.} \qquad \| p \| \leq \Delta_k, \nonumber
\end{align*}
using the CG solver (see \cite{nocedal_book}); and, (3) in the calculation of $\rho_k$.

In the remainder of this section we describe the steps for checking whether the curvature pairs constructed by  Algorithm 1 satisfy \eqref{eq:cond_slsr1}. This is by no means a trivial task; several researchers have proposed mechanisms for doing this \cite{brust2017solving,lu1996study} by using spectral decomositions of $B_k$. We propose to do this in a dynamic manner leveraging \eqref{eq:compactFormVec}. 

The condition that we want to check \eqref{eq:cond_slsr1} has the following form:
\begin{align}   \label{eq:cond_slsr1tmp}
    | s_{k,i}^T(y_{k,i}-B_k^{(i-1)}s_{k,i}) | \geq \epsilon \| s_{k,i}\| \| y_{k,i}-B_k^{(i-1)}s_{k,i} \|,
\end{align}
for $i = 1,\dots,m$, where $B_k^{(i)}$ are constructed recursively via,
\begin{align}   \label{eq:compactFormVecS1}
    B_k^{(i)} = B_k^{(0)} + (Y_k^i- B_k^{(0)}S_k^i)\Big(D_k^i+L_k^i +(L_k^i)^T -(S_k^i)^TB_k^{(0)}S_k^i\Big)^{-1}(Y_k^i-B_k^{(0)}S_k^i)^T,
\end{align}
$S_k^i = [s_{k,1},s_{k,2},\dots, s_{k,i}]$, $Y_k^i = [y_{k,1},y_{k,2},\dots, y_{k,i}]$, and $D_k^i$ and $L_k^i$ are defined in equations \eqref{eq:hessFreeD} and \eqref{eq:hessFreeL}, respectively, using $S_k^i$ and $Y_k^i$.
Of course we want to check condition \eqref{eq:cond_slsr1tmp} without explicitly forming the matrices $B_k^{(i)}$, and instead construct $B_k^{(i)}s_{k,i+1}$ directly. To this end, by using \eqref{eq:compactFormVec}, for any $i=0,\dots, m-1$ we have:
\begin{align}   \label{eq:compactFormVecS}
 B_k^{(i)} s_{k,i+1} = B_k^{(0)} s_{k,i+1} + (Y_k^i- B_k^{(0)}S_k^i)\Big(D_k^i+L_k^i +(L_k^i)^T -(S_k^i)^TB_k^{(0)}S_k^i\Big)^{-1}(Y_k^i-B_k^{(0)}S_k^i)^Ts_{k,i+1}
\end{align}
where the matrices $S_k^i$, $Y_k^i$, $D_k^i$ and $L_k^i$ are defined as above. 

The steps for checking \eqref{eq:cond_slsr1} are as follows:
\begin{enumerate}
    \item Consider the $k$th iteration of  Algorithm 3, where the pairs $\bar{S}_k=[s_{k,1},\dots,s_{k,m}]$ and $\bar{Y}_k=[y_{k,1},\dots,y_{k,m}]$ and $B_k^0$ are constructed by  Algorithm 1. Let $S_k = [\;]$ and $Y_k= [\;]$ be two empty matrices. 
    \item For any $i=1,\dots,m$, consider the pair $(s_{k,i},y_{k,i})$  and compute $B_{k}^{(i-1)}s_{k,i}$ by using \eqref{eq:compactFormVecS} and the updated lists $S_k$ and $Y_k$. Note that for $i=1$ the matrices $S_k$ and $Y_k$ are empty and $B_k^{(0)}$ is the initial Hessian approximation, thus condition \eqref{eq:cond_slsr1tmp} can be checked directly.
    \begin{enumerate}
        \item If condition \eqref{eq:cond_slsr1tmp} is satisfied for this pair, then $S_k = [S_k \; s_{k,i}]$ add $Y_k = [Y_k \; y_{k,i}]$
        \item Else, discard the pair $(s_{k,i},y_{k,i})$
    \end{enumerate}
\end{enumerate}

Using the mechanism described above, we recursively check condition \eqref{eq:cond_slsr1tmp}, and construct well defined pairs $S_k$ and $Y_k$ which are used for the calculation of $p_k$ and $\rho_k$. As mentioned above, we can use the same idea for implementing the LSR1 method. Finally, we should note that the idea described above was successfully implemented in a distributed variant on the S-LSR1 method in \cite{jahani2019scaling}.

\clearpage 

\subsection{Cost of communication}
\label{sec:A:Communication}

In this section, we show experiments conducted on a HPC cluster using a Cray Aries High Speed Network. The bandwidth ranges depending on the distance between nodes. We compiled the C++ code with the provided cray compiler.

\begin{figure}[h!]
\centering
\includegraphics[scale=0.3]{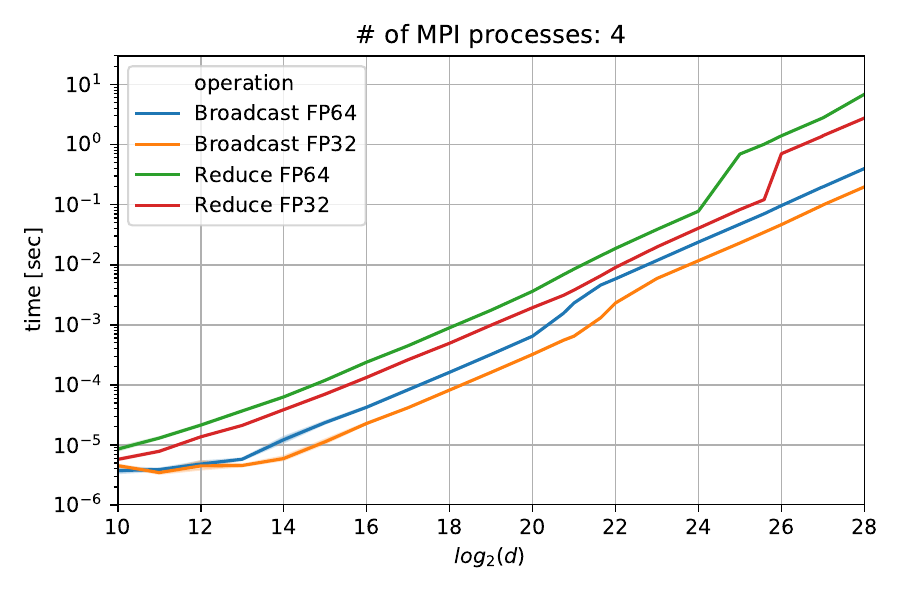}
\includegraphics[scale=0.3]{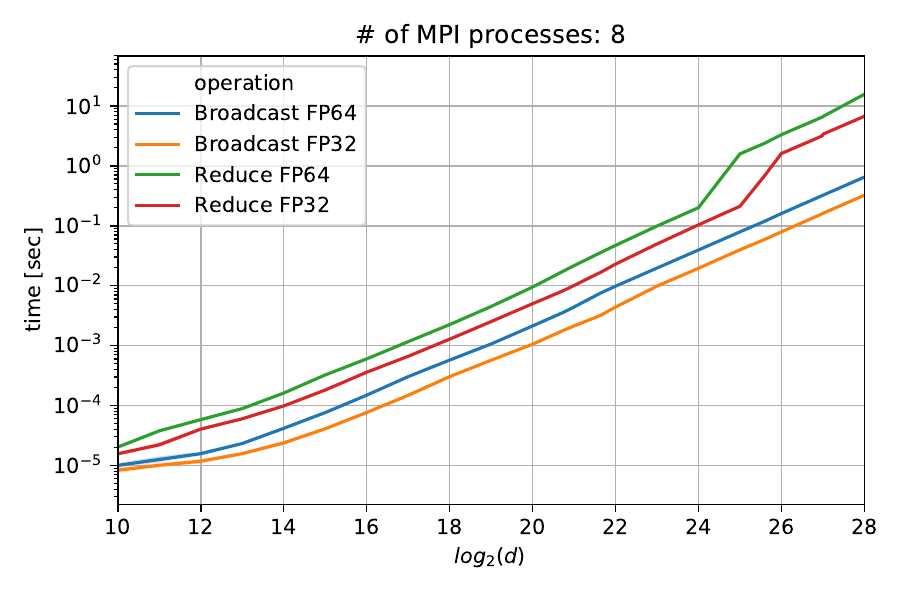}
\includegraphics[scale=0.3]{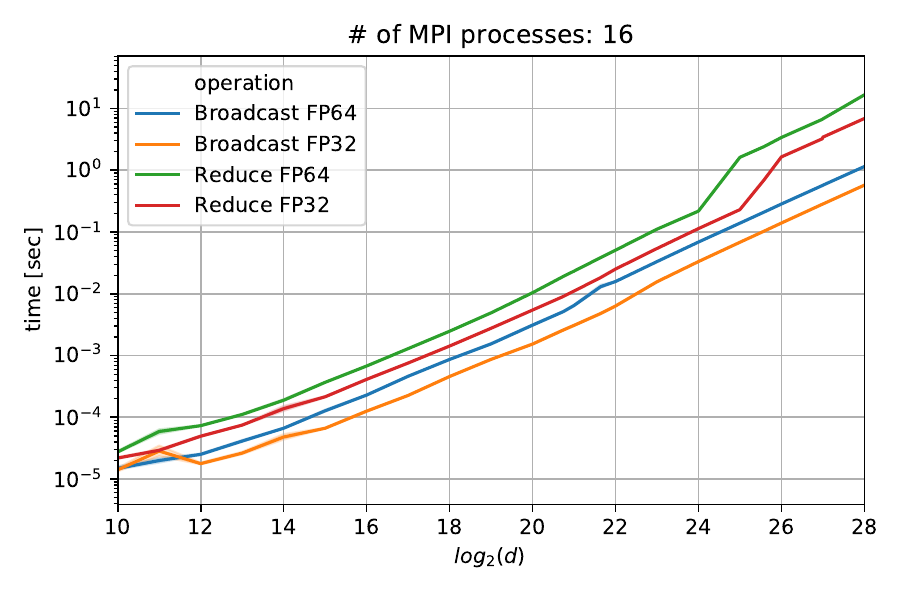}

\includegraphics[scale=0.3]{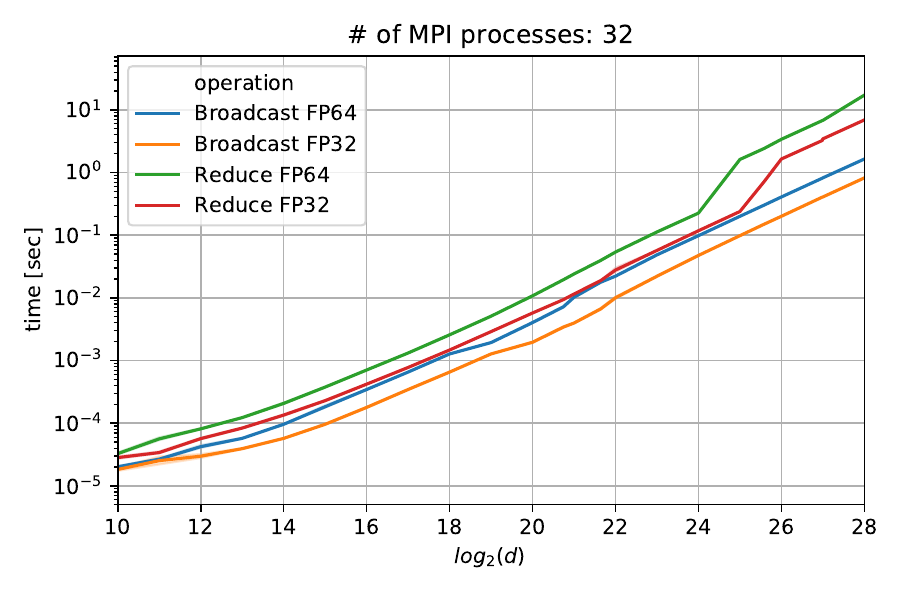}
\includegraphics[scale=0.3]{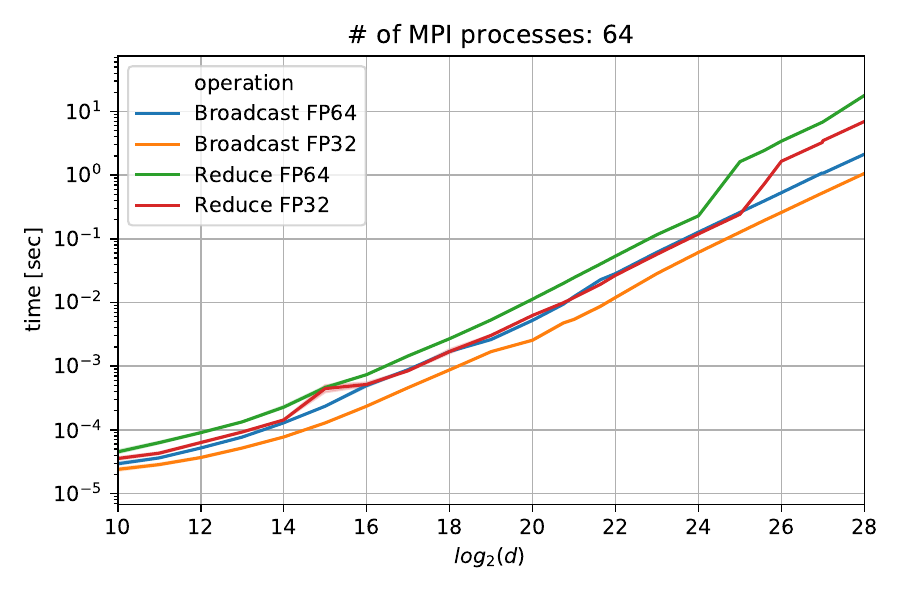}
\includegraphics[scale=0.3]{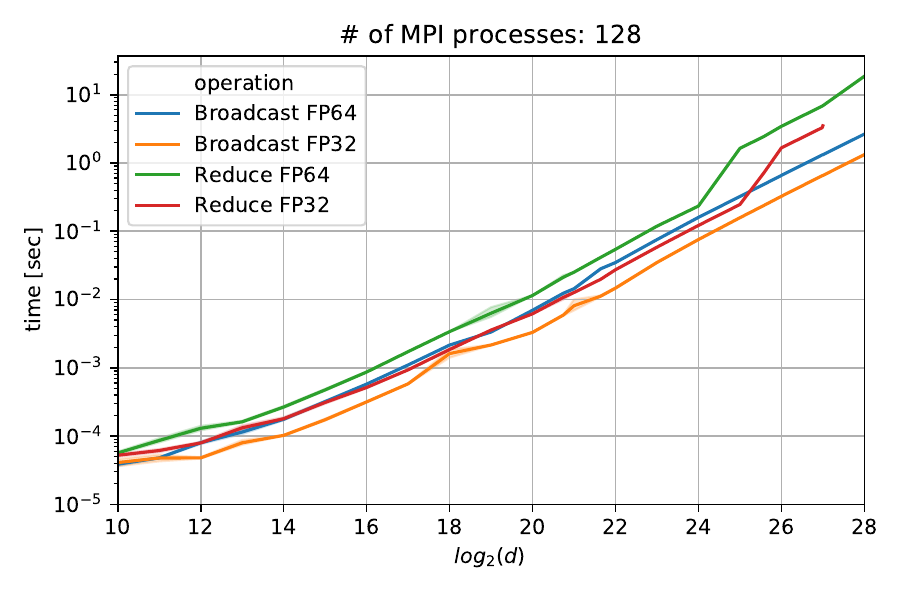}

\caption{Duration of Broadcast and Reduce for various number of MPI processes and different length of the vector.}
\label{fig:A:communicationCost}

\end{figure}

In Figure~\ref{fig:A:communicationCost}, we show how the duration (seconds) of Broadcast and Reduce increases when vectors of longer length are processed.

\begin{figure}[h!]

\centering
 
\includegraphics[scale=0.3]{figures_appendix/scaling_15.pdf}
\includegraphics[scale=0.3]{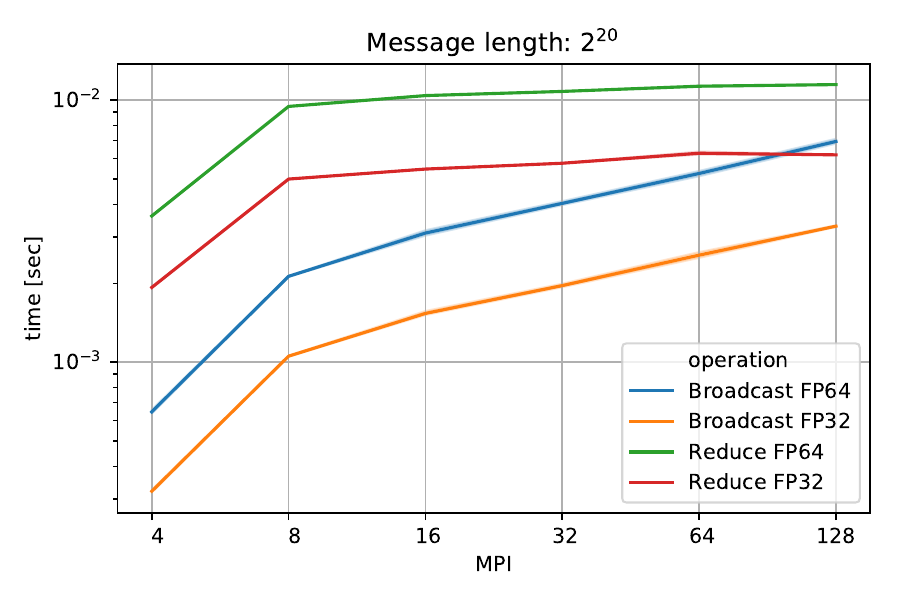}
\includegraphics[scale=0.3]{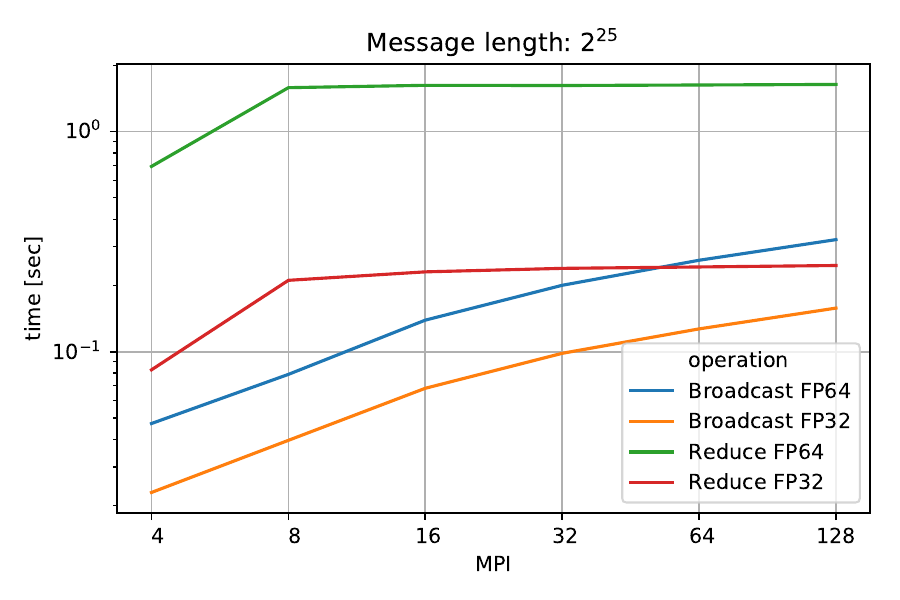}

\caption{Duration of Broadcast and Reduce as a function of \# of MPI processes for various lengts of vectors.}
\label{fig:A:scaling}

\end{figure}

In Figure~\ref{fig:A:scaling},
we show how long it takes (seconds) to perform Broadcast and Reduce operations for vectors of a given length if performed on different numbers of MPI processes. We have performed each operation 100 times and are showing the average time and 95\% confidence intervals.

\clearpage
\subsection{Toy Example}
\label{sec:app_toy}

In this section, we present additional numerical results for the toy classification problem described in Section 7.2. In the following experiments, we ran each method from 100 different initial points.

\subsubsection{Performance of Methods on small, medium and large toy classification problems - Box-plots}

The following box-plots show the accuracy achieved by different methods for different budgets (epochs) and iterations.

\begin{figure}[H]
	\centering
	\includegraphics[width=0.47\textwidth]{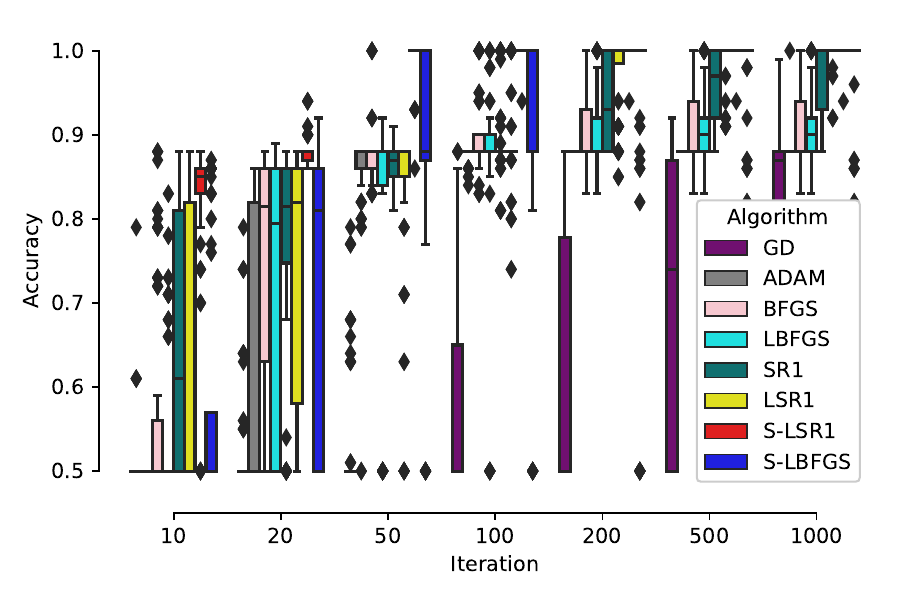}
	\includegraphics[width=0.47\textwidth]{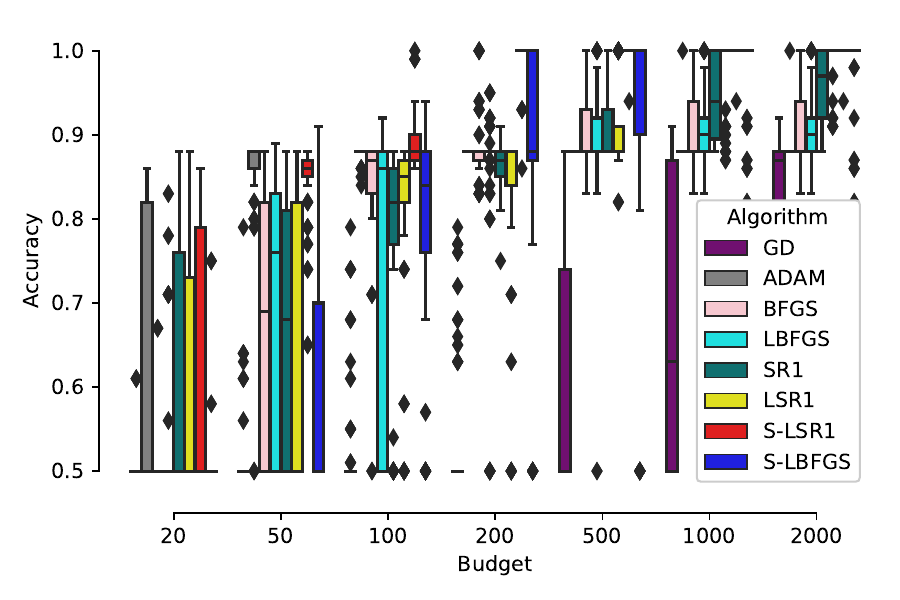}
	
	\includegraphics[width=0.47\textwidth]{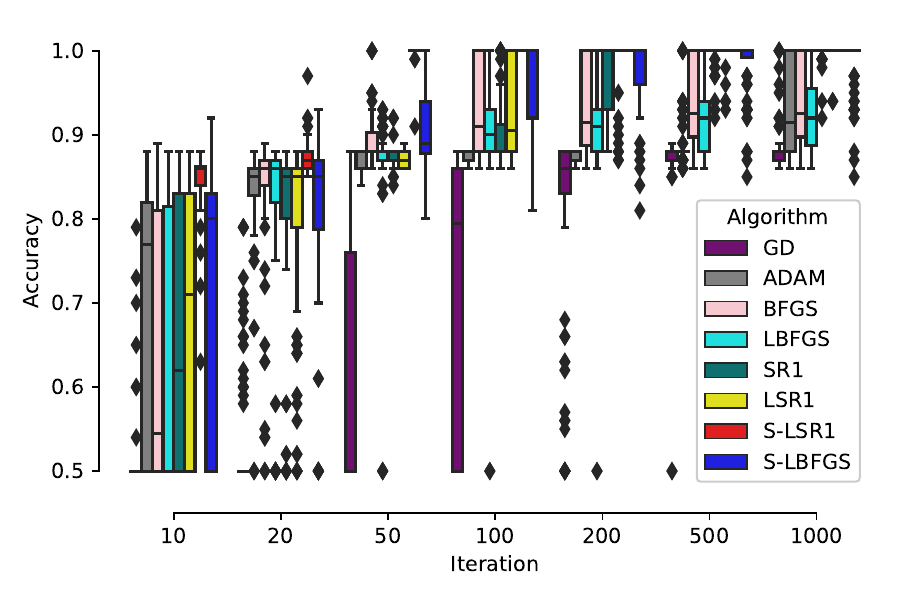}
	\includegraphics[width=0.47\textwidth]{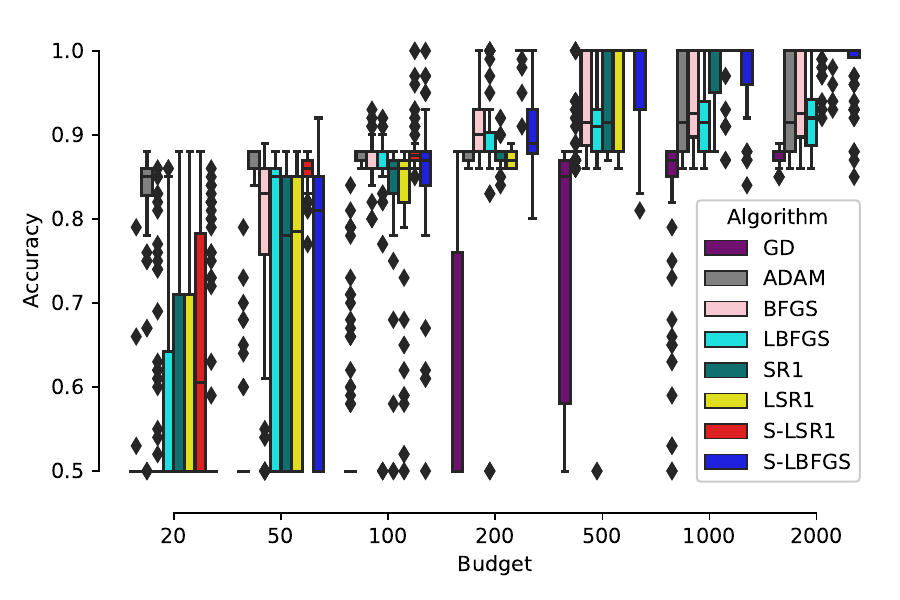}
	
	\includegraphics[width=0.47\textwidth]{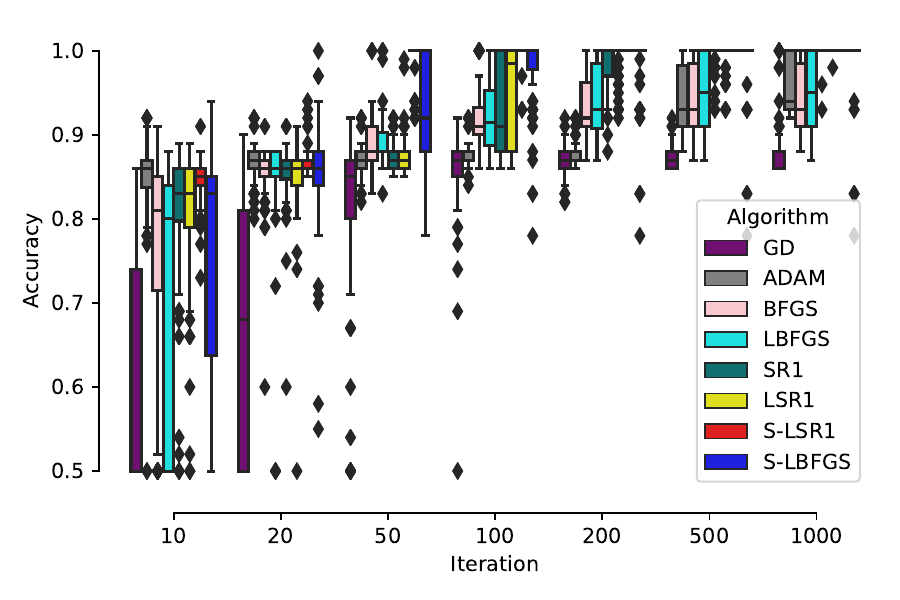}
	\includegraphics[width=0.47\textwidth]{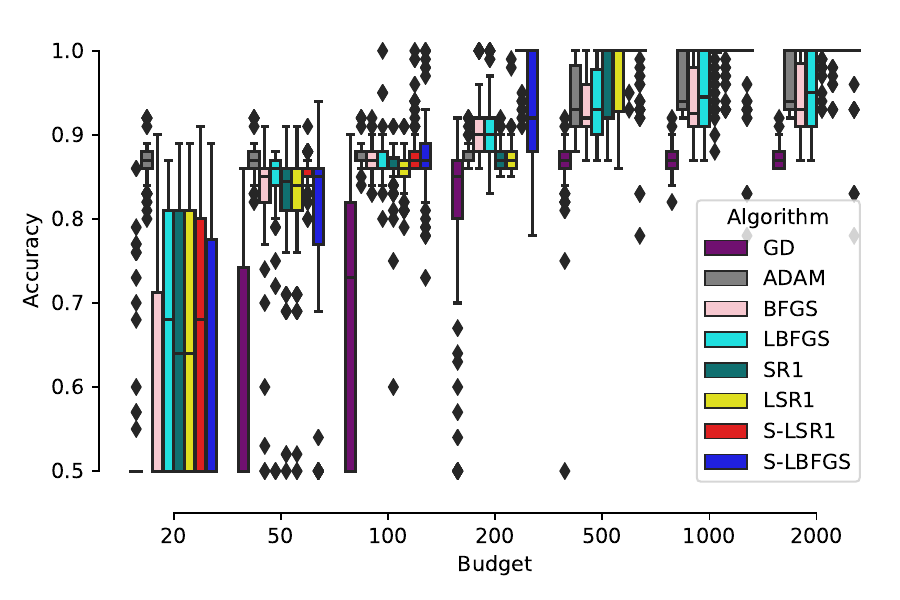}
	\caption{Performance of GD, ADAM, BFGS, LBFGS, SR1, LSR1, S-LSR1 and S-LBFGS on toy classification problem over \texttt{small network} (first row), \texttt{medium network} (second row), and \texttt{large network} (third row).}
	\label{SmalllNet32Altogether}
\end{figure}

\clearpage

\subsubsection{Performance of Methods on small, medium and large toy classification problems}
In this section, we present more experiments showing accuracy vs. iterations and accuracy vs. epochs for different methods on toy classification problem.
\begin{figure}[H]
	\centering
	\includegraphics[width=0.47\textwidth]{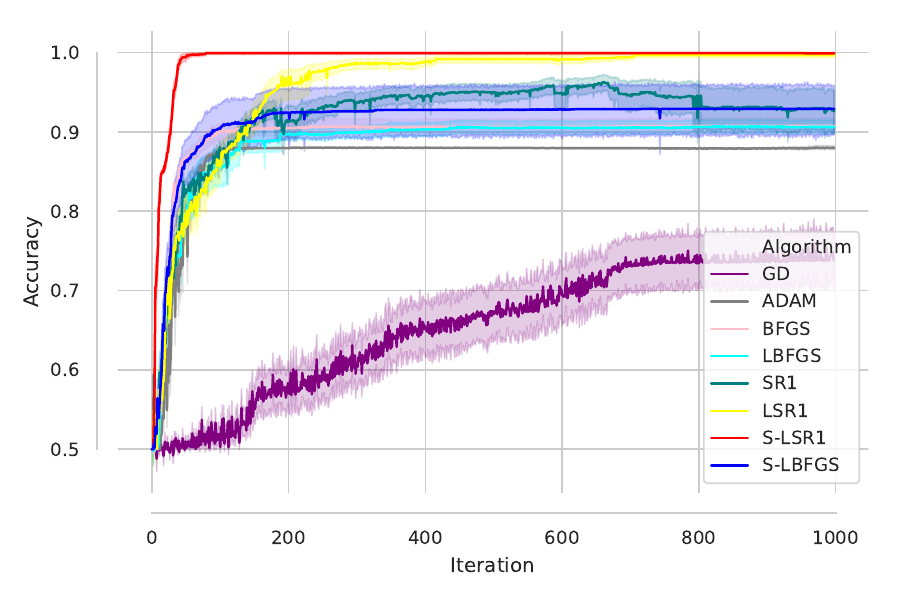}
	\includegraphics[width=0.47\textwidth]{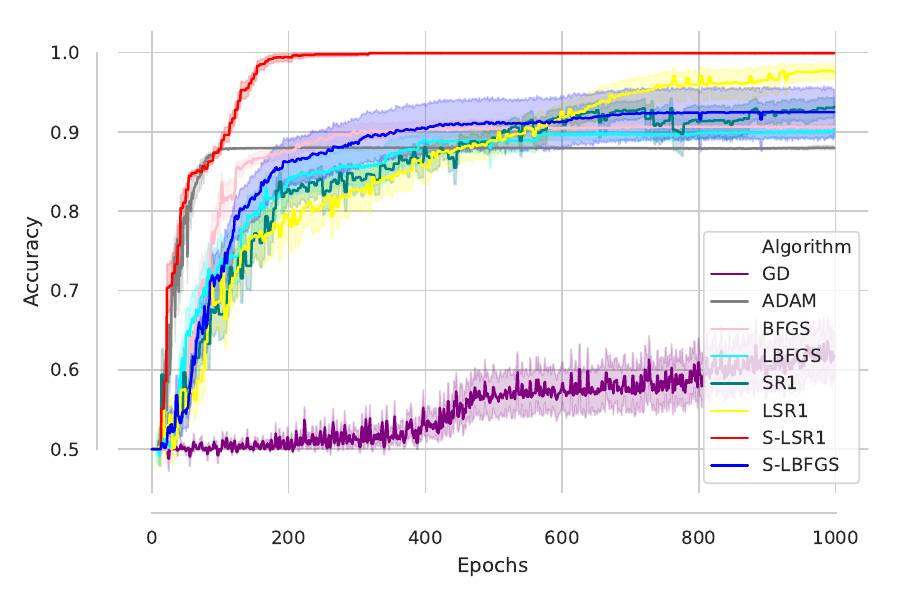}
	
	\includegraphics[width=0.47\textwidth]{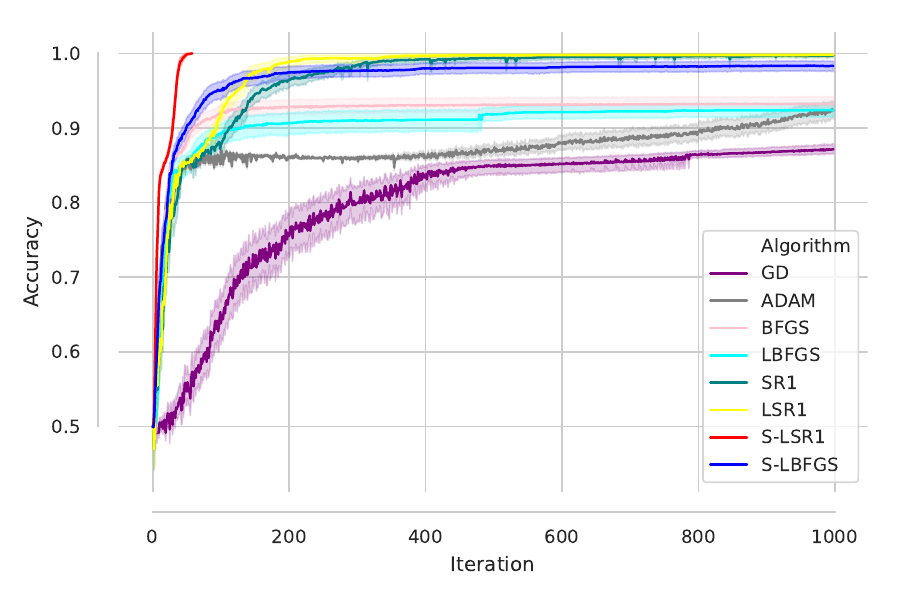}
	\includegraphics[width=0.47\textwidth]{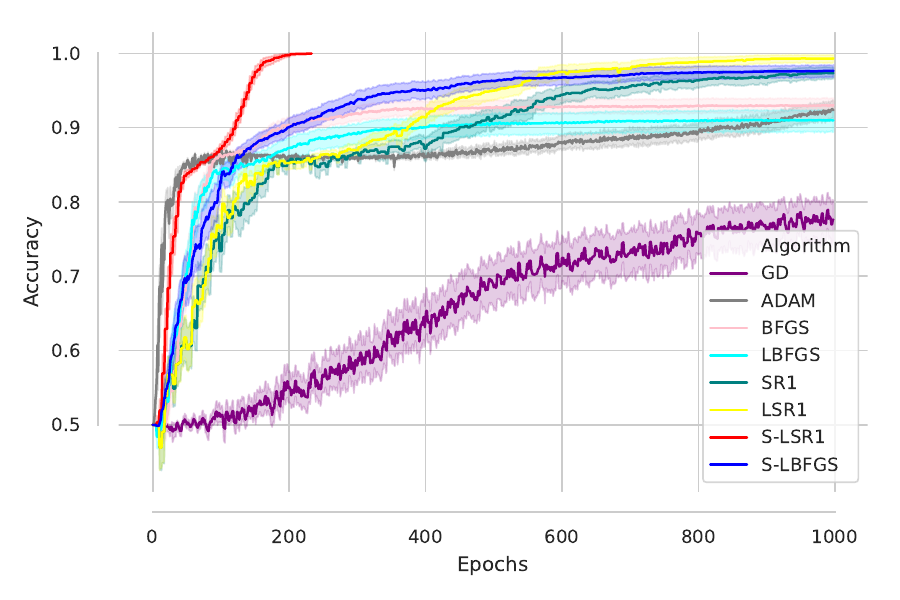}
	
	\includegraphics[width=0.47\textwidth]{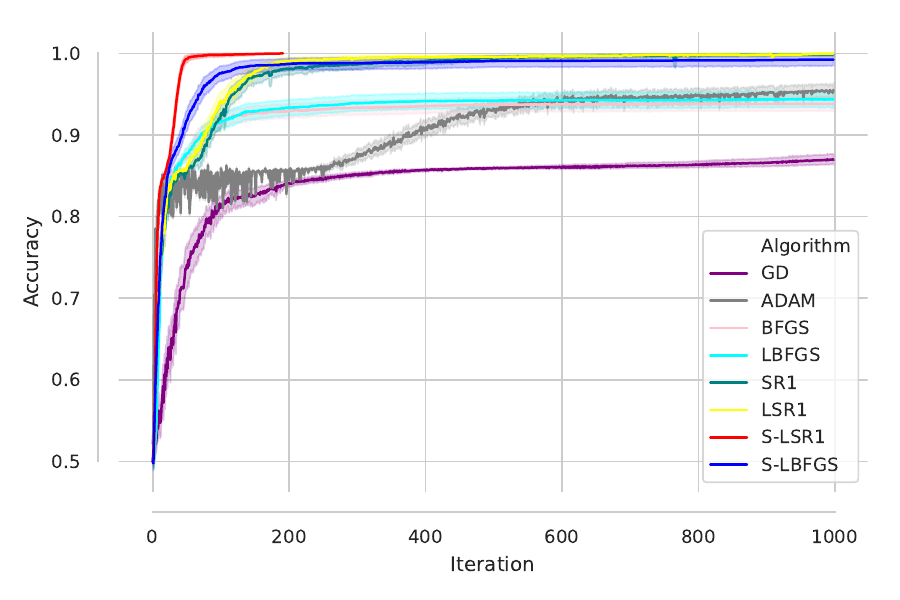}
	\includegraphics[width=0.47\textwidth]{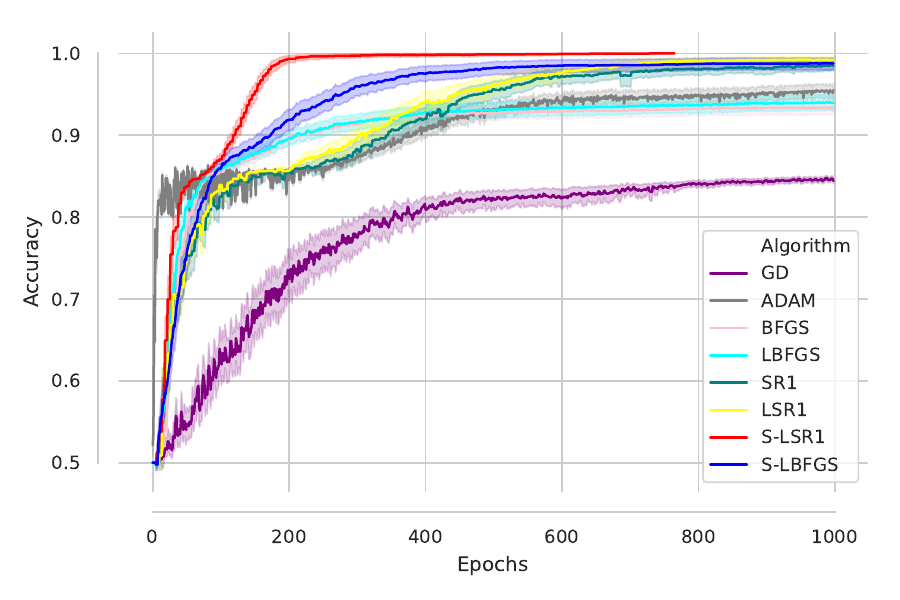}
	\caption{Performance of GD, ADAM, BFGS, LBFGS, SR1, LSR1, S-LSR1 and S-LBFGS on toy classification problem over \texttt{small network} (first row), \texttt{medium network} (second row), and \texttt{large network} (third row).}
	\label{SmalllNet422Altogether}
\end{figure}

\clearpage

\subsubsection{Comparison of BFGS-type methods}
In this section, we present more experiments showing the accuracy achieved in terms of iterations and epochs for BFGS-type methods  on toy classification problem.
\begin{figure}[H]
	\centering
	\includegraphics[width=0.47\textwidth]{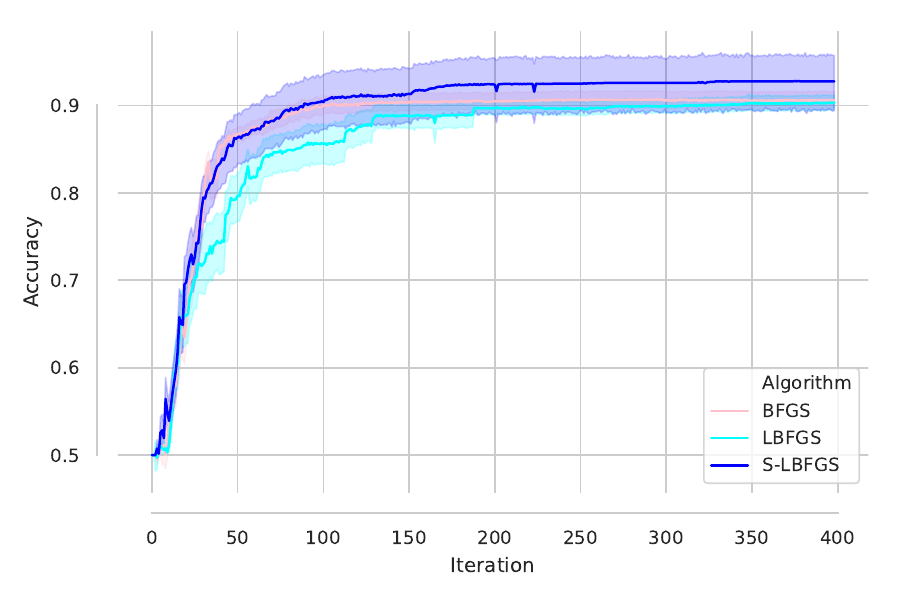}
	\includegraphics[width=0.47\textwidth]{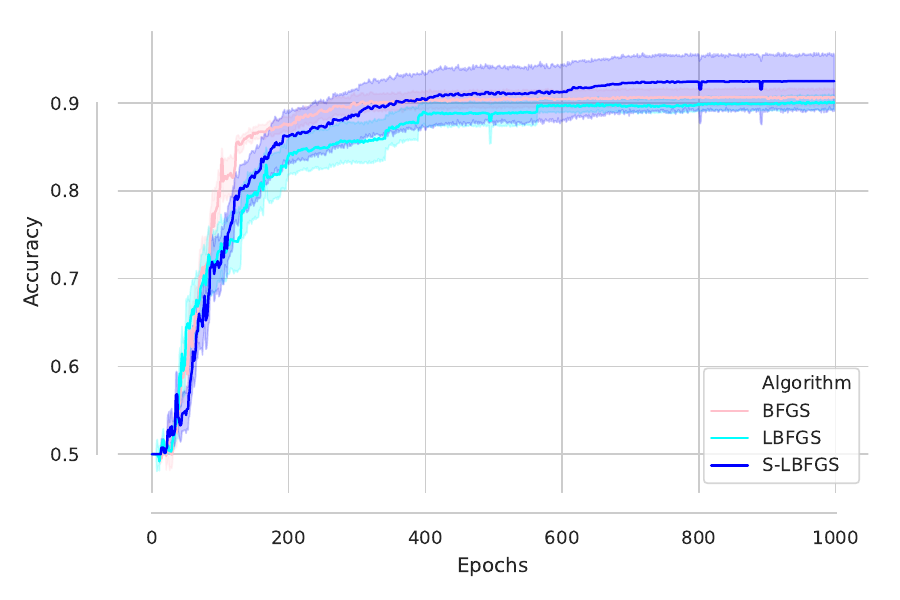}
	
    \includegraphics[width=0.47\textwidth]{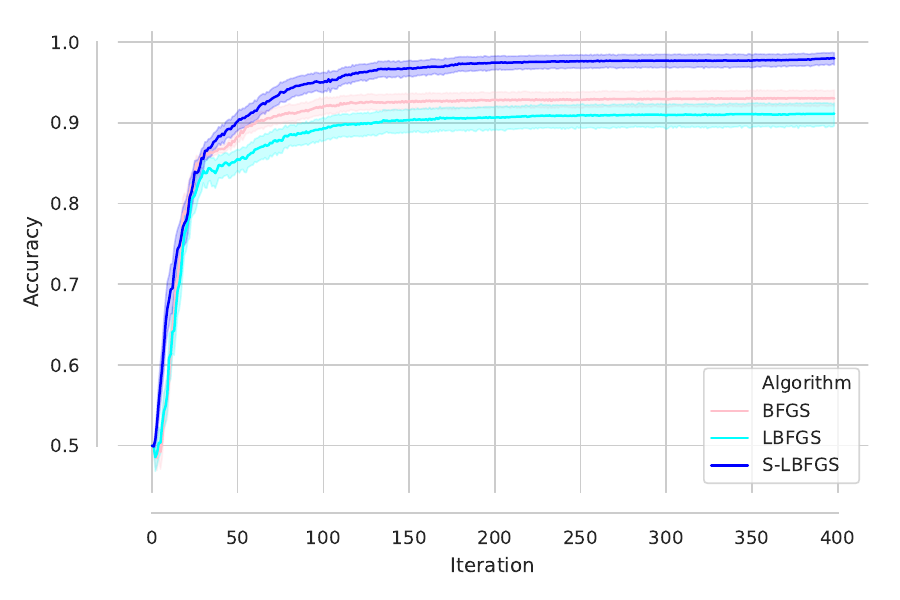}
	\includegraphics[width=0.47\textwidth]{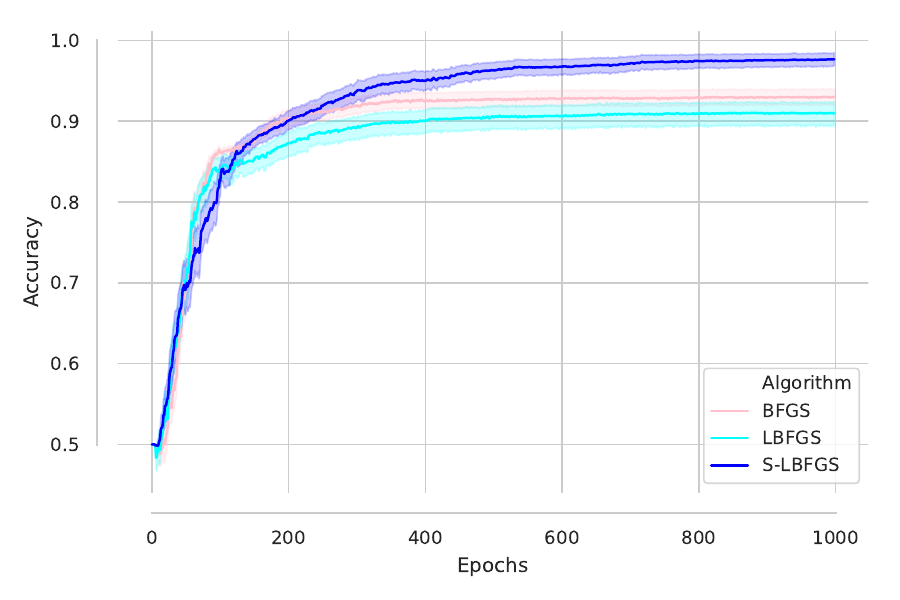}
	
	\includegraphics[width=0.47\textwidth]{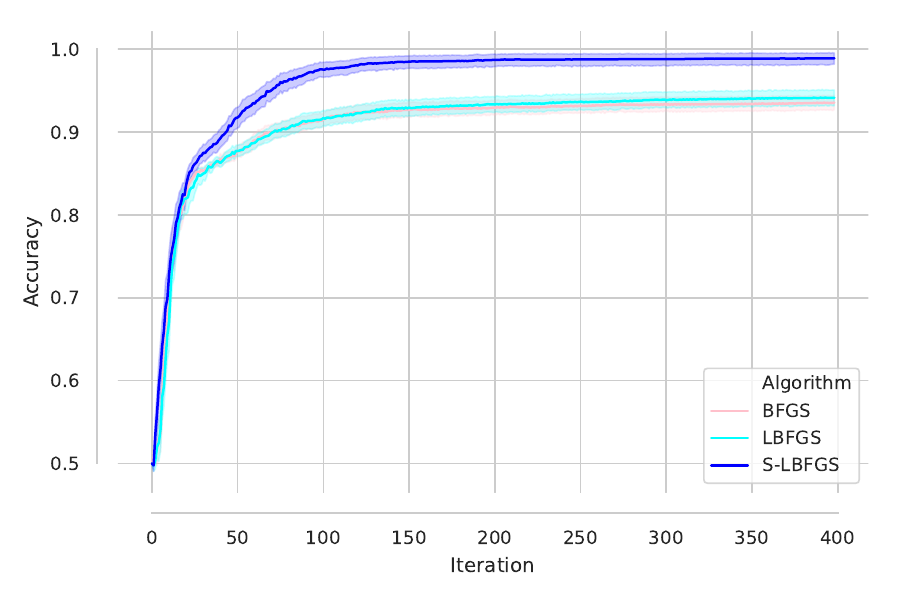}
	\includegraphics[width=0.47\textwidth]{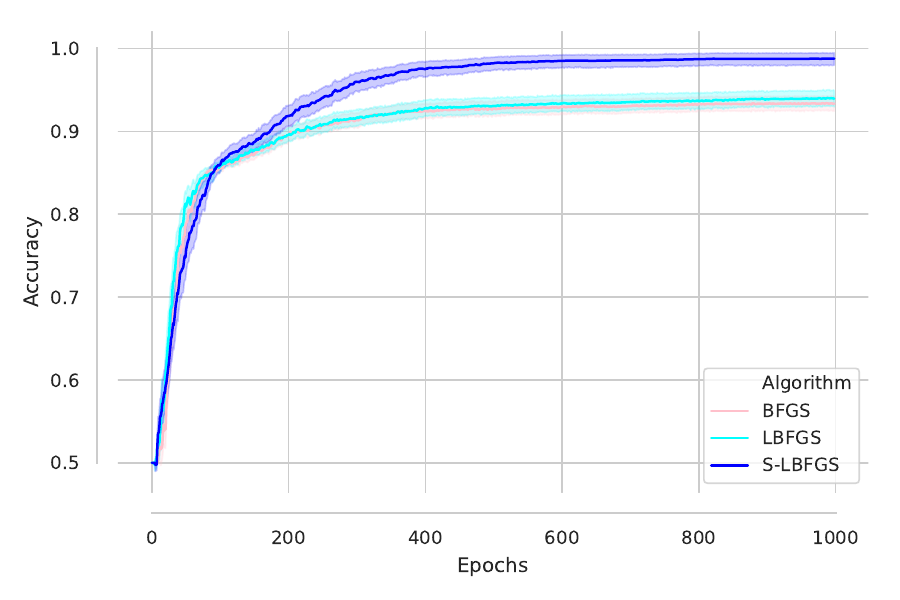}
	\caption{Performance of BFGS-type methods on toy classification problem over \texttt{small network} (first row), \texttt{medium network} (second row), and \texttt{large network} (third row).}
	\label{BFGSfamilyAltogether}
\end{figure}

\clearpage
\subsubsection{Comparison of SR1-type methods}
In this section, we present more experiments showing the accuracy achieved in terms of iterations and epochs for SR1-type methods  on toy classification problem.
\begin{figure}[H]
	\centering
	\includegraphics[width=0.47\textwidth]{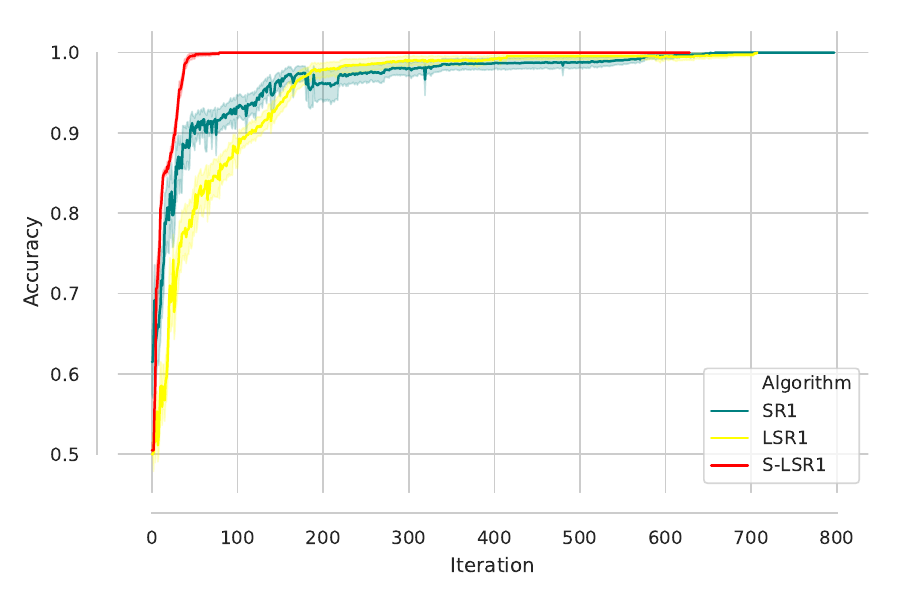}
	\includegraphics[width=0.47\textwidth]{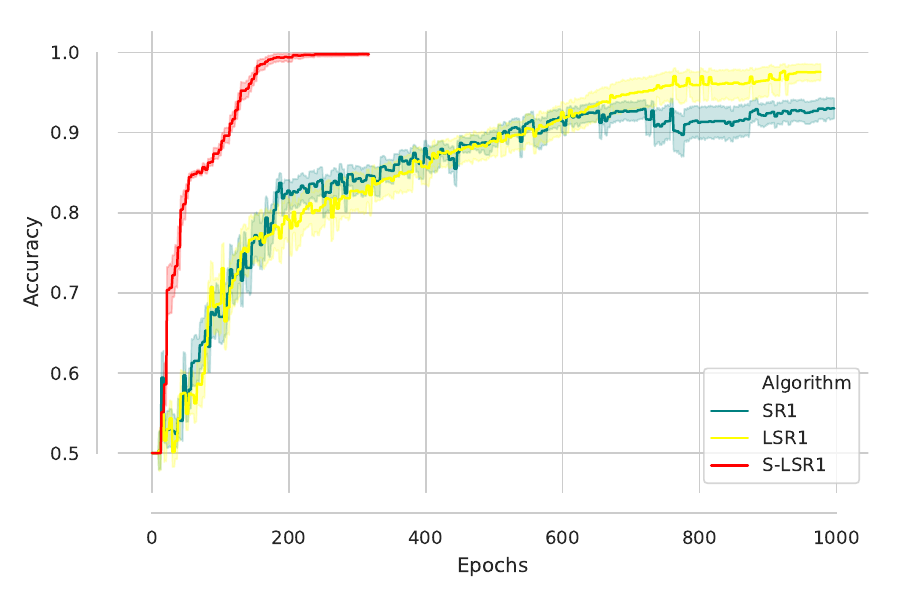}

	\includegraphics[width=0.47\textwidth]{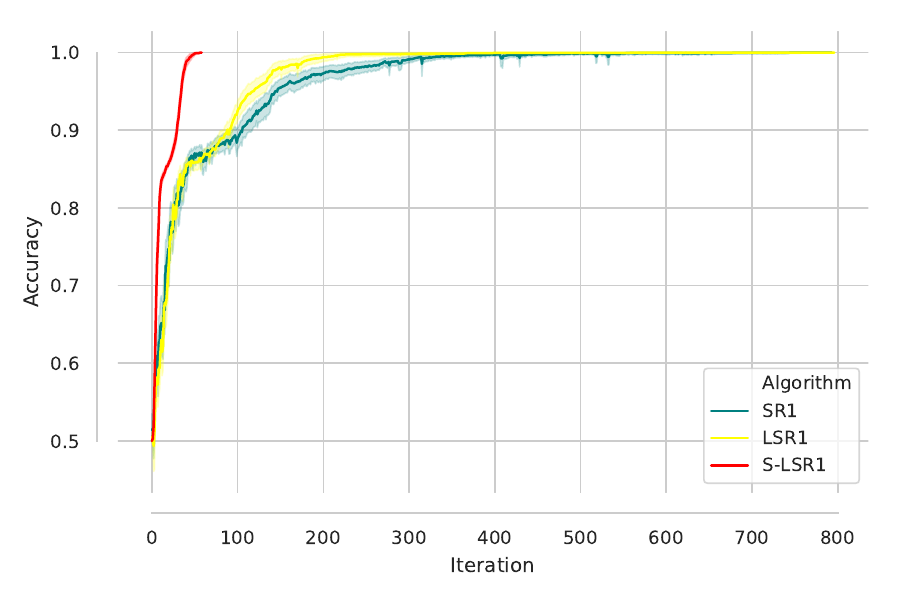}
	\includegraphics[width=0.47\textwidth]{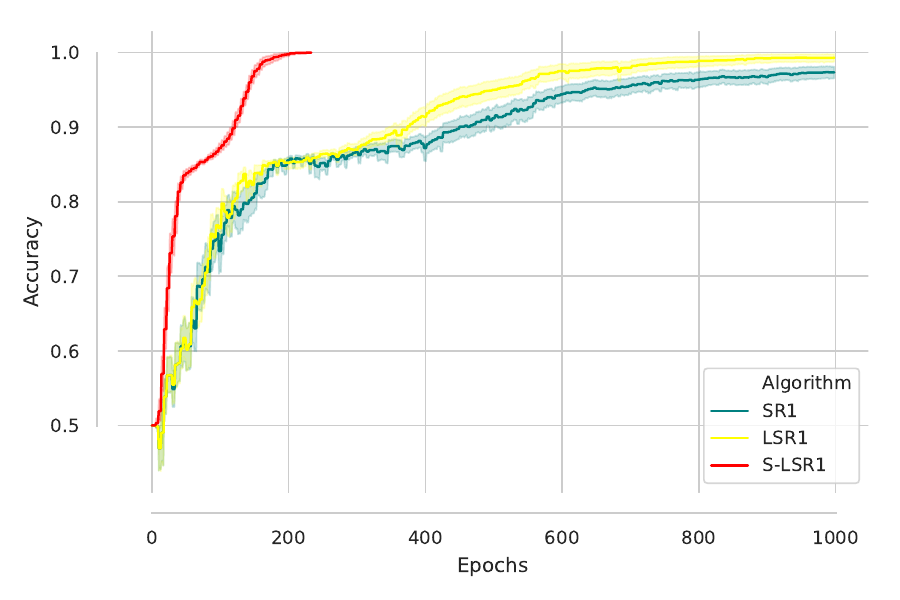}	

	\includegraphics[width=0.47\textwidth]{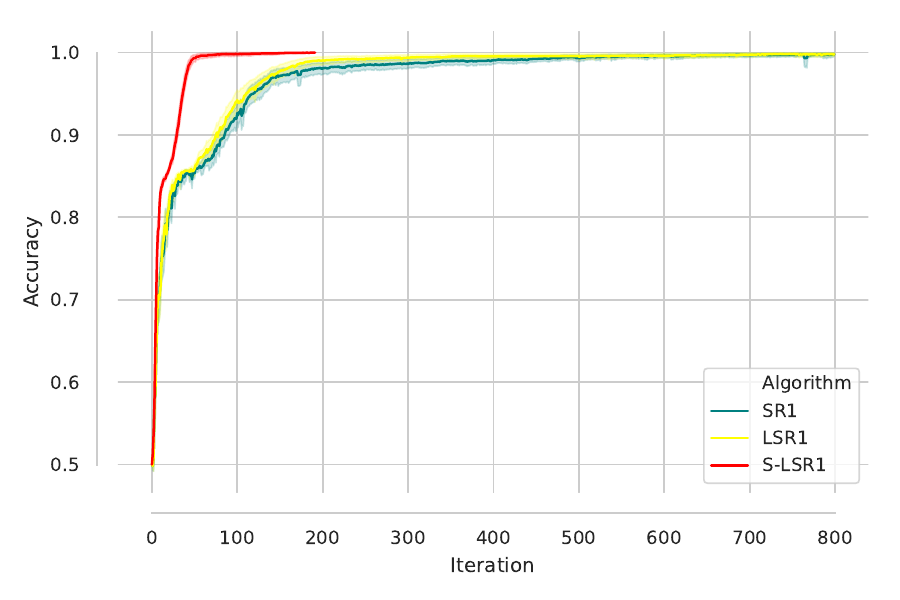}
	\includegraphics[width=0.47\textwidth]{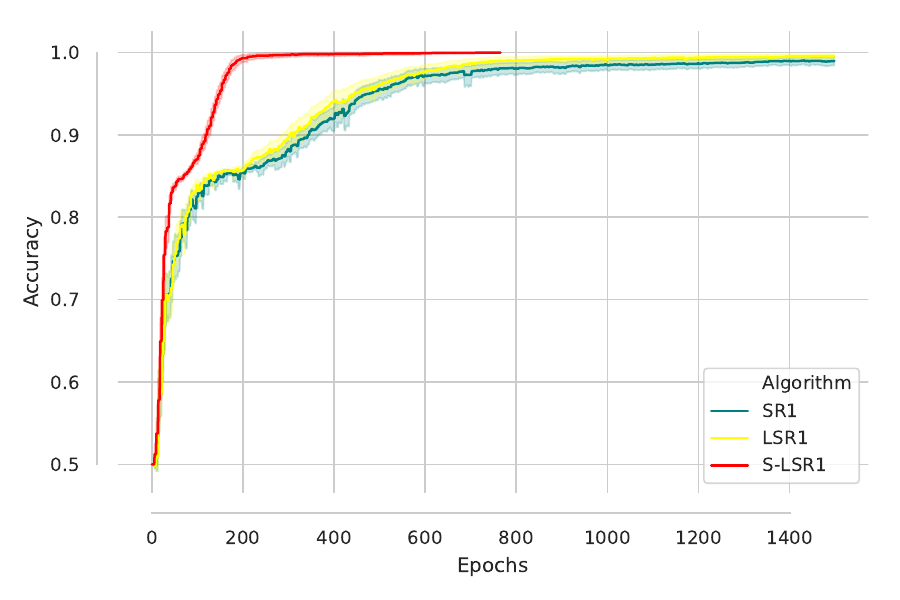}
	\caption{Performance of SR1-type methods on toy classification problem over \texttt{small network} (first row), \texttt{medium network} (second row), and \texttt{large network} (third row).}
	\label{SR1familyAltogether}
\end{figure}

\clearpage
\subsection{Binary Classification}\label{app.binary}

In this section, we present additional numerical results for the binary classification tasks Sections 7.3 \& 7.4. We show results on the following datasets\footnote{Datasets are available at \texttt{https://www.csie.ntu.edu.tw/~cjlin/libsvmtools/datasets/}}.

\begin{table}[htb]
	\centering
	\caption{Summary of datasets.}
	\begin{tabular}{lccc}
		\toprule
		\textbf{dataset} & \textbf{\# of samples } & \textbf{\# of features}& \textbf{\# of classes} \\
		\midrule
		\texttt{ijcnn1}& 35,000 & 22& 2\\\hdashline
		\texttt{rcv1}& 20,242& 47,326 & 2\\\hdashline
	\texttt{gisette}& 6,000& 5,000 & 2\\\hdashline
		\texttt{w8a}& 49,749& 300& 2\\\hdashline
		\bottomrule
	\end{tabular}
	
	\label{tab:datasets}
\end{table}
\subsubsection{Logistic Regression}\label{sec:logistic}

In this section, we present additional numerical results for the logistic regression problems. The figures in this section show: the optimality gap, norm of the gradient and training/testing accuracy.

\begin{figure}[H]
\centering
    \includegraphics[width=0.85\textwidth]{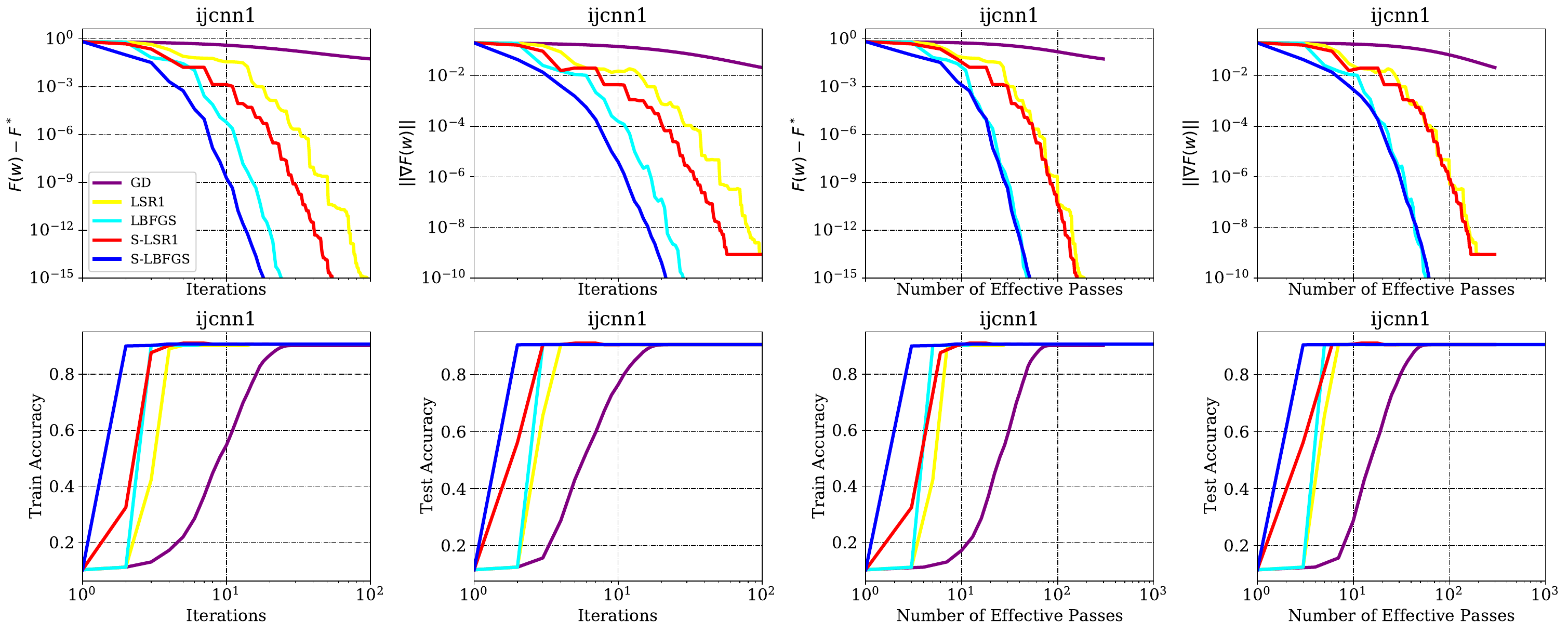}
    \caption{\small \texttt{ijcnn1}:
    Performance of GD, LBFGS, LSR1, S-LSR1 and S-LBFGS on Logistic Regression ($\lambda = 10^{-3}$).}
	\label{LR_appndx_ijcnn1_3}
\end{figure}

\begin{figure}[H]
\centering
    \includegraphics[width=0.85\textwidth]{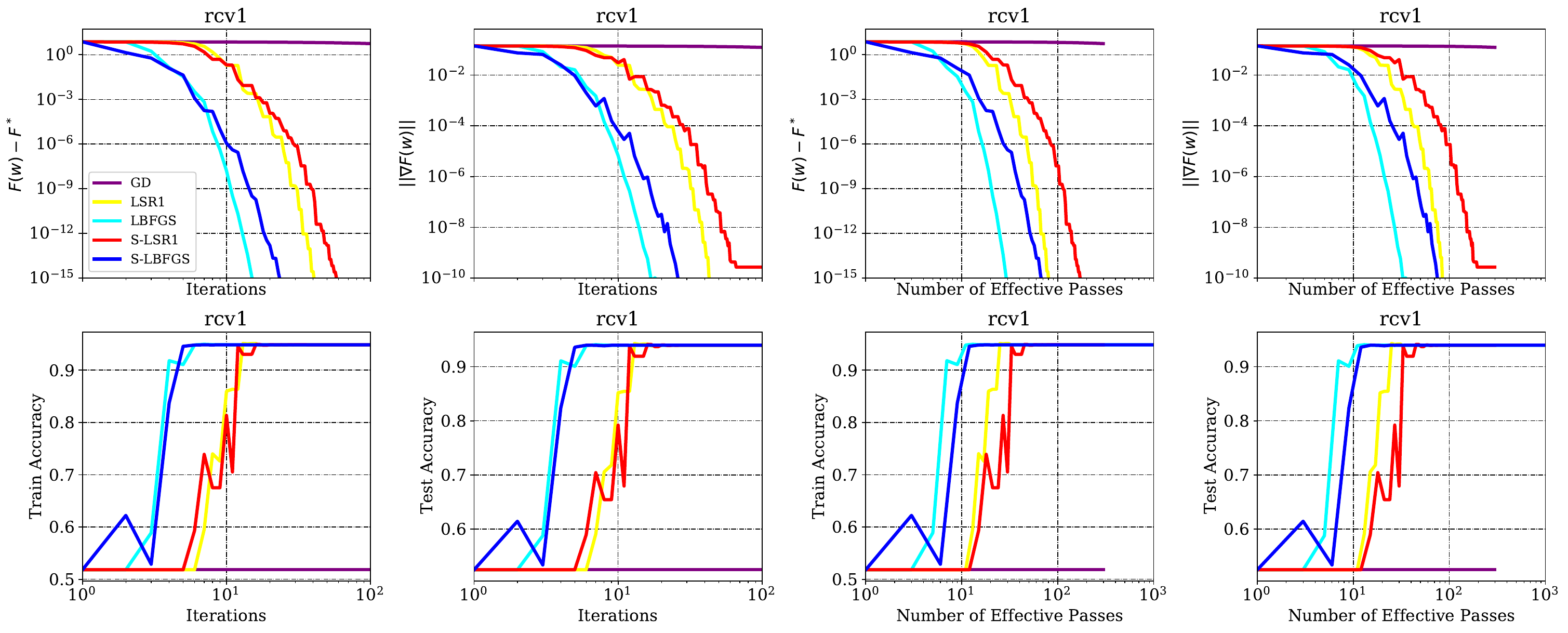}
    \caption{\small \texttt{rcv1}:
    Performance of GD, LBFGS, LSR1, S-LSR1 and S-LBFGS on Logistic Regression ($\lambda = 10^{-3}$).}
	\label{LR_appndx_rcv1_3}
\end{figure}

\begin{figure}[H]
\centering
    \includegraphics[width=0.85\textwidth]{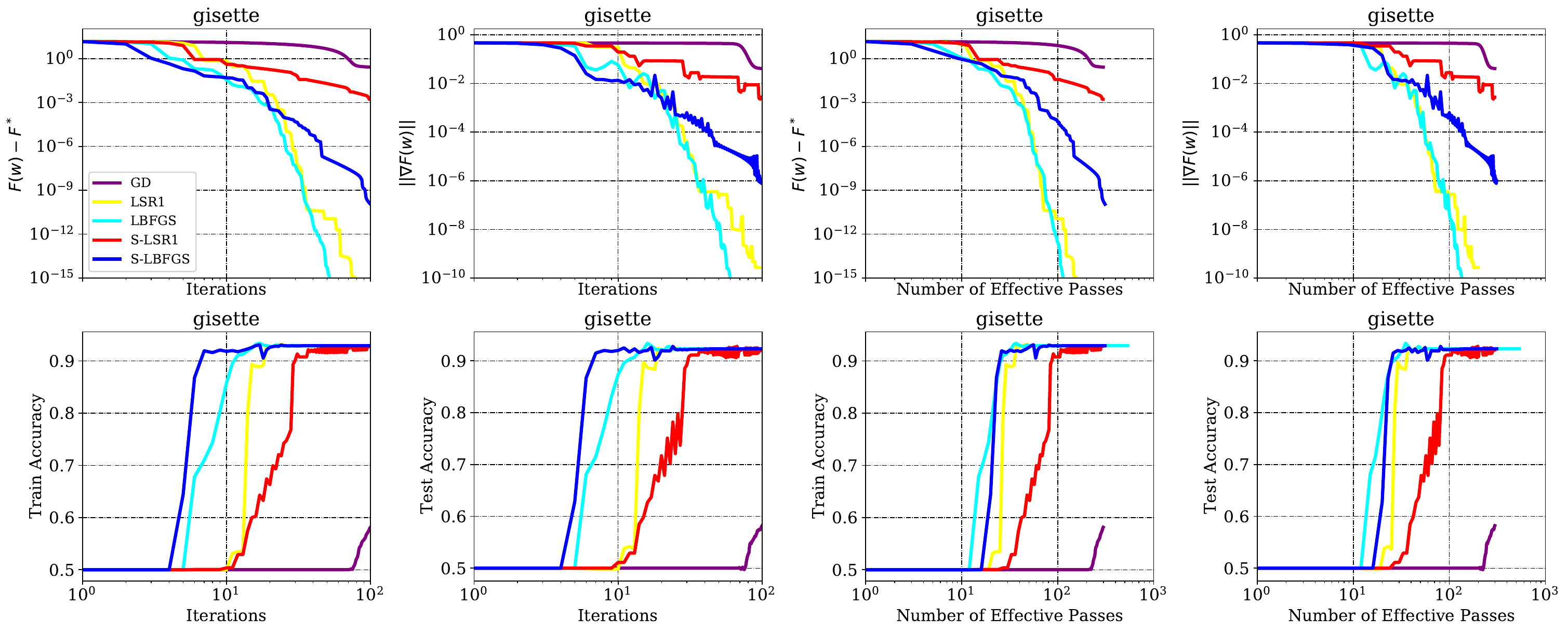}
    \caption{\small \texttt{gisette}:
    Performance of GD, LBFGS, LSR1, S-LSR1 and S-LBFGS on Logistic Regression ($\lambda = 10^{-3}$).}
	\label{LR_appndx_gisette_3}
\end{figure}

\begin{figure}[H]
\centering
    \includegraphics[width=0.85\textwidth]{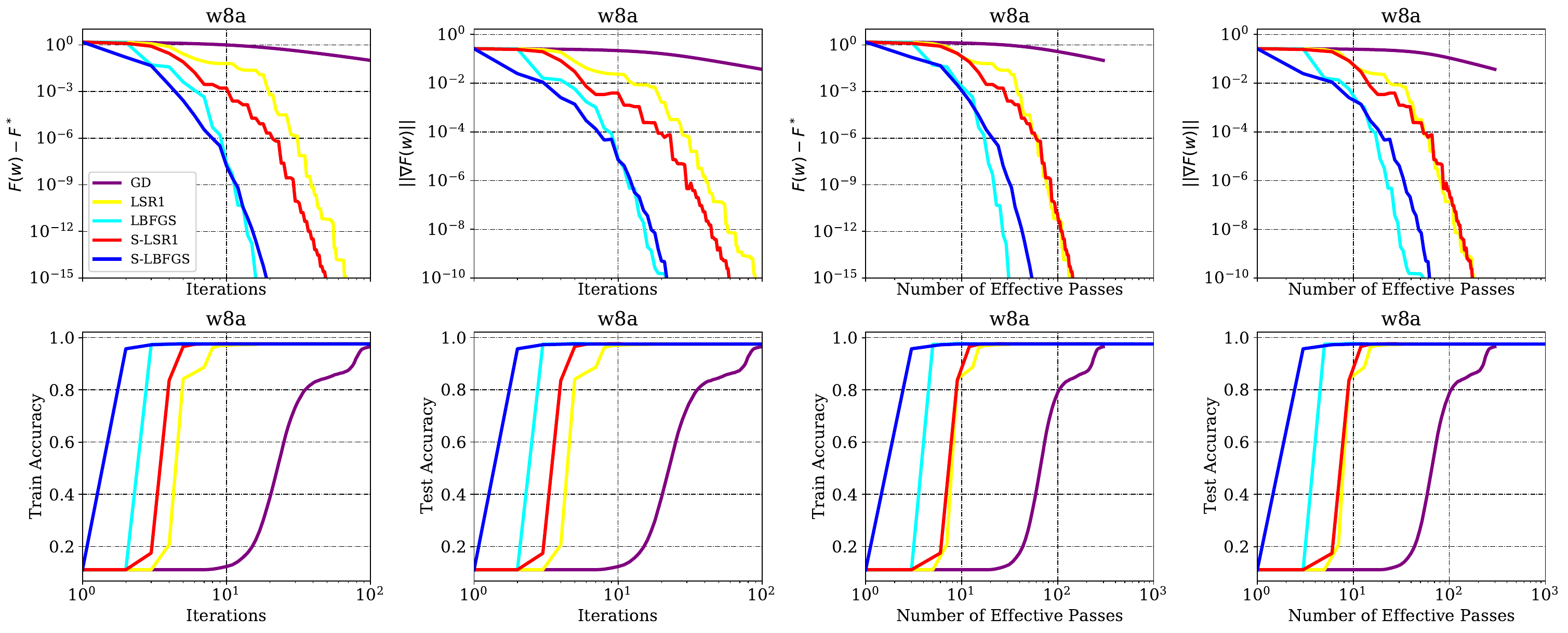}
    \caption{\small \texttt{w8a}:
    Performance of GD, LBFGS, LSR1, S-LSR1 and S-LBFGS on Logistic Regression ($\lambda = 10^{-3}$).}
	\label{LR_appndx_w8a_3}
\end{figure}

\clearpage

\subsubsection{Non-Linear Least Squares}\label{sec:nonLS}

In this section, we present additional numerical results for the nonlinear least squares problems.

\begin{figure}[H]
\centering
    \includegraphics[width=0.85\textwidth]{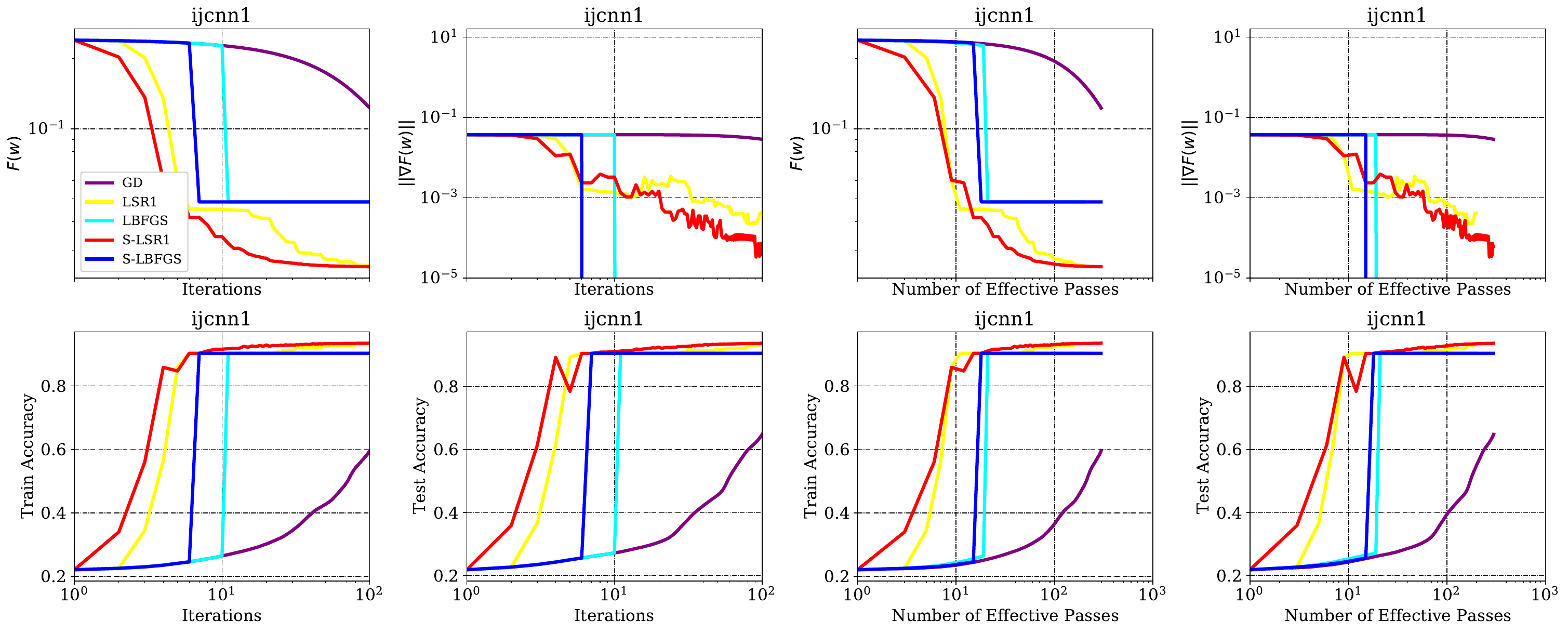}
    
    \caption{\small \texttt{ijcnn1}: Performance of GD, LBFGS, LSR1, S-LSR1 and S-LBFGS on Non-Linear Least Square.}
	\label{NLS_appndx_ijcnn1}
\end{figure}

\begin{figure}[H]
\centering
    \includegraphics[width=0.85\textwidth]{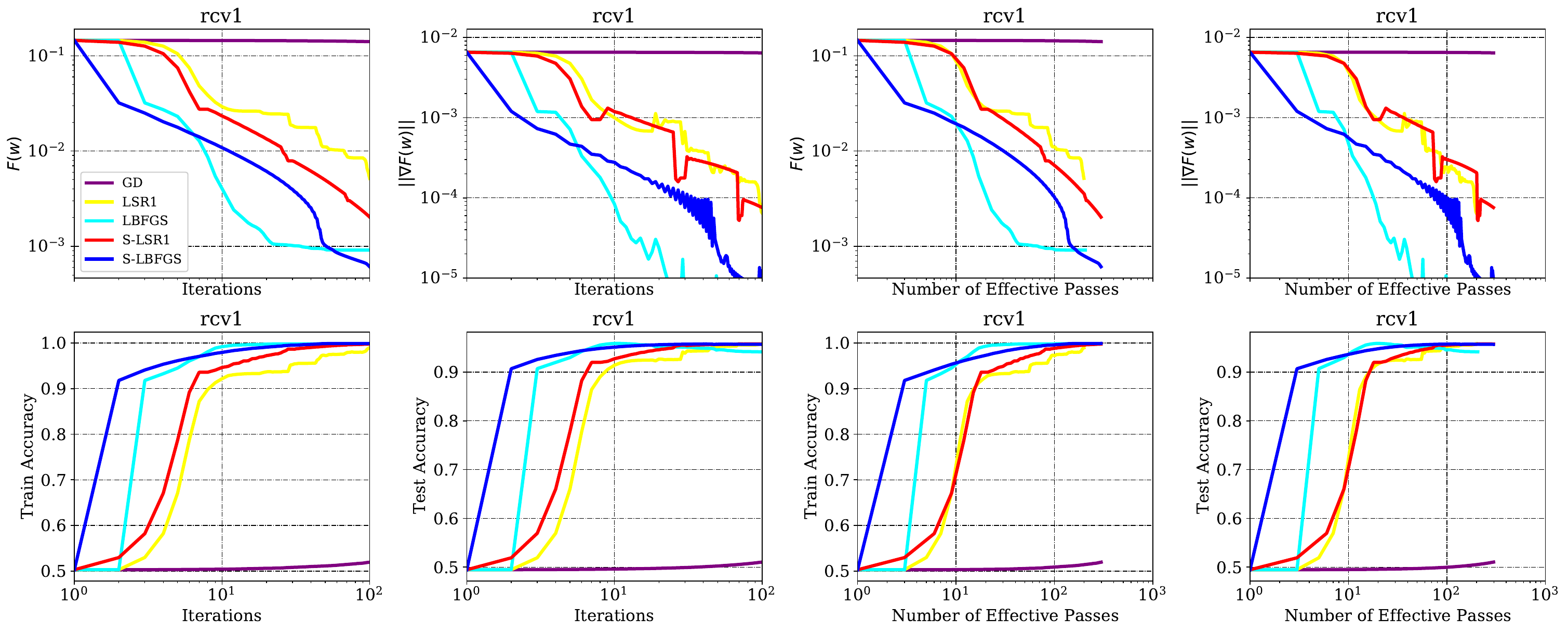}
    \caption{\small \texttt{rcv1}:
    Performance of GD, LBFGS, LSR1, S-LSR1 and S-LBFGS on Non-Linear Least Square.}
	\label{NLS_appndx_rcv1}
\end{figure}

\begin{figure}[H]
\centering
    \includegraphics[width=0.85\textwidth]{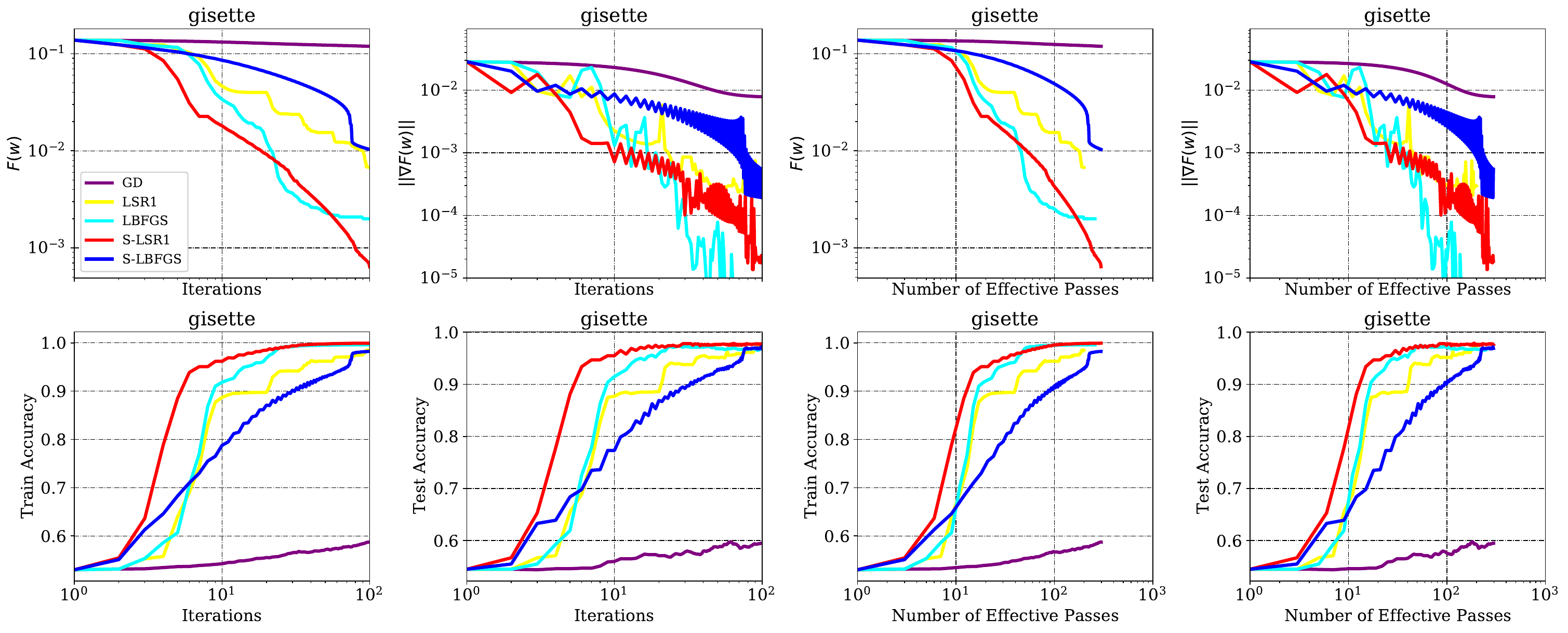}
    \caption{\small \texttt{gisette}: Performance of GD, LBFGS, LSR1, S-LSR1 and S-LBFGS on Non-Linear Least Square.}
	\label{NLS_appndx_gisette}
\end{figure}

\begin{figure}[H]
\centering
    \includegraphics[width=0.85\textwidth]{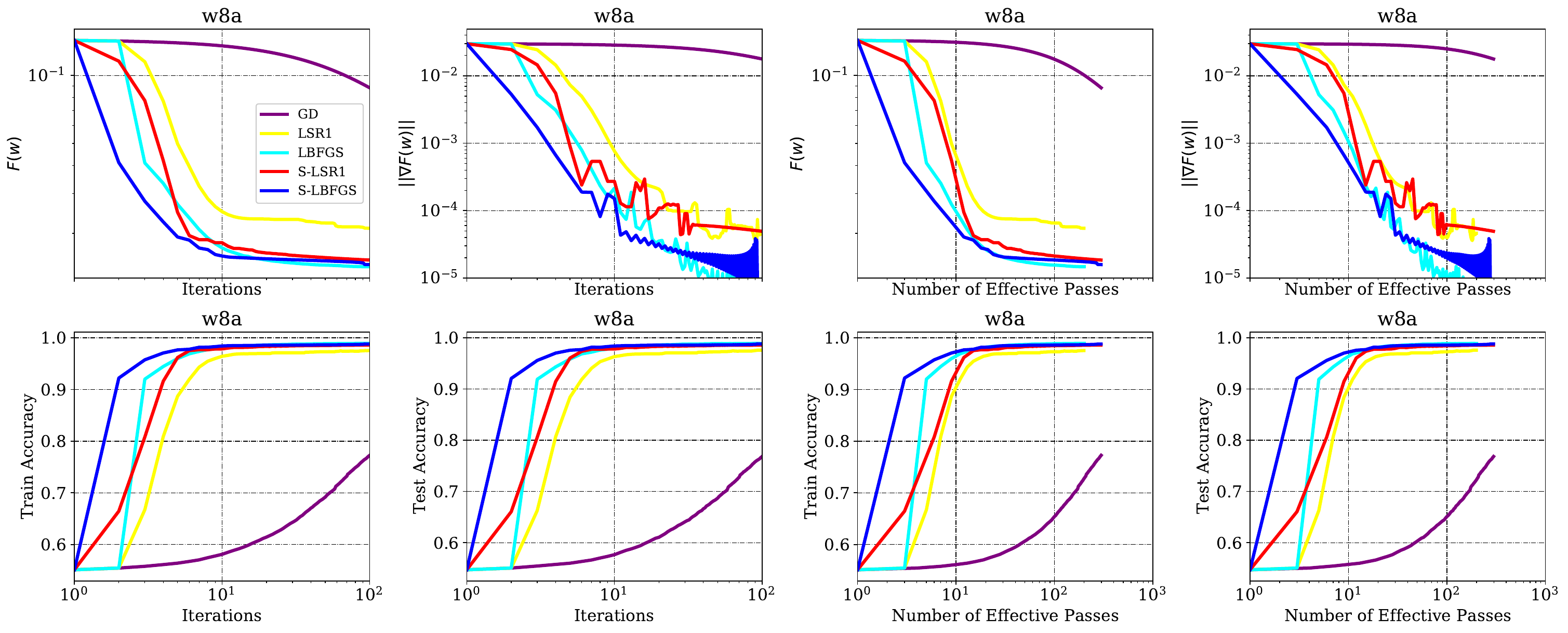}
    \caption{\small \texttt{w8a}: Performance of GD, LBFGS, LSR1, S-LSR1 and S-LBFGS on Non-Linear Least Square.}
	\label{NLS_appndx_w8a}
\end{figure}


\subsection{Performance of ADAM on MNIST}
\label{sec:app_MNIST}
\begin{figure}[ht]
	\centering
 	
 	\includegraphics[width=0.32\textwidth]{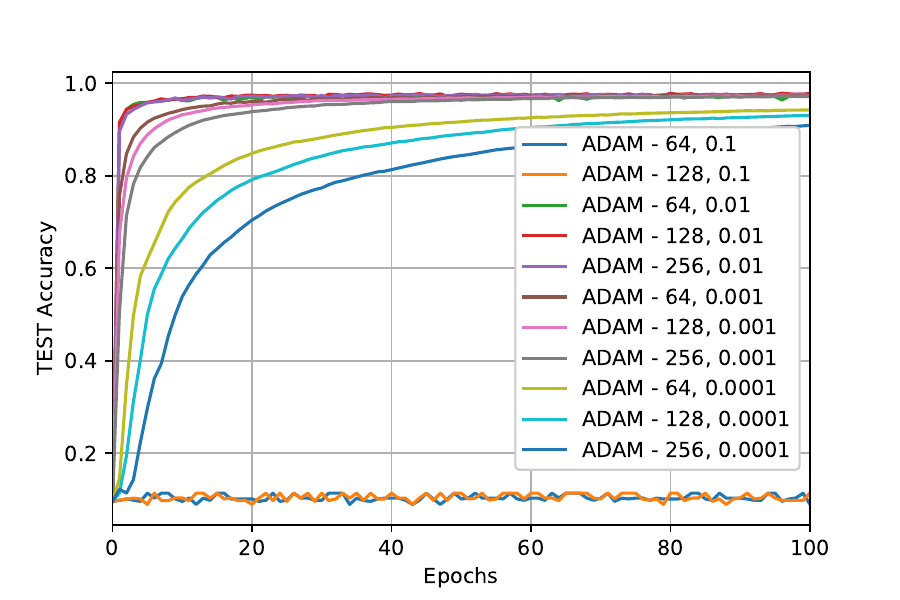}
 	\includegraphics[width=0.32\textwidth]{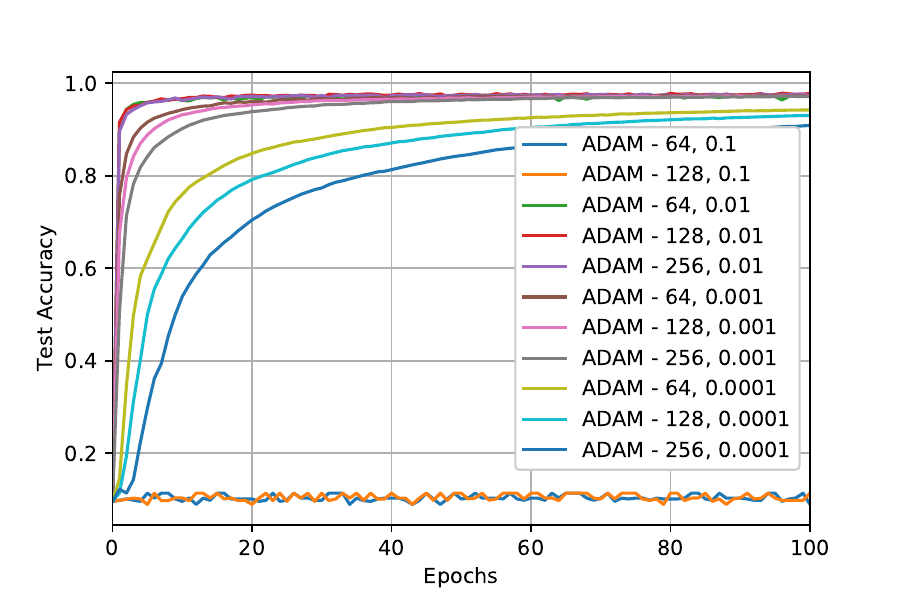}
 	\includegraphics[width=0.32\textwidth]{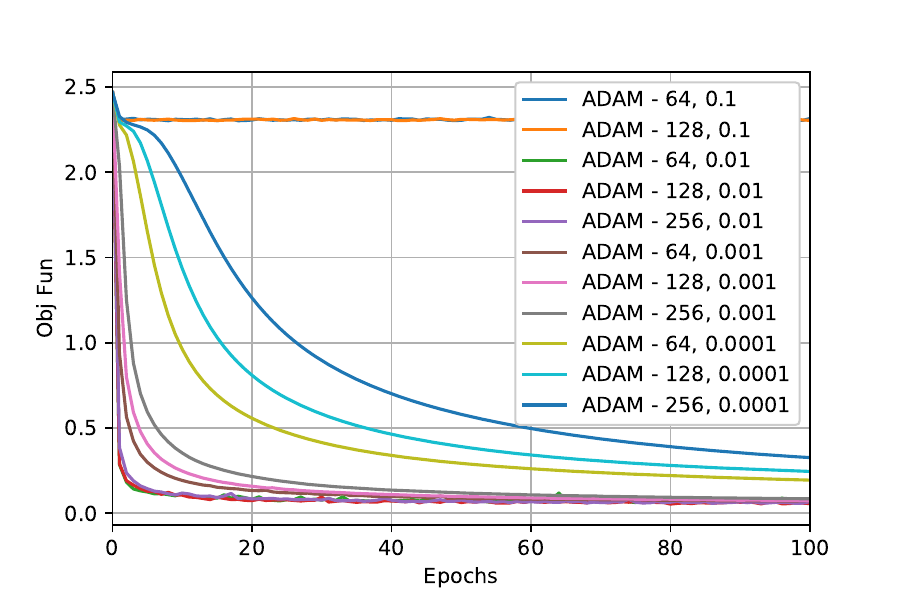}	
 	
	\caption{Performance of ADAM with different steplengths on MNIST over \texttt{Net1} architecture.}
	\label{accuarcyMNIST_ADAM}
\end{figure}

In this section, we show the performance of ADAM with different steplenghts on the MNIST problem. As is clear from the results in Figure \ref{accuarcyMNIST_ADAM}, the performance of well-tuned ADAM is very good, however, when the steplength is not chosen correctly, the performance of ADAM can be terrible. Note, we have omitted runs for which ADAM diverged (i.e., when the steplength was chosen to be too large).

\bibliographystyle{tfs}
\bibliography{refs_sampledQN}

\end{document}